\documentstyle[twoside]{article}
\input{math.mac}
\renewcommand{\beware}{}

\newtheorem{Thm}{Theorem}[subsection]
\newtheorem{Prop}[Thm]{Proposition}
\newtheorem{Cor}[Thm]{Corollary}
\newtheorem{Main}{}

\newtheorem{Lem}[Thm]{Lemma}

\newtheorem{Def}[Thm]{Definition}
\newtheorem{rmk}[Thm]{Remark}
\newenvironment{Rmk}{\begin{rmk}\em}{\end{rmk}}
\newtheorem{exm}[Thm]{Example}
\newenvironment{Exm}{\begin{exm}\em}{\end{exm}}

\newcommand{\qed}{\par {\EM QED} }
\newtheorem{prf}{Proof}

\newenvironment{Prf}{\begin{prf}\em}{\qed\end{prf}}
\newtheorem{prff}{}

\newenvironment{Prff}{\begin{prff}\em}{\qed\end{prff}}

\newcommand{\spn}{{\mbox{\it satur}\,}}
\newcommand{\supp}{{\mbox{supp}\,}}
\newcommand{\dist}{{\mbox{distance}}}
\newcommand{\diam}{{\mbox{diameter}}}
\newcommand{\Pre}{{\mbox{Pre}\,}}
\newcommand{\cnt}{{\mbox{count}\,}}
\newcommand{\In}{{\mbox{\it in}\,}}
\newcommand{\Ex}{{\mbox{\it ex}\,}}
\newcommand{\Arr}{{\mbox{\Arr}\,}}

\newcommand{\fn}{function}  
\newcommand{\fns}{functions}
\newcommand{\cts}{continuous}
\newcommand{\nbd}{neighbourhood}

\title{Kac's Formula, Vertices of Chains 
and Equidecomposability of Functions and
Enhanced Functions}
 
\author{Eliahu Levy\\
Department of Mathematics,
Technion IIT, 32000 Haifa, Isreal}

\begin{document}

\pagestyle{empty}
\setlength{\oddsidemargin}{+0.1in}
\setlength{\evensidemargin}{-0.3in}

\newcommand{\mytitle}{
       KAC'S FORMULA, VERTICES OF CHAINS,\\
       INFINITESIMAL MEASURES,\\ 
       AND EQUIDECOMPOSABILITY OF FUNCTIONS\\
       AND ENHANCED FUNCTIONS\\}

\NOT{
\title{Generalizations of Kac's Formula,\\
       Vertices of Chains, Enhanced Functions,\\
       Equidecomposability\\
       - a memorandum}

\title{Generalizations of Kac's Formula,\\
       Vertices of Chains, Infinitesimal Measures,\\
       Equidecomposability
       \NOT{\\- a memorandum}}
}
\begin{center}
{\LARGE\bf\mytitle}
\end{center}
\begin{flushright}
\vspace{18cm}
{\LARGE\bf Eliahu Levy~~~~~~~\\}
\end{flushright}

\cleardoublepage

\begin{center}
{\LARGE\mytitle}
\vspace{3cm}
{\Large RESEARCH THESIS\\
\vspace{1.5cm}
        SUBMITTED IN PARTIAL FULFILLMENT OF THE REQUIREMENTS\\
        FOR THE DEGREE OF 
        DOCTOR OF PHILOSOPHY\\
\vspace{3cm}
        Eliahu Levy\\
\vspace{6cm}
        SUBMITTED TO THE SENATE OF\\
        THE TECHNION -- ISRAEL INSITUTE OF TECHNOLOGY\\
        TISHREI 5760~~~~~~~~~~HAIFA~~~~~~~~~~OCTOBER 1999\\}
\end{center}
        
\cleardoublepage
{
\setlength{\parindent}{0cm}
THE RESEARCH THESIS WAS DONE UNDER THE SUPERVISION OF PROF.\
VITALY BERGELSON IN THE FACULTY OF MATHEMATICS.\\
\par\vspace{6cm}
I THANK PROF.\ VITALY BERGELSON FOR HIS SUPERVISION, ENCOURAGEMENT AND HELP.\\
\par\bigskip
I ALSO THANK PROF.\ YOAV BENYAMINI FOR HIS ENCOURAGEMENT AND HELP.\\
\par\bigskip
THE GENEROUS FINANCIAL HELP OF THE TECHNION IS GRATEFULLY ACKNOWLEDGED.\\
\par\vspace{6cm}
THIS WORK IS DEDICATED TO THE MEMORY OF MY DEAR FATHER,
WHO DID NOT LIVE TO SEE IT COMPLETED.\\
}



\cleardoublepage

\tableofcontents

\cleardoublepage
\pagestyle{headings}
\setcounter{page}{1}
\section*{Synopsis} \label{S:SYN}
\addcontentsline{toc}{section}{Synopsis}

The present work hinges on the five parts of the title: Kac's formula,
equidecomposability, vertices of chains, infinitesimal measures,
enhanced functions.

The notion of chains and their vertices is proposed as a language
to define and treat equidecomposability of \fns\ on measure spaces acted
measure-preservingly by groups. Thus, equidecomposable \fns\ are vertices
of the same chain\NOT{(for details see the beginning of Section
\ref{S:DS})}. The fact that vertices of the same chain
have the same integral is proposed as a way to get many formulas,
including and generalizing Kac's formula\NOT{\cite{Kac}},
which states that when the
acting group is $\bZ$ and the measure is probability, the relative
expectation of the return time to a measurable subset $E$ with conull
saturation is the reciprocal of the measure of $E$. This formula is
viewed as a consequence of the fact that
the return time (defined as zero outside $E$)
and the characteristic \fn\ of the ``set of points that visited $E$'',
differing from the saturation in a set of measure $0$, are equidecomposable,
a fact derived from an argument in the language of ``chains'' akin to
the well-known Wright's proof using the ``Kakutani skyscraper''\NOT{
(see \cite{Petersen}, pp.\ 45--46))}.

When $G$ is discrete (which we shall refer to as the ``discrete case''),
the equality in integral of equidecomposable \fns\ is trivial, so is the
passage to the language of vertices of chains, yet this language is an
easy path to find many formulas, developed in Section \ref{S:DS}.
These deal with $\bZ$ as acting group, where one gets modifications of
Kac's formula, as well as with multidimensional groups (speaking,
for example, on ``partially ordered Time''). One finds connections among
attributes of (repeated) return and arrival in the $\bZ$ case
(cf.\ results of Kastelyn\NOT{\cite{Kastelyn}}).
One has here a kind of ``toy'' where possibly interesting
formulas can emerge from the ``playing''.

It should be noted that ergodicity is usually irrelevant to the above
considerations and assuming it seems of little help. Thus in Kac's Formula
what matters is only that the saturation of $E$ is conull.

What gives more emphasis to this language is the fact that it extends
to the ``continuous'' case, to which Section \ref{S:CT} is devoted, and
which for our purpose is Borel measure-preserving action of a 2nd-countable
locally compact group $G$ on a standard $\sigma$-finite measure space
(one cannot go beyond locally compact groups, since Haar measure is crucial
for our treatment).
This allows us to formulate ``continuous'' Kac-like formulas in complete
analogy with the discrete case (and with the prospect of getting more
by ``playing''), {\em provided $G$ is unimodular}.
But in order to do this and circumvent the obvious
fact that simple-minded Kac formula fails for $\bR$-actions (i.e.\ flows)
-- for $\bR$ acting on the circle $\bT$ by rotation and $E$ composed of
intervals the return time to $E$ is $\ne0$ only on a finite set --
one uses the notions of ``infinitesimal measures'' and ``enhanced \fns''
dealt with in Section \ref{S:CT}. (These notions are akin to\NOT{linked with}
the concept of {\em Palm measure}, standard in the theory of point
processes\NOT{-- see \S\ref{SS:Palm}}.)
Thus, ``infinite constants'' which are mock Radon-Nikodym derivatives
w.r.t.\ Haar measure of non-Radon invariant measures in $G$ such as the
counting measure, are uses to ``enhance'' \fns\ with zero integral to
``reveal'' their ``infinitesimal integral'' -- their integral w.r.t.\
an ``infinitesimal measure'' induced by the original (say, probability)
measure. In the above example of rotation of the circle, the return time,
having zero usual expectation, will nevertheless have ``infinitesimal
expectation'' which relates to the measure of $E$ by a Kac-like formula.
Some examples of infinitesimal measures are given, some of them having
differential-geometric aspects. As mentioned above, formulas analogous
to the discrete case, but involving ``enhanced \fns'' and ``infinitesimal
measures'' are stated with proofs analogous to the proofs in the discrete
case (modulo a ``foundation'' involving measure theory and descriptive
set theory -- see below).
A kind of ``multiple-dimensional continuous Kac arena'' for $G$ a
unimodular Lie group is given by replacing the $\bZ$- or $\bR$- case
``future until first return'' by ``the nearest point'' w.r.t.\ a
right-invariant Riemannian metric in $G$. This is applied in Section 4.\
(see below).

Working with chains in the ``continuous case'' consists of working with
measures on $G$ depending on parameters, in particular on the points of the
standard Borel space acted by the group. These measures need not be
$\sigma$-finite.
Indeed, in most examples the most important deviation from the analogy with
the discrete case is the participation of given non-$\sigma$-finite invariant
measures on $G$, such as the counting measure. Thus the treatment splits
into two parts of distinct flavour: the ``formulas'' part, with close
analogy with the discrete case, the main difference being the appearance of
``infinitesimal measures'', and the ``foundation'' part, involving
measure-theoretic and descriptive set-theoretic considerations, to which
some paragraphs are devoted.

A feature of the language of vertices of chains is that its notions
are independent of the measure (as long as the action is measure-preserving,
in other words, the measure is invariant), and are defined using just the
action of the group on the Borel space
(usually assume standard). For example, Kac's formula --
the fact that the integral of the return time to a Borel set $E$ is equal
to the measure of the saturation of $E$, holds for any invariant measure,
since the return time (defined as zero outside $E$) and the characteristic
\fn\ of the set of points that visited $E$
(a set differing from the saturation of $E$ in a set null for every
invariant probability measure) are equidecomposable.

Section \ref{S:ED} tries to formulate (in the discrete case) reverse
implications: from equality of integral w.r.t.\ a comprehensive collection
of invariant measures to being equidecomposable. We insist on the
equidecomposability being via {\em non-negative} \fns.
Some theorems answering this are formulated and proved.
In these theorems topological
assumptions are made: it is assumed that the group acts continuously
on a compact or locally compact space and the \fns\ and chains are assumed
upper or lower semi-continuous.
(It is a well-known fact, proved by Varadarajan,
 that any standard Borel space acted in a Borel
manner by a 2nd-countable locally compact group $G$ can be
embedded as an invariant Borel subset in a compact metric space on
which $G$ acts \cts{ly}\NOT{ -- see \S\ref{SS:CP}}.) Special
attention is given to the special case of equidecomposability when
one \fn\ is the finite average of translates of another. While
most results are formulated for any acting group, in some matters
{\em amenable} groups behave better. This connects with the
Banach-Tarski paradox and Tarski's theorem, and also with
``accumulating averages'' ergodic theorems for general acting
groups related to weak compactness and Ryll-Nardzewski's
fixed-point theorem\NOT{(see \S\ref{SS:Av} and \S\ref{SS:ERG})}.
The proofs in Section \ref{S:ED} use functional-analytic methods,
mainly convex separation, sometimes taking the form of infinite
dimensional minimax theorems analogous to Von Neumann's minimax
theorem in Game Theory.

Section \ref{S:EDCT} deals with some aspects of a \cts\ counterpart of
the subject of Section \ref{S:ED} -- equidecomposable
{\em enhanced \fns}
as independent of the original invariant measure. While many questions are
raised, it is shown that the relation of equidecomposability is transitive,
if some restrictions on the chains are imposed (``tame chains'', which
include most chains that we deal with in this work).
It is also proved that in the case of unimodular {\em Lie groups} the
original measure can be reconstructed from the infinitesimal measure and
a related fact about finding \fns\ enhanced by a given measure and
equidecomposable with, say, a usual \fn.
This is done using, for Lie groups $G$, sets with
``discrete intersection with orbits'' generalizing the Ambrose-Kakutani
way of making a general flow a ``flow under a function''.
To this known generalization of Ambrose-Kakutani (proved by
Kechris\NOT{\cite{Kechris}}, and by Feldman, Hahn and Moore\NOT{
\cite{FeldmanHahnMoore}}) a proof is given (for Lie acting groups) using
differential-geometric notions.

in Section \ref{S:APL} it is shown how our language may be used in proofs of
results essentially the same as\NOT{in \cite{Helmberg} and
\cite{AaronsonWeiss}} those of Helmberg and of Aaronson and Weiss,
which relate to Kac's formula. Our treatment of the latter shows how the
elementary (in the discrete case) method of vertices of chains can
sometimes replace the pointwise ergodic theorem for multi-dimensional
groups.\NOT{which is both deep and tricky.} This language is also used
to get a proof of classical limit theorems of Renewal Theory%
\footnote{I am indebted to Prof.~Jon~Aaronson for his suggestion to apply
the language of this work to Renewal Theory.}%
, which proceeds
completely analogously for discrete and continuous Time.

In the appendices some small theories needed or connected to the above are
expounded.

\newpage
\section*{Notations} \label{S:NT}
\addcontentsline{toc}{section}{Notations}
\NOT{The following abbreviations will be used:}

\NOT{\fn\ for function, \cts\ for continuous \nbd\ for neighborhood}

A group $G$ (generally nonabelian)
will be acting on a set $\OM$, always on the left, the action denoted by
$(x,\om)\mapsto x\om=T^{x}\om,\:x\in G,\om\in\OM$. Thus it is assumed
that $e\om=\om$, $x(y\om)=xy\om$. A set endowed with such an action of $G$
will be referred to as a $G$-set.

The unit element of $G$ is denoted by $e$, but we
shall also write $0$ for the unit element (even in the non-abelian case)
when it is considered as a ``point'' in $G$, not as an ``acting agent''.

Naturally, $G$ acts (on the left) on functions $f$ from $\OM$ to, say, $\bR$,
by $xf(\om)=f(x^{-1}\om)$; Thus $\langle xf,x\om\rangle=\langle
x,\om\rangle$. Further, $G$ acts (on
the left) on functionals on invariant spaces of such functions, e.g.\ on
probability measures $\mu$ on a compact space on which $G$ acts by
homeomorphisms. Thus $\langle x\mu,f \rangle=\langle \mu,x^{-1}f\rangle$


   The {\em orbit} of an $\om\in\OM$ is the set
    $$\{x\om:x\in G\}\subset\OM$$
   For $E \subset \OM$, $f$ a function on $\OM$ and $\om \in \OM$,
\NOT{ the {\em orbit} of $E$ and $f$ are defined as:}%
   define:
   $$\begin{array}{c}
    \cO_\om E \NOT{= E_\om} := \{x \in G : x\om\in E\}\subset G\\
    \cO_\om f := x\mapsto f(x\om), \mbox{ thus } \cO_\om f:G\to
    \end{array}
   $$

Let $E\st\OM$. The {\EM saturation}%
\NOT{ or {\EM span}}
of $E$ is defined by $\spn E:=\{\om\in\OM:\cO_\om E\ne\es\}=\cup_{x\in G}xE$.

Often the group will act on a measure space $(\OM,\cB,\mu)$, usually in a 
{\EM measure-preserving} way, i.e.\ $\forall x\in G$ $\forall E\in\cB$ also
$xE\in\cB$ and $\mu(xE)=\mu(E)$.

A {\EM probability} (measure) space is a measure space with total mass 
$1$. When dealing with such a space, probabilistic terminology will be
used. For example, the {\EM expectation} of a scalar- or vector-valued
\fn\ (alias {\EM stochastic variable} alias {\EM random variable}) on
$\OM$ is just its integral.

A {\EM null set} in a measure space $\OM$ 
is a measurable subset $E\subset\OM$ with measure $0$.
A {\EM conull set} is a set whose complement is null.

A measure is {\EM complete} if every subset of a null set is a measurable,
hence a null set. The {\EM completion} of a measure is the unique (complete)
measure extending $\mu$ to the $\sigma$-algebra generated by the
$\mu$-measurable sets and the subsets of $\mu$-null sets.

For a group $G$ acting on a measure space $(\OM,\cB,\mu)$,
a measurable $E\subset\OM$ is {\EM almost-invariant} if 
$E\triangle xE$ is null for every $x\in G$. If $G$ is countable, a set
is almost-invariant iff it differs from an invariant set by a null set.

A {\EM Borel structure} in a set $\OM$ is a $\sigma$-algebra of subsets of
$\OM$. A set with a Borel structure in it is called a {\EM Borel space} and
members of the $\sigma$-algebra are referred to as {\EM Borel sets}.
When a topological space is considered as a Borel space, it is understood,
unless otherwise stated, that the Borel sets are the usual ones, i.e. the
members of the $\sigma$-algebra generated by the open sets.

A mapping between two Borel spaces is {\EM Borel} if the preimage of every
Borel set is Borel.

Our Borel spaces will usually be {\EM standard} -- see \S\ref{SS:PRD}.

A {\EM standard measure space} is a standard Borel space with a measure
on the Borel $\sigma$-algebra, sometimes the {\em completion} of such
measure is understood. When the measure is probability, we speak of a
{\EM standard probability space}.

As usual, in a product of two standard Borel spaces $\OM_1$ and $\OM_2$
one takes as Borel structure the $\sigma$-algebra generated by the
rectangles with Borel sides, i.e.\ the products $E_1\times E_2$,
$E_1\st\OM_1$, $E_2\st\OM_2$ Borel sets. 

An action of a group $G$ with a standard Borel structure
(say, a 2nd-countable locally compact group) on a standard Borel space
$\OM$ will be called {\EM Borel action} if the mapping 
$(x,\om)\mapsto x\om:G\times\OM\to\OM$ is Borel.
It will be said then that $G$ acts {\EM in a Borel manner} on $\OM$,
and we shall speak of a {\EM $G$-standard space} or a
{\EM standard $G$-space}.

{\EM Discrete} subsets of topological spaces are assumed {\em closed}.

A {\EM clopen} set is a simultaneously closed and open set.

A {\EM meager} set in a topological space is a set contained in a countable
union of closed sets with empty interior.

A {\EM compact (locally compact) space} is a Hausdorff compact
(locally compact) topological space. When a topological group (possibly
discrete) $G$ acts on it so that $(x,\om)\mapsto x\om:G\times\OM\to\OM$
is \cts, we speak of a {\EM $G$-compact space ($G$-locally compact
space)} or a {\EM compact $G$-space}.

\NOT{
A {\EM convex compact} space is a {\em non-empty} convex compact subset of
a (Hausdorff) locally convex space. When a group $G$ acts on it by affine
ransformations, we speak of a {\EM $G$-convex compact space} or a
{\EM convex compact $G$-space}.
}

Similarly for other categories.

We adopt the Bourbaki notation for open and half-open intervals in a
totally ordered set. Thus
$$]a,b[\,=\{x:a<x<b\}\qquad]a,b]=\{x:a<x\le b\}$$

Some other notations will be
$$\begin{array}{rl}
     \bZ &\mbox{ the set of integers}\\
     \bZ^+&\mbox{ the set }\{0,1,2,\dots\}
           \mbox{ of nonnegative integers}\\
     \bN &\mbox{ the set }\{1,2,\dots\}
          \mbox{ of positive integers}\\
     \bR &\mbox{ the set of real numbers}\\
     \bR^+&=\{x\in\bR:x\ge 0\}\\
     \BAR{\bR^+}&=\bR^+\cup\{+\I\}\\
      \bT &\mbox{ the circle group }\bR/\bZ\\
    \#\,E &(\in\bZ^+\cup\{\I\})
           \mbox{ the number of elements of a set }E\\
      E^c &\mbox{ the complement of the set }E\\
   \triangle &\mbox{ symmetric difference of sets}\\   
      1_E &\mbox{ the characteristic function of the set }E\\
   \BAR{E}&\mbox{ is the closure of the set }E\\
   E^\circ&\mbox{ is the interior of the set }E\\
   \partial{E}&\mbox{ is the boundary of the set }E\\
     \delta_x &\mbox{ the Dirac measure at }x\\
     \cnt=\cnt_X &\mbox{ the counting measure on the set }X\\
     \Pr &\mbox{ the probability of an event}\\
     \bE &\mbox{ the expectation of a stochastic variable}\\
      {}  *  &\mbox{ convolution}\\
     \land   &\mbox{ logical ``and''; exterior multiplication;
                     lattice operation }\min\\
     \lor    &\mbox{ logical ``or''; lattice operation }\max\\
     \&      &\mbox{ logical ``and''}\\
         V^* &\mbox{ the dual of a (normed) vector space }V\\
     \LA,\RA &\mbox{ the result of an element of a dual space
                     applied to an element of a space}\\
            \cC(\OM)&\mbox{ the space of \cts\ \fns\ on a compact }\OM\\ 
          \cC^+(\OM)&\mbox{ the cone of nonnegative \fns\ in }\cC(\OM)\\
\cC_{00}=\cC_{00}(\OM)&\mbox{ the space of \cts\ \fns\ with compact support on
                       a locally compact }\OM\\
          \cC^+_{00}&\mbox{ the cone of nonnegative \fns\ in }\cC_{00}\\
    \L^1(\OM)&\mbox{ the space of integrable \fns\ on a measure space }\OM\\
 \L^{\I}(\OM)&\mbox{ the space of bounded measurable \fns\ 
                     on a measure space }\OM\\ 
      \hat{G}&\mbox{ the dual group of a locally compact abelian group }G\\
      \hat{f}&\mbox{ the Fourier transform of the \fn\ }f\\
       \Delta&\mbox{ the modular \fn\ of a locally compact group}
\end{array}$$

A notational convention, explained in the text, is:

When (see \S\ref{SS:EnhInf})
$\la$ is a left Haar measure on the acting group $G$ and $\La$ some other
(not necessarily $\sigma$-finite) $\Delta$-right invariant measure,
$f\frac{d\La}{d\la}$ denotes the \fn\ $f$ ``enhanced'' by $\frac{d\La}{d\la}$;
$\frac{d\La}{d\la}\mu$ denotes the ``infinitesimal measure'' equal to
the measure $\mu$ in $\OM$ enhanced by $\frac{d\La}{d\la}$

\NOT{
In \S\ref{SS:Nearest}, \S\ref{SS:CntCls} and \S\ref{SS:Rnwl}
For $\om\in\OM$, let $\pi_E(\om)\in G$ be the unique nearest point
(w.r.t.\ the Riemannian metric) to $0$ in $\cO_\om E$, if it exists.
If there is no unique nearest point, $\pi_E(\om)$ is undefined.

For $\om\in E$, let $P(\om)$ be the set of $x\in G$ s.t.\ $0$ is the unique
nearest point (w.r.t.\ the Riemannian metric) to $x$ in $\cO_\om E$.
}

We use the following abbreviations:
 {\EM resp.} -- {\em respectively};
 {\EM s.t.} -- {\em such that} or {\em so that};
 {\EM w.r.t.} -- {\em with respect to};
 {\EM w.l.o.g.} -- {\em without loss of generality};
 {\EM t.f.a.e.} -- {\em the following are equivalent};
 {\EM a.e.} -- {\em almost everywhere};
 {\EM a.a.} -- {\em almost all};
 {\EM a.s.} -- {\em almost surely};
 {\EM p.m.} -- {\em probability measure};
 {\EM i.i.d.} -- {\em independent identically distributed};
 {\EM u.s.c.} -- {\em upper semi-continuous};
 {\EM (b.)l.s.c.} -- {\em (baire) lower semi-continuous}.

\newpage
\setcounter{section}{0}
\section{The Discrete Case: Chains, Hypergraphs, Expectation of Vertices}
\label{S:DS}
\subsection{Equidecomposability and Chains, The VE and HG Theorems}
\label{SS:EqiChn}
The impetus for this work has been an attempt to generalize Kac's formula
\cite{Kac} to more general acting groups, in particular to ``multi
dimensional'' groups and ``continuous'' groups.

Kac's formula can be stated as follows: for a measure-preserving invertible
transformation $T$ acting on a probability space $(\OM,\cB,\mu)$, and for
a measurable $E\subset\OM$, the following two \fns\ have the same
integral:
\begin{itemize}
\item
$\rho_E(\om)$, {\EM the return time to $E$}, defined for $\om\in E$ as
the first $k>0$ s.t.\ $T^k\om\in E$, as $\I$ if $\om\in E$ and
$T^k\om\notin E, k=1,2,\ldots$ and as $0$ outside $E$;
\item
The characteristic \fn\ $1_{\spn E}$ of the smallest (measurable)
invariant set containing $E$ (note that for ergodic action and $\mu(E)>0$
$\spn E$ is conull. Yet $\spn E$ may be conull without ergodicity, e.g.\
take $\OM=\bT\times\bT=$ the torus, $\bZ$ acting by irrational rotation
{\em of the first coordinate} and $E=I\times\bT$, $I$ an interval).
\end{itemize}

Note that the sets
$$\{\om:\{k:T^k\om\in E\} \mbox{ has a biggest (resp.\ smallest)
element $n$}\}$$
form an infinite sequence of disjoint sets with the same measure,
hence are all null (this is Poincar\'e's Recurrence Theorem).
Therefore $\rho_E$ is finite a.e. This consideration also allows us
to replace $\spn E$ by $\spn'E$, to be defined as the set of all $\om$
that has visited $E$ in the past or present, i.e.\
\begin{equation} \label{eq:satur}
\spn'E:=\{\om\in\OM:\cO_\om E\cap\,]-\I,0]\neq\es\}
\end{equation}

It is helpful to view Kac's formula as a consequence of $\rho_E$ and
$1_{\spn'E}$ being, as we shall show, (infinitely, via measurable \fns)
{\EM equidecomposable} i.e.\ there exists a family $(f_k)_{k\in\bZ}$ of
nonnegative measurable \fns\ s.t.\
\begin{equation} \label{eq:RhoSatur}
\forall\om\in\OM \quad\rho_E(\om)=\sum_{k\in\bZ}f_k(\om), \quad
1_{\spn'E}(\om)=\sum_{k\in\bZ}f_k(T^{-k}\om)
\end{equation}


The way we choose to demonstrate this equidecomposability, which is
sometimes the best way to describe equidecomposability, is as follows:

Expand any ``sequence'' $f_k$ to a ``matrix''
\begin{equation} \label{eq:Fklom}
F_{k,l}(\om)=f_{l-k}(T^{k}\om)
\end{equation}
This matrix has the invariance property:
\begin{equation} \label{eq:Finvr}
F_{k+n,l+n}(\om)=F_{k,l}(T^n\om)
\end{equation}
and one easily sees that every ``matrix'' satisfying (\ref{eq:Finvr}) comes
from some sequence $f_k=F_{0,k}$.

Call an ordered pair $(k,l)$ of elements of the group $\bZ$ a {\em
$1$-simplex}, and consider the {\em invariant (in the above sense) $1$-chain}
$\sum_{k,l}F_{k,l}\cdot(k,l)$. This may be thought of as a chain, (i.e.\
formal sum of simplices) with
\fns\ on $\OM$ as coefficients (with the group $\bZ$ acting on
these \fns), or as a \fn\ from $\OM$ to the space of chains with
scalar coefficients. This is in accordance with one of the ways to treat
group cohomology (compare the treatment in \cite{Weil}, Ch. IX \S 3). {\em
Invariant $m$-chains} $(m=0,1,2,\ldots)$ are defined analogously, as
``sums'' of $m$-simplices $(k_0,k_1,\ldots,k_n),\,\,k_i\in\bZ$ with
$m+1$-dimensional ``matrices'' $F_{k_0,k_1,\ldots,k_m}(\om)$ of
coefficients, satisfying the invariance property
\begin{equation} \label{eq:mInvr}
F_{k_0+n,\ldots,k_m+n}(\om)=F_{k_0,\ldots,k_n}(T^n\om)
\end{equation}
or in words:
\begin{equation} \label{eq:Invr}
\begin{array}{l}
\mbox{for a simplex in the chain of }T^n\om,\\
\mbox{we have the } n\mbox{-shifted simplex in the chain of }\om.
\end{array}
\end{equation}

In particular, a single \fn\ $g$ is
``equivalent'' to an invariant $0$-chain
 $\sum_k G_k(\om)\cdot(k) = \sum_k g(T^k\om)\cdot(k)$.

Now, one easily checks that the above-mentioned property of \fns\
$g(\om)=\sum_k f_k(\om)$ and $g'(\om)=\sum_k f_k(T^{-k}\om)$ to be
equidecomposable via the sequence $f_k(\om)$ can be expressed using the
$1$-chain $F$ and the 0-chains $G$ and $G'$ corresponding to $f$, $g$ and $g'$
resp. Namely, the coefficients in $G'$ are the {\em row-sums} and
those in $G$ the {\em column-sums} of the coefficient matrix in $F$. In
other words, $G'$ and $G$ are the two {\em vertices} - the {\em source} and
the {\em target} - of $F$. To define the vertices of a $1$-chain, note that
every $1$-simplex $(k,l)$ has two vertices (these being, of course,
$0$-simplices): the source $(k)$ and the target $(l)$, and this is extended
to chains by linearity and summability. (By invariance, to check that an
invariant $0$-chain is a vertex of an invariant $1$-chain it is enough to check
one coefficient of the $0$-chain, say, the coefficient at $(0)$. This
coefficient is just the \fn\ corresponding to the $0$-chain.)

One should be warned that this use of the word ``vertices'' is at variance
with the common use in Graph Theory, where the vertices of a graph are all
the $0$-simplices occuring in the graph (thus in the above chains these
``vertices'' are all the elements of the acting group $\bZ$). Still, we stick
to our use of this term and it seems not to cause misunderstanding.

Now, the definition of invariant $m$-chain carries over to {\EM general
discrete groups} - not necessarily abelian, where if the group acts on
the left {\em right shifts} are required in (\ref{eq:Invr}), i.e.\ one
requires
\begin{equation} \label{eq:gInvr}
F_{x_0y,\ldots,x_my}(\om)=F_{x_0,\ldots,x_m}(y\om)\quad
x_0,\ldots,x_m,y\in G\;\om\in\OM
\end{equation}

Treating equidecomposable \fns\ as vertices of an invariant chain,
the (trivial) fact that they have the same integral ($=$ expectation)
is encoded in the following formulation, which, as we shall see, 
is amenable to various applications in the discrete case and can be
carried over to the ``\cts'' case:%
\NOT{We have Kac's formula as a special case of the following:}

\begin{Main}
{\EM The Vertices Expectation (VE) Theorem (Discrete Case)}
Suppose a countable discrete group $G$ acts measure-preservingly on a
probability space $\OM$, and a nonnegative (right-)invariant $m$-chain
depending on $\om\in\OM$ measurably is given.
(right-)invariance means that (\ref{eq:gInvr}) is satisfied.
Then all the vertices of the $m$-chain have the same expectation.
\end{Main}
\begin{Prf} 
an ``invariant'' $m$-chain is invariant only ``globally'' - the chain of a
single $\om$ is not invariant (shifting the chain by $x^{-1}\in G$ replaces
$\om$ by $T^x\om$) but its expectation (which is an $m$-chain,
of course not dependent on $\om$),
enjoys ``genuine'' invariance - for every simplex in the chain we have all
its right-shifts in the chain with the same coefficient. Therefore it is a
countably infinite linear combination of 
``right-diagonal chains'', (these being the
sums of all the right-shifts of one simplex). Since all the vertices of a
right-diagonal chain are the same (being equal to the $0$-chain $\sum_{x\in
G}(x)$) we are done.
\end{Prf}

Note that in the VE thm.\ the assumption that $\OM$ is probability is
superfluous. $\OM$ can have $\sigma$-finite measure. Only then one has to
speak of ``integral'' instead of ``expectation''.

\NOT{This gives another way to formulate the proof of Kac's.}

Let us return to the Kac case:

The $1$-chain $F$ we take here is, for each $\om$,
\begin{txeqn}{eq:KacHg}
\mbox{the sum (=collection) of the $1$-simplices (=arrows) $(k,l)$}\\
\mbox{emanating from some $k\in\bZ$}\\
\mbox{to the last $l\in\bZ$ s.t.\ $l\le k$ and $T^l\om\in E$.}
\end{txeqn}
The invariance (\ref{eq:Invr}), i.e. (\ref{eq:Finvr}) is clear. Also, it is
evident that the source of $F$ corresponds to $1_{\spn'E}(\om)$  and the
target to $\rho_E(\om)$ (check the coefficient at $(0)$ !
recall that the coefficient at $(0)$ of the source (target) is the \# of
arrows with source (target) $(0)$). 
Consequently, by VE,%
\NOT{these two \fns\ are equidecomposable, which proves}
we have Kac's formula.

This proof is very close to
Wright's proof which uses the Kakutani skyscraper (see
\cite{Petersen}, pp.\ 45--46).) Indeed, the arrows of our graph
``do the construction work'' of the Kakutani skyscraper by
``pushing'' to $0\in\bZ$ all the $k\in\bZ$ s.t.\ $0$ is the last
$l\in\bZ$ s.t.\ $l\le k$ and $T^l\om\in E$, thus ``piling'' all of
$\OM$ over $E$. 

In\NOT{this treatment of} the Kac case, the chain was of a special kind: for
every $\om$ we had a set of simplices, and the chain was their sum, i.e.
the coefficient matrix was the characteristic \fn\ of this set. A set
of $m$-simplices in $G$ (possibly depending on $\om$)
is called a (directed) {\EM $m$-hypergraph} on $G$
(for $m=1$ -- a (directed) {\EM graph}), and (in the discrete case)
{\em every hypergraph is (identified with) a chain}.

So, as a particular case of the VE thm., which has Kac's as a special case,
we have

\begin{Main}
{\EM The Hypergraph (HG) Theorem}
Let a countable discrete group $G$ act
measure-preservingly on a probability space $\OM$.

Suppose we are given an $m$-hypergraph on $G$, depending on $\om\in\OM$
measurably and (right-)invariantly. (invariance means: the hypergraph for
$\om$ is the right $x$-shift of the hypergraph for $x\om$, $\forall x\in G,
\om\in\OM$).

Then for $i=0,1,\ldots,m$, the expectation of the number of simplices
having $(0)$ as their $i$-th vertex is the same for all $i$'s.
\end{Main}
Again, the HG thm.\ holds when $\OM$ has $\sigma$-finite measure, speaking
of ``integral'' instead of ``expectation''.

\begin{Rmk}
In some particular cases (for general $G$), one refers, as in the Kac case,
to a measurable $E\st\OM$, and the graph (for each fixed $\om$)
has the property that every arrow ends in an element of $\cO_\om E$ and every
$x\in\cO_\om E$ is the end of a {\em unique} arrow. Then one has a
``Kac formula'': the expectation of the number of arrows ending in $0$
(that being a \fn\ of $\om$) is $\mu(\spn E)$. One such case is the treatment
of Aaronson and Weiss' ``Kac \fns'' in \S\ref{SS:AW}.
\end{Rmk}

\begin{Rmk} \label{Rmk:Groupoid}
\NOT{\subsubsection{Formulation for Groupoids, the Induced Groupoid}}%
The notion of an invariant $m$-chain is adapted to viewing a $G$-set as a
{\em groupoid}.

A set $\OM$ acted by a group $G$ may be treated as a {\em groupoid} $\Gamma$
with base $\OM$ (see \cite{Weinstein}), A groupoid $\Gamma$ with base $\OM$
being just a small category with set of objects $\OM$ in which every morphism
is invertible, $\Gamma$ being the set of morphisms. In the case of a $G$-set
$\OM$, $\Gamma$ is the set of ``arrows''
$(\om_1,x,\om_2),\:\:\om_1,\om_2\in\OM,\:x\in G,\:\om_2=x\om_1$. 
In this way the ``$\om$-independent'' character of an $x\in G$ as an ``acting
agent'' is lost, thus if $G$ acts freely the groupoid structure encodes just
the orbit equivalence relation and is isomorphic to a {\em subgroupoid} of
the groupoid $\OM\times\OM$, where a unique morphism for every
$\om_1\to\om_2$ is understood (the subgroupoids of $\OM\times\OM$
are just the equivalence relations in $\OM$).%
\footnote{%
Investigations initiated in \cite{FeldmanMoore1} and \cite{FeldmanMoore2}
have shown that for many purposes
it suffices, instead of the structure of a
$G$-Borel space ($G$ countable), to consider just the groupoid, i.e.,
in case of a free action, the orbit equivalence relation.
}
On the other hand, the groupoid structure is more flexible: for instance,
for $E\st\OM$ we always have the {\EM induced groupoid} $\Gamma_E$, 
which in case of a $\bZ$-action will correspond to the 
{\em induced transformation} (see \cite{Petersen} p.12). 

Thus, instead of $m$-simplices in $G$ one may consider $m$-simplices in
$\Gamma$, defined as morphisms from the groupoid
$\{0,1,\ldots,m\}\times\{0,1,\ldots,m\}$ to $\Gamma$. Such an $m$-simplex
corresponds to an orbit of $G$ acting on pairs consisting of an $\om\in\OM$
and an $m$-simplex in $G$. Hence, invariant $\om$-dependent $m$-chains on $G$,
as defined above, are in one-one correspondence with $m$-chains on $\Gamma$.
The {\em vertices} of such chains turn up to be\NOT{$0$-chains on $\OM$}
functions on $\OM$, corresponding to the vertices of invariant
$\om$-dependent $m$-chains on $G$. 

This gives another merit to the treatment of equidecomposability via
$1$-chains: it carries over straightly to groupoids. 

In fact, since equidecomposable \fns\ can be defined in a groupoid, although
translates of a \fn\ have no meaning, the holding of VE for $m=1$
becomes the {\em definition} of measure-preservingness of the
groupoid.
\end{Rmk}

\subsection{An Assortment of Kac-like Theorems, Discrete Case}
\label{SS:DsAsso}

In this \S, the setting is a that of group $G$ acting measure-preservingly
on a probability space $(\OM,\cB,\mu)$.

\subsubsection{One-Dimensional (Discrete) Examples}
\label{SS:Asso1D}

The VE Thm., in particular the HG Thm., specialize,
for $\bZ$ as well as for ``multi-dimensional'' groups, to many cases
which are of interest. Let us list some of them.

Recall that the source of an invariant $1$-chain $F$ is an invariant
$0$-chain, determined by its coefficient at $(0)$ which is the sum of
coefficients of all edges (arrows) in $F$ with source $(0)$; similarly for
targets of $1$-chains and vertices of $m$-chains.

\begin{enumerate}
\item
\label{it:InduFn}
$G=\bZ$; $E\in\cB$; given a measurable nonnegative \fn\ $f$ on $E$.

the 1-chain $F =$ the sum of arrows $(k,l)$
s.t.\ $T^k\om\in E$ and $l$ is the first $l<k$ with $T^l\om\in E$,
with coefficient $f(T^k\om)$.

Then: If $\om\in E$, the (coefficient of) the source at $(0)$ is $f(\om)$,
while the coefficient of target at $(0)$ is $f(T_E(\om))$
where $T_E$ is the {\EM induced transformation} (see \cite{Petersen} p.12).
For $\om\notin E$ both the source and the target at $(0)$ are $0$.

By VE $\int_E f(\om)=\int_E f(T_E(\om))$, so we have proved that
{\em the induced transformation is measure-preserving}.

\item
\label{it:KacDist}
$G=\bZ$; $E\in\cB$. Let $n\in\bZ^+$. As a variation on (\ref{eq:KacHg}),
take the graph:
 $$F = \{(k,l):k-l=n\; ,
               T^l\om\in E\; ,
               \forall j\in]l,k]\; T^j\om\notin E\}
 $$
Then: there is one arrow with {\em target} $(0)$ if $\om\in E$ and
$n<\rho_E(\om)$, otherwise there is no such arrow.

Define the {\EM arrival time} $\xi_E(\om)$ of an $\om\in\OM$ as the
minimal $k=0,1,\ldots$ s.t.\ $T^k\om\in E$ (and as $+\I$ if $\exists$
no such $k$). In particular
$\xi_E(\om)=0\Leftrightarrow\om\in E$.
$\xi_E$ is finite a.e.\ on $\spn E$ and $+\I$ outside $\spn E$.

Then there is one arrow with {\em source} $(0)$ if $\xi^{(-)}_E(\om)=n$,
otherwise none, when {\EM a superscript $(-)$} denotes entities referring to
{\EM the inverse action} $T^{-1}$.

so we conclude the following strengthening of Kac's:

\begin{txeqn}{eq:KacDist1}
\mbox{The probability that $\rho_E>n$ ($n=0,1,\ldots$)}\\
\mbox{is equal to the probability that $\xi_E^{(-)}=n$.}
\end{txeqn}

So, knowing the distribution \fn\ of one of $\rho_E$, $\xi_E^{(-)}$ gives us
that of the other one.

\item
\label{it:InvDist}
Invoking the inverse action in (\ref{eq:KacDist1}) can be avoided:

Consider the following two graphs:
$$
\begin{array}{c}
F'=\{(k,l)\in\bZ^2: k-l=n,\; T^k\om,T^l\om\in E,
      \forall j\in]k,l[ \; T^j\om\notin E\}\quad n=1,2,\ldots\\
F''=\{(k,l)\in\bZ^2: k-l=n,\; \forall j\in[k,l]\; T^j\om\notin E\}\quad
    n=0,1,2,\ldots
\end{array}
$$
The target and source of $F'$ at $(0)$ give the probabilities
for $\rho_E=n$
and for $\rho^{(-)}_E=n$, resp., while the target and source of
$F''$ at $(0)$ give the probabilities for $\xi_E>n$ and for
$\xi^{(-)}_E>n$, resp. By HG, we have:

\begin{Prop}\label{Prop:InvDist}
$\rho_E$ and $\rho^{(-)}_E$ have the same distribution, similarly
$\xi_E$ and $\xi^{(-)}_E$ have the same distribution.
\end{Prop}

Hence we may replace (\ref{eq:KacDist1}) by the following strengthening
of Kac's:

\begin{Prop}\label{Prop:KacDist}
The probability that $\rho_E>n$ ($n=0,1,\ldots$)
is equal to the probability that $\xi_E=n$.
\end{Prop}

\item
\label{it:Ss}
To obtain another formulation of (\ref{eq:KacDist1}),
let $s:\bZ^+\to\BAR{\bR^+}$ be a nonnegative sequence,
let $S(n)=\sum_{0\le k<n}s(k)$.

Take the chain with same arrows as in (\ref{eq:KacHg}), i.e. the $(k,l)$
with $l$ the last $l\le k$ s.t.\ $T^l\om\in E$ , but with coefficients
$s(k-l)$.

The source at $(0)$ $=$ $s(\xi^{(-)}_E(\om))$; the target at $(0)$ $=$
$S(\rho_E(\om))$.

Using VE and arguing as in \ref{it:InvDist}, one obtains:

\begin{Prop} (another formulation of Prop. \ref{Prop:KacDist} and
\ref{Prop:InvDist}) \label{Prop:SsDist}
Let $s:\bZ^+\to\BAR{\bR^+}$ be a nonnegative sequence.
Let $S(n)=\sum_{0\le k<n}s(k)$. Then
$s(\xi_E(\om))$, $s(\xi^{(-)}_E(\om))$ $S(\rho_E(\om))$ and
$S(\rho^{(-)}_E(\om))$
have the same expectation.
\end{Prop}

Kac's is the case $s\equiv 1$.

\begin{Cor}
If $\spn E$ is conull, then for any $p\ge 0$,
$\rho_E\in L^{p+1}(\OM)\Leftrightarrow \xi_E\in L^p(\OM)$.
\end{Cor}

\begin{Rmk} \label{Rmk:1Line}
Since $\xi_E$ is obviously in $L^0(\OM)$, we have
$\rho_E$ is in $L^1(E)$ which is included, of course, in Kac's.
One cannot say anything further, because for any $g:[0,1]\to\bZ^+$ with
integral $1$ one can construct $\OM$, $E$ with $\spn E=\OM$ and $\rho_E$
having the same distribution as $g$. This is done using the discrete
``flow under a \fn'' construction
(see \cite{Petersen}, p.~11). So $\xi_E$ need not be integrable.

In case $\xi_E$ is integrable, equivalently,
by the above Corollary, $\rho_E$ is in $L^2$,
one can give the following proof to Kac's formula \NOT{($T=T^1$)}:
\begin{Prf}
Assuming $E$ conull and $\xi_E$ integrable,
Kac's formula follows from the equality:
   $$\rho_E(\om)=\xi_E(T\om)+1-\xi_E(\om)$$
\end{Prf}
\end{Rmk}

\item
\label{it:sf}
A combination of items \ref{it:Ss} and \ref{it:InduFn}:

Let $E\in\cB$ with $\spn E$ conull,
$f:\OM\to\bR^+$ measurable and $s:\bZ^+\to\BAR{\bR^+}$,
as in item \ref{it:Ss}.
 
The chain: our ``Kac'' set of simplices, namely the $(k,l)$ with 
$l$ the last $l\le k$ s.t.\ $T^l\om\in E$, but with coefficients
$f(T^k\om)s(k-l)$.

The source at $(0)$: $f(\om)s(\xi^{(-)}\om)$.

The target at $(0)$: for $\om\in E$, $\sum_{0\le k<\rho(\om)}f(T^k\om)s(k)$;
$0$ outside $E$.

And one obtains:
\begin{equation}
\int_{\OM}f(\om)s(\xi^{(-)}\om)=
\int_E\sum_{0\le k<\rho(\om)}f(T^k\om)s(k)
\end{equation}

\item
\label{it:sSf}
Another combination of items \ref{it:Ss} and \ref{it:InduFn}:

Let $E\in\cB$ with $\spn E$ conull, $f:E\to\bR^+$ measurable and
let $s:\bZ^+\to\BAR{\bR^+}$ and $S(n)=\sum_{0\le k<n}s(k)$ be as in
item \ref{it:Ss}.

The chain: similarly to item \ref{it:sf} -- the ``Kac'' set of simplices:
the $(k,l)$ with $l$ the last $l\le k$ s.t.\ $T^l\om\in E$, with
coefficients $f(T^l\om)s(k-l)$.

The source at $(0)$: $f\LP T^{-\xi^{(-)}\om}(\om)\RP s(\xi^{(-)}\om)$.

The target at $(0)$: for $\om\in E$, $f(\om)S(\rho(\om))$;
$0$ outside $E$.

And one obtains:
\begin{equation}
\int_{\OM}f\LP T^{-\xi^{(-)}\om}(\om)\RP s(\xi^{(-)}\om)=
\int_E f(\om)S(\rho(\om))
\end{equation}

\begin{Rmk} \label{Rmk:DsFlow}
For $\bZ$-action and $\spn E$ conull, the system $(\om,\mu,T)$ can
be recovered from the system $(E,\mu(\;|E),T_E)$ (the measure is conditional
probability w.r.t.\ $E$) and the function $\rho_E:E\to\bN$, via the discrete
analog of the ``flow under a \fn'' construction
(see \cite{Petersen}, p.~11). Note that by Kac's, $\rho_E$ has integral
$1/\mu(E)$ on $E$, so we have a kind of reciprocity between the integral of
the given function on one side and the measure of the given set on the other
side.
\end{Rmk}

\item \label{it:TwoSets}
\begin{Prop} \label{Prop:TwoSets}
-- {\EM return and arrival to two sets}
Let $\bZ$ act in a measure-preserving manner on $(\OM,\mu)$.
Let $E_1\st\OM$ and $E_2\st\OM$ be measurable. 

Then:

\begin{equation} \label{eq:TwoSets}
\begin{array}{l}
\int_{E_2\cap\{0<\xi_{E_1}\le\rho_{E_2}\}}\LQ\xi_{E_1}-1\RQ=
\int_{E_1\cap\{0<\xi^{(-)}_{E_2}\le\rho^{(-)}_{E_1}\}}\LQ\xi^{(-)}_{E_2}-1\RQ
=\\=
\mu\LQ\spn E_1\cap\spn E_2\cap
\{0<\xi_{E_1}\le\xi_{E_2}\}\cap\{0<\xi_{E_2}^{(-)}\le\xi_{E_1}^{(-)}\}\RQ$$
\end{array}
\end{equation}
\end{Prop}

\begin{Prf}
Take the ($\om$-dependent) $2$-hypergraph consisting of the
$2$-simplices
$(k_0,k_1,k_2)\in\bZ^3$ s.t.\ $k_2<k_0<k_1$,
$[k_2,k_1[\cap\cO_\om E_2=\{k_2\}$ and $]k_2,k_1]\cap\cO_\om E_1=\{k_1\}$.
The three expressions in (\ref{eq:TwoSets}) are the expectations of its
$2$-, $1$- and $0$-th vertices resp.
\end{Prf}

One can try to formulate other ``variations'' on this theme.

\end{enumerate}

\subsubsection{The Common Distribution of Repeated Return and Arrival}

We remain in the case of $\bZ$-action on a probability space $(\OM,\mu)$.

Let $E\subset\OM$ be measurable.
Let $T_E$ be the induced transformation (see \S\ref{SS:Asso1D}).
For $\om\in E$, we have its
{\em return time}
$\rho^{(1)}(\om)=
\rho_E(\om)$, its {\em 2nd return time} $\rho^{(2)}(\om)=
\rho_E(T_E\om)$, its {\em 3rd return time} $\rho^{(3)}(\om)=
\rho_E(T_E^2\om)$ etc.
The same for the inverse transformation $\rho^{(-)(1)(\om)}=
\rho^{(-)}_E(\om)$, $\rho^{(-)(2)}(\om)=
\rho^{(-)}_E(T_E^{-1}\om)$ etc.

For any $\om\in\OM$, we have its {\em arrival time} $\xi^{(1)}(\om)=
\xi_E(\om)$, then the {\em subsequent return time} $\xi^{(2)}(\om)=
\rho_E(T^{\xi_E(\om)}\om)$, the {\em 2nd subsequent return time} etc.
The same for the inverse transformation.

What can be said about the joint distribution of these four sequences of
stochastic variables? This is the theme of the following theorem,
which is a direct application of HG (except item f.)
(\cite{Kastelyn} gives formulas describing various aspects of these
distributions. The repeated return
times are discussed in \cite{Breiman}, Ch.~6).

\begin{Thm} \label{Thm:JointDist}
Retaining the above notations,
\begin{itemize}
\item[a.]
Let $m>0$, $r_1,\ldots,r_m>0$, $0\le p\le m$, be integers.
\BER{l}
\mu\{\om\in E: \rho^{(1)}=r_1,\rho^{(2)}=r_2,\ldots,\rho^{(m)}=r_m\}=\\
=\mu\{\om\in E:
\rho^{(-)(1)}=r_m,\rho^{(-)(2)}=r_{m-1},\ldots,\rho^{(-)(m)}=r_1\}=\\
=\mu\{\om\in E:
\rho^{(1)}=r_{p+1},\rho^{(2)}=r_{p+2},\ldots,\rho^{(m-p)}=r_m,
\rho^{(-)(1)}=r_p,\rho^{(-)(2)}=r_{p-1},\ldots,\rho^{(-)(p)}=r_1\}
\EER
Denote this value by $P[r_1,r_2,\ldots,r_m]$.

\item[b.]
Let $m>0$, $r_1\ge0,r_2\ldots,r_m>0$, be integers.
\BER{l}
\mu\{\om\in E: \rho^{(1)}>r_1,\rho^{(2)}=r_2,\ldots,\rho^{(m)}=r_m\}=
\mu\{\om\in\OM: \xi^{(1)}=r_1,\xi^{(2)}=r_2,\ldots,\xi^{(m)}=r_m\}
\EER
\item[c.]
Let $m,m'>0$, $r_1\ge0,r_2\ldots,r_m>0$, $r'_1\ge0,r'_2,\ldots,r'_{m'}>0$, 
be integers.
\BER{l}
\mu\{\om\in\OM: \xi^{(1)}=r_1,\rho^{(2)}=r_2,\ldots,\xi^{(m)}=r_m,
\xi^{(-)(1)}=r'_1,\xi^{(-)(2)}=r'_2,\ldots,\xi^{(-)(m')}=r'_{m'}\}=\\
=P[r'_{m'},\ldots,r_2',r_1+r'_1,r_2,\ldots,r_m]
\EER
(note that this depends only on the {\em sum} of the values $r_1$ and $r_1'$
for the arrival and inverse arrival times.)

\item[d.]
Each one of the joint distribution of the $\rho^{(j)}$ or the joint
distribution of the $\xi^{(j)}$ determines the joint distribution of all
the $\rho^{(j)}$'s, $\xi^{(j)}$'s, $\rho^{(-)(j)}$'s and $\xi^{(-)(j)}$'s.

\item[e.]
Assume $\mu(E)>0$, $\mu(\spn E)=1$.

The joint distribution of $\rho^{(1)},\rho^{(2)},\ldots$, w.r.t.\ the 
conditional probability in $E$, is stationary (that is, the image measure
that they define in $\bN^\bN$ is invariant w.r.t.\ the shift 
$(\alpha_j)\mapsto(\alpha_{j+1})$). Note that by Kac', each $\rho^{(j)}$
has expectation $1/\mu(E)$.

\item[f.] (\cite{Breiman} {\S}6.10)
The statement of e.\ is the only restriction on the joint distribution of the
$\rho^{(j)}$, namely, any shift-invariant measure on $\bN^\bN$ with 
integrable coordinates can be realized as the image of 
$\rho^{(j)}$'s as in e.
\end{itemize}
\end{Thm}

\begin{Prf}
Formulas a.-c.\ follow from HG:

For a., the
hypergraph (depending on $\om$) consists of all $m$-simplices
$(k_0,k_1,\ldots,k_m)$,
$k_0<k_1<\ldots<k_m$ with $k_1-k_0=r_1, k_2-k_1=r_2,\ldots, k_m-k_{m-1}=r_m$,
$T^{k_j}\om\in E, 0\le j\le m$, and $T^j\om\notin E$ for $j$ between any
two adjacent $k_j$'s. Now apply HG.

The hypergraph for b.\ is defined exactly the same, only here
it is assumed $T^{k_i}\om\in E, 0<i\le m$ but $T^{k_0}\om\notin E$.
HG gives:
\BER{l}
\mu\{\om\in E: \rho^{(-)(1)}>r_1,\rho^{(1)}=r_2,\ldots,\rho^{(m-1)}=r_m\}=
\mu\{\om\in\OM: \xi^{(1)}=r_1,\rho^{(2)}=r_2,\ldots,\rho^{(m)}=r_m\}
\EER
and to get b.\ use a.

For c., the hypergraph consists of the $(m+m'+1)$-simplices
$(k_{-m'},\ldots,k_{-1},k_0,k_1,\ldots,k_m)$, with
$k_j-k_{j-1}=r_j, k_{-j+1}-k_{-j}=r'_j, j>0$, $T^j\om\notin E$ for $j$
between any two adjacent $k_j$'s,
$T^{k_j}\om\in E, j\ne0$ but $T^{k_0}\om\notin E$.

d.\ follows from a.-c.

To prove e., take $p=1$ in a.\ and sum over $r_1\in\bN$, to obtain:
\BER{l}
\frac{1}{\mu(E)}
\mu\{\om\in E: \rho^{(1)}=r_1,\rho^{(2)}=r_2,\ldots,\rho^{(m)}=r_m\}=\\
=\frac{1}{\mu(E)}\mu\{\om\in E:
\rho^{(1)}=r_2,\rho^{(2)}=r_3,\ldots,\rho^{(m-1)}=r_m\}
\EER

To prove f.\ proceed by the method of \cite{Neveu}:

We are given a probability measure on $\bN^\bN$,
and a measure-preserving action of $\bZ^+$ on $\bN^\bN$ by the shift
$\LP T(\al)\RP_j:=\al_{j+1}$, $\al\in\bN^\bN$.
First replace it by the unique shift-invariant measure on $\bN^\bZ$ whose
image by the projection $\pi:\bN^\bZ\to\bN^\bN$ is the given measure.
(To get the mass of a cylinder in $\bN^\bZ$, i.e.\ a set depending on a
finite number of coordinates, write it as $T^{-k}\pi^{-1}C$ for a cylinder
$C$ in $\bN^\bN$.) 
Now, following \cite{Neveu}, note that on the space $X$ of all strictly
increasing sequences of intgers $(x_n,n\in\bZ)$ $\bZ$ acts in two commuting
ways:
first by the shift $\LP T(x)\RP_j:=x_{j+1}$ and second by the translation
$\LP S(x)\RP_j:=x_j+1$.
Each action has a Borel section, hence a standard Borel space
of orbits: for $S$ the space of orbits is identified with $\bN^\bZ$ via
$x\mapsto\al$ where $\al_j:=x_{j+1}-x_j$, and the identification commutes
with the shifts, while for $T$ the space of orbits is identified with
the space $\OM$ of all subsets of $\bZ$ which are unbounded above and below,
this identification transferring $S$ into the shift in $\OM$
$\sigma\mapsto\sigma+1$, $\sigma\in\OM$.
Now the given probability measure on $\bN^\bZ$, viewed as the space of
orbits of $X$ w.r.t.\ $S$, induces a probability measure on $X$
(the integral of a measurable \fn\ on $X$
being equal to the integral of its sums on orbits), this measure is
shift-invariant on $X$, and similarly this measure in its turn induces a
probability measure on the space $\OM$ of $T$, the latter being
shift-invariant. 
For this space $\OM$, with $T=$ the shift and
$E:=\{\sigma\in\OM,0\in\sigma\}$,
the $\rho^{(j)}_E$ will be distributed as the $\al_j$.

\NOT{
The proof of f.\ (see \cite{Breiman} {\S}6.10) can be put as follows:

We are given a measure in $\bN^\bN$,
\NOT{(with $\cB=$Borel -- $\bN^\bN$ is a
complete separable metric space),} and a measure-preserving action of $\bZ^+$
by $T=$ the shift $(\al_j)\mapsto(\al_{j+1})$.
by identifying each $\bN$-sequence $(\al_j)_j$ with the $(0,1)$ sequence
with $1$'s at the $\al_1$'th, $(\al_1+\al_2)$'th, $\ldots$,
$(\al_1+\ldots+\al_n)$'th, $\ldots$ places,
$\bN^\bN$ is identified with the set $2^\bN$,
\NOT{ (again $\cB=$Borel),}
$T$ acting by deleting the initial $0\ldots01$ segment.
Let $f:2^\bN\to\bN$ be the position of the first $1$.
Apply the discrete ``flow under a \fn'' construction
(see \cite{Petersen}, p.~11). to $2^\bN$ and $f$, to get
a measure-preserving $T$ on the measure space 
$$\OM=\LB ((\be_j)_{j>0},k), 0\le k<f(\be_j)=\min_{\be_j=1}j\RB$$
with the appropriate measure.
Replace $((\be_j)_{j>0},k)$ by a sequence
$(\be_j)_{j>-k}$ be adding a segment $0\ldots01$ of length $k$ before
$(\be_j)_{j>0}$. This makes $T$ a ``truncated'' right-shift
$(\be_j)\mapsto(\be_{j-1})$.
$T$ is not invertible. One finds an extension with invertible $T$
by taking the inverse limit $\tilde{\OM}$
of the inverse system $(\OM_q)_{q<0}$ with $\OM_q=\OM$ and the map
$\OM_{q-1}\to\OM_q$ is $T$. Clearly, $\tilde{\OM}$ can be identified with
$2^\bZ$ with right shift and some measure.
For this space $2^\bZ$, with $T=$ left shift and $E=$ the sequences with
$\be_0=1$, the $\rho^{(j)}_E$ will be distributed as the $\al_j$.
}

\end{Prf}

\begin{Rmk} \label{Rmk:KacDec}
Take $m=m'=1$ in c.\ (Thm. \ref{Thm:JointDist}). This gives for
$k>0$ and any $r,r'\ge0, r+r'=k$ that
\begin{equation} \label{eq:KacDec}
\mu\{\om\in\OM :\xi=r,\xi^{(-)}=r'\}=
\mu\{\om\in E:\rho=k\}
\end{equation}
Summing over $r>0,r'>0$ one gets $\int_E(\rho-1)=\mu((\spn E)\setminus E)$,
i.e.\ Kac's. (\ref{eq:KacDec}) is a kind of ``decomposition'' of Kac's. 

\end{Rmk}

\subsubsection{``Multidimensional'' Preordered Groups (``Relativistic'' Time)}
\label{SS:DsPOG}
Natural generalizations of Kac's formula and its ``refinements'' such as
Prop.\ \ref{Prop:KacDist} arise by letting ``past'' and ``future'' refer
to a shift-invariant preordering $\le$ (i.e.\ a reflexive and transitive
relation $\le$ with $x\le y\Rightarrow x+z\le y+z$)
in the acting group $G$ (for simplicity assume $G$ abelian).
One may say that we have ``relativistic'' time, similar to
$G=R^4=$ Space-Time of Special Relativity.
A typical example is $\bZ^d$ with the usual partial ordering (a $d$-tuple
is $\ge 0$ iff all the coordinates $\ge 0$).

So assume {\EM $G$ Abelian and $\le$ a shift-invariant preordering in $G$}.

Define for $x,y\in G$:
\BER{l}
[x,y] := \{z\in G | x\le z\le y\},
]x,y[\,:= [x,y]\setminus \{x,y\},
[x,y[\,:= [x,y]\setminus \{y\},
]x,y] := [x,y]\setminus \{x\}.
\EER

\begin{Def} \label{Def:EpoDur}
Let $\om\in\OM$.

The {\EM return epoch (r.ep.)} at $\om$ is the set\par
$\LB x\in G :
\om\in E\:,\: x>0\:,\: x\om\in E \:,\: \forall z\in]0,x[\; z\om\notin E\RB$\par
(where $y>x$ means $y\ge x\:,\: y\neq x$)

The {\EM return duration (r.du.)} at $\om$ is the set\par
$\LB x\in G :
\om\in E\:,\: x\ge0\:,\: \forall z\in]0,x]\; z\om\notin E\RB$

Thus the return epoch and duration are empty for $\om\notin E$.

The {\EM arrival epoch (a.ep.)} at $\om$ is the set\par
$\LB x\in G :
x\ge0\:,\: x\om\in E \:,\: \forall z\in[0,x[\; z\om\notin E\RB$\par

The {\EM arrival duration (a.du.)} at $\om$ is the set\par
$\LB x\in G :
x\ge0\:,\: \forall z\in[0,x]\; z\om\notin E\RB$
\end{Def}

\begin{Rmk} \label{Rmk:EpoDurZ1}
For the $\bZ$-case, the return epoch is the singleton $\{\rho_E(\om)\}$
for a.a.\ $\om\in E$ and the arrival epoch is the singleton $\{\xi_E(\om)\}$
for a.a.\ $\om\in\spn E$. The return duration and arrival duration are
the intervals $[0,\rho_E(\om)[$ and $[0,\xi_E(\om)[$, resp.
\end{Rmk}

In analogy with {\S}\ref{SS:Asso1D}, apply HG to the following graphs
depending on $\om\in\OM$ (fix $z\in G$):

$$
\begin{array}{c}
F'=\{(x,y)\in G^2: x-y=z,\; x\om,y\om\in E,\;
      \forall u\in]y,x[ \; u\om\notin E\}\\
F''=\{(x,y)\in G^2: x-y=z,\; \forall u\in[y,x]\; u\om\notin E\}\\
F'''=\{(x,y)\in G^2: x-y=z,\; y\om\in E,\; 
      \forall u\in]y,x]\; u\om\notin E\} 
\end{array}
$$

One obtains the following version of Kac's for preordered time:

\begin{Prop} \label{Prop:EpoDur}
Consider the case of the acting group $G$ being
discrete countable preordered Abelian
(with the preordering assumed, of course, shift-invariant).
Let $E\st\OM$ be measurable.

Consider the following pairs of
sets, depending on $\om$
\begin{itemize}
\item[(i)] the return epoch, and the return epoch
   w.r.t.\ the inverse action
\item[(ii)] the return duration, and the arrival epoch
   w.r.t.\ the inverse action
\item[(iii)] the arrival duration, and the arrival duration
   w.r.t.\ the inverse action
\end{itemize}

Then, for each of these pairs:

Every $z\in G$ has the same probability to belong to both members of the pair;

Consequently, the sums of a positive \fn\ $G\to\bR^+$ on both members of
the pair have the same expectation;

In particular, the cardinalities of both members of the pair have the same
expectation.
\end{Prop}
\qed

Kac's is the $\bZ$-case of the equality of the expectations of the
cardinalities of the r.du.\ and the a.ep.\ for the inverse action.

\begin{Rmk} \label{Rmk:EpoDurZ2}
In the ``one-dimensional'' $\bZ$-case, the durations determine the epochs and
vice-versa (Remark \ref{Rmk:EpoDurZ1}),
so we have complete ``symmetry in expectation'' between these
attributes of the action and of the inverse (Prop.\ \ref{Prop:InvDist}).
In the partially-ordered case, there is no such symmetry for return durations
(or arrival epochs), as is shown in the following example.
\end{Rmk}

\begin{Exm} \label{Exm:Tor}
Let $\bZ^2$ (with the usual partial ordering) act on
the ``discrete torus'' $(\bZ/n\bZ)^2$ by addition modulo $n$.

Let $E$ be the ``upper triangle''
$$E=\LB(s,t)\in\bZ^2 : 1\le s\le n,\, 1\le t\le n,\, n\le s+t\RB$$
considered as a subset of $(\bZ/n\bZ)^2$.

Let us find the cardinality $\rho(\om)$ of the return duration at
$\om=(s,t)$, $1\le s\le n$, $1\le t\le n$ and the average of $\rho(\om)$
over the ``torus'' (asymptotically as $n\to\infty$):\par
  if $(s,t)\notin E$,\quad $\rho(\om)=0$;\par
  if $(s,t)\in E$:\par
  if $s<n,\,t<n$,\quad $\rho(\om)=1$, (the r.du.\ is $\{0\}$), adding to
  $\sim \frac12 n^2$;\par
  if $s=n,\,t=n$,\quad $\rho(\om)=1$, (the r.du.\ is $\{0\}$);\par
  if $s=n,\,t<n$,\quad $\rho(\om)=n-t$ (the r.du.\ is a horizontal
  ``segment''), adding to $\sim \frac12 n^2$;\par
  similarly for $s<n,\,t=n$\par
  Thus $\sum_{\om}\rho(\om)\sim\frac32 n^2$. The average is $\sim\frac32$
  \par\medskip

Now let us do the same for $\rho^{(-)}$ -- the return duration for the
inverse action:\par
  if $(s,t)\notin E$,\quad $\rho^{(-)}(\om)=0$;\par
  if $(s,t)\in E$:\par
  if $s+t>n$,\quad $\rho^{(-)}(\om)=1$, (the r.du.\ is $\{0\}$),
   adding to $\sim \frac12 n^2$;\par
  if $s+t=n$,\quad $\rho^{(-)}(\om)=st$ (the r.du.\ is a ``rectangle''),
  adding to $\sim\frac16 n^3$;\par
  and the average is $\sim\frac16 n$.
\par\medskip
It might be interesting to estimate the averages of the cardinalities
$\zeta(\om)$, $\zeta^{(-)}(\om)$ of the arrival epochs:\par
  if $s+t\ge n$.\quad $\zeta(\om)=1$ (the a.ep.\ is $\{0\}$), adding to
  $\sim \frac12 n^2$;\par
  if $s+t<n$,\quad $\zeta(\om)=n-(s+t)+1$ (the a.ep.\ is part of a
  ``diagonal''), adding to $\sim\frac16 n^3$;\par
  and the average is $\sim\frac16 n$.
\par\bigskip

For the inverse action:\par
  if $s+t\ge n$,\quad $\zeta^{(-)}(\om)=1$ (the a.ep.\ is $\{0\}$),
   adding to $\sim \frac12 n^2$;\par
  if $s+t<n$,\quad $\zeta^{(-)}(\om)=2$ (the a.ep.\ is composed of one point
  ``below'' and one point ``to the left''), adding to $\sim n^2$;\par
   and the average is $\sim\frac23$\par
Thus our approximations agree with Prop.\ \ref{Prop:EpoDur}

\end{Exm}

\newpage
\section{The Continuous Case; Infinitesimal Measures}
\label{S:CT}
\subsection{Introduction}
For us, the ``\cts'' case means a Borel measure-preserving action of a
2nd-countable locally compact group%
   \footnote{This reference to topology is, in fact, inessential.
   A.~Weil has shown (see \cite{HalmosMEA} \S 59) that the structure of
   a 2nd-countable locally-compact topological group with Borel subsets and
   Haar measure can be equivalently defined purely measure-theoretically, as
   a Standard Measurable Group, i.e.\ a Group $G$ which is also a standard
   Borel space with a $\sigma$-finite measure on the Borel sets, 
   satisfying: $(x,y)\mapsto (x,xy)$ is a measure-preserving automorphism of
   $G\times G$ (in particular, the locally-compact topology in $G$ can be
   defined as the weakest topology making all convolutions of two $L^2$
   \fns\ \cts).}
$G$ on a standard measure space, i.e.\ a standard Borel space with a
completion of a Borel measure (for comments about standard spaces see
\S\ref{SS:PRD}).%
   \footnote{If a Polish topology is given in the standard space (compatible
   with the Borel structure) then a $\sigma$-finite complete measure on
   this Polish space is a completion of a Borel measure iff every open set
   is measurable and the measure is regular -- see \cite{BourbakiINT}.} 
As a prototype one may think of such a {\em flow} -- an $\bR$-action on a
standard probability space $(\OM,\cB,\mu)$.
One may try to formulate a simple-minded generalization of Kac's formula, but
as is well known, that would fail. For example, if $\OM$ is the circle,
our $\bR$-action is rotation, and $E$ is an interval, the return time is
$0$ except at one point, and its integral is $0$.

We shall try to remedy this situation by letting $\mu$ induce an
``infinitesimal measure'' (in the above case this turns out to be the counting
measure on the circle), and consider the integral\NOT{expectation} of the
return time w.r.t.\ this ``infinitesimal measure'' (in our example, this
integral equals $1-\mu(E)$).  
 To do this and to formulate ``\cts'' VE, HG and ``Kac'' theorems,
we apply the notion of invariant chains in the \cts\ case and use the
multi-faceted way in which a \fn\ on $\OM$ determines an invariant $0$-chain. 

\begin{Rmk}
Kac-like assertions in the ``\cts'' case appear in \cite{Helmberg}
(see \S\ref{SS:HE}) and even as early as \cite{Birkhoff} where one
works with a ``lower-dimensional'' measure alongside the usual measure.
Our ``infinitesimal measure'' has strong links to the Palm measures,
standard in the theory of stationary random measures on a locally-compact
group, in particuar stationary point processes (see \S\ref{SS:Palm}). 
\end{Rmk}

\subsection{Chains in the Continuous Case}
\label{SS:CtsChn}

\NOT{Before saying how I propose to deal with the above problem}

While simple-minded ``Kac'' fails in the \cts\ case,
The ``invariant chain'' approach can be generalized,
as follows: instead of {\em summing} over simplices, $m$-chains are formed
by {\em integrating} over them w.r.t.\ a measure on $G^{m+1}$. Thus,
while in the discrete case $m$-chains were, in fact, $m+1$-dimensional {\em
matrices} depending on $\om\in\OM$, in the \cts\ case they will be
{\em measures} $\phi$ over $G^{m+1}$ again depending on $\om\in\OM$.

For instance, a simplex may be identified with the chain that is the
$\delta$-measure on a simplex, and a general chain may be obtained by
integrating such entities over measures. (Here I deliberately ignore
possible restrictions on the measures.)

Thus, in the discrete case, the chain given by
$\sum_{k,l} F(k,l)\cdot(k,l)$ is viewed now as the measure
$\sum_{k,l} F(k,l)\delta_{(k,l)}$ (on the discrete $G^2$).

\NOT{
In order to make things clear, let me comment on the needed
measure-theoretic framework.
This subsection serves mainly as a reference.
may be skipped on first reading and returned to when needed.
}

In the sequel, a {\EM measure} on a set $X$ is merely
an $\BAR{\bR^+}=[0,\I]$-valued $\sigma$-additive \fn\ $\mu$ on
a $\sigma$-algebra $\cB$ of subsets of $X$ (whose members are called
``measurable sets''), with $\mu(\es)=0$.
A {\EM null set} is a member $E$ of $\cB$ with $\mu(E)=0$.
The measure is {\EM complete} if every subset of a null set is a null set.
Unless stated otherwise, measures are assumed complete. 

If $X$ is a 2nd-countable locally compact space, we shall usually assume
that every Borel set is measurable, but we will not assume
$\sigma$-finiteness nor that compact sets have finite measure.
So, {\em our measures need not be Radon measures}.%
\footnote{On a 2nd-countable locally compact space, a measure is Radon
iff it is finite on compact sets.}

Returning to a $G$ acting on $\OM$, 
The requirement for a chain (i.e.\ measure) $\phi(\om)$ to be {\EM invariant}
is, {\em for unimodular $G$}, the analog of (\ref{eq:gInvr}):
\begin{equation} \label{eq:PhiInvr}
\forall y\in G\quad \phi(\om)=\phi(y\om)*\delta_{(y,\ldots,y)}
\end{equation}
(recall: it is assumed that $G$ acts on the left).

On the other hand, for general $G$ we require, in view of the \cts\ VE Thm.\
in the next \S, that the chain {\em share the right-invariance of the left
Haar measure on $G$}. The left Haar measure
$\la$ on $G$ satisfies:
$$\la=\Delta(y)\,\la*\delta_y$$
Where $\Delta$ is called the {\EM modular \fn} of $G$
and is a \cts\ homomorphism
from $G$ to the multiplicative group $\bR^{+\times}$ (see \cite{Loomis}).

Thus we say that the chain $\phi$ is {\EM $\Delta$-right-invariant} if
\begin{equation} \label{eq:PhiInvr1}
\forall y\in G\quad \phi(\om)=\Delta(y)\,\phi(y\om)*\delta_{(y,\ldots,y)}
\end{equation}
Clearly for uniomodular $G$ this is just invariance (\ref{eq:PhiInvr}).

A chain (i.e.\ measure) $\phi(\om)$ depending on $\om$ will be called
{\EM measurable}, if the integral of any
{\em test \fn}, i.e.\ a nonnegative \cts\ \fn\ with compact support,
w.r.t.\ $\phi(\om)$ is measurable in $\om$.
If these integrals are Borel-measurable in $\om$, the chain will be called
{\EM Borel-measurable}.

\NOT{(Here we have in mind measures
One may single out measures
that are determined by their value on test \fns, but these need not be Radon
measure -- one may consider e.g.\ measures giving mass $\I$ to points.)
}

Also, the notion of {\EM vertices} of chains generalizes readily:
these are just the projections of the measures from $G^{m+1}$ on the $m+1$
coordinates $G$.
Indeed, they are obtained by integrating the $\delta$-measure of the
relevant vertex of the simplex w.r.t.\ the measure that the chain gives
on the simplices
(where integration of measure-valued \fns\ may be defined via test-\fns).
Similarly, one may speak of lower-dimensional {\EM faces} of an $m$-chain,
which will be projections of the measure from $G^{m+1}$ on some $G^{k+1}$,
$k\le m$.

We may also speak about the {\EM expectation} or {\EM integral w.r.t.\ $\om$}
of a measure $\phi(\om)$ on
$G^{m+1}$ depending on $\om$: this will be the measure $\eps$ on $G^{m+1}$
(now {\em not} a \fn\ of $\om$) s.t.\ the integral of any test-\fn\,
i.e.\ nonnegative \cts\ \fn\ with compact support, w.r.t.\
$\eps$ is the expectation of its integral on $\phi(\om)$.
the issue of existence and uniqueness of this expectation 
will be made more precise when theorems are formulated. Since Radon
measures on $G^{m+1}$ are determined by their value on test \fns, a Radon
measure expectation is indeed unique. Such expectation will be referred
to as a {\EM Radon expectation} of the $m$-chain.

Note that the Radon expectation of a $\Delta$-right-invariant $0$-chain must
be $\Delta$-right-invariant, i.e.\ a {\em left Haar measure}. 

If an $\Delta$-right-invariant $0$-chain has a Radon expectation, necessarily
left Haar, one readily sees (by taking monotone increasing and bounded
monotone decreasing limits with common compact support) that any nonnegative
Borel \fn\ on $G$ will have the property of test \fns, namely integrating
it commutes with taking the expectation.

\begin{Rmk} \label{Rmk:CntTstFns}
Note that one may choose a countable collection of
test-\fns\ -- non-negative \cts\ with compact support -- on $G^{m+1}$
s.t.\ {\em every test-\fn\ is a non-decreasing limit of a sequence of members
of $\cF$}. (Take as $\cF$, e.g., all positive rational combinations
of the union of non-decreasing sequences of test-\fns\ that converge to the
characteristic \fns\ of finite unions of members of a countable base to the
topology in $G^{m+1}$.)
Thus for our usual purposes it suffices to test on this countable
collection of test-\fns. Therefore if something about the measure (that
depends, say, on $\om$) holds for a.a.\ $\om$ for every fixed test-\fn\
it will hold for a.a.\ $\om$ for the measure.

Note that we may impose on $\cF$ that all its members are finite positive
combinations of \fns\ of the form $h_0(x_0)\cdots h_m(x_m)$,
the $h_i$ being test-\fns\ on $G$
(use the Stone-Weierstrass Approximation Theorem),
or that all its members are finite positive combinations of convolutions
of two test-\fns. This is sometimes useful.
\end{Rmk}

\NOT{
\begin{Rmk}
Note that in our setting, {\em Fubini's theorem need not apply} when
integrations on $G^{m+1}$ and on $\OM$ are interchanged,
since $\sigma$-finiteness is not assumed.
So one has to be cautious.
In fact, this non-Fubini will be later exploited to
our advantage, in dealing with ``infinitesimal measures''.
\end{Rmk}
}

\NOT{
\subsection{Continuous Versions, First Steps} \label{SS:CVE}
}

\subsection{The Continuous VE (CVE) Theorem}
\label{SS:CVE}

In the previous \S\ we extended the notion of chain from the discrete to the
\cts\ case.
Yet there are some important differences: 
\begin{itemize}
\item
In the \cts\ case, a hypergraph is not automatically a chain, since
there is no ``natural'' measure on a set (except, of course, the highly
massive ``counting measure'' which is usually not suitable).

\item
One cannot substitute in a measure, so the fact, holding in the discrete
case, that e.g.\ any invariant 0-chain $\sum F_k(\om)\cdot(k)$ comes from
a \fn\ on $\OM$, namely $F_0(\om)$, has no \cts\ analogue.
In fact, as we shall see,
{\em invariant $0$-chains are quite richer than \fns}.
\end{itemize}

These facts make the transition from a ``VE'' theorem to ``HG'' theorems
more involved. Yet one may formulate readily a \cts\ VE theorem.
In this theorem one does not need the assuption that $(\OM,\mu)$
is probability.

\begin{Main}
{\EM The Continuous Vertices Expectation (CVE) Theorem}
Let a 2nd-countable locally compact group $G$ act
in a\NOT{Borel and} measure-preserving manner on a\NOT{standard}
measure space $(\OM,\cB,\mu)$ ($\mu$ need not be $\sigma$-finite).

Let us be given a $\Delta$-right-invariant $m$-chain. This is a measure $\phi(\om)$
on the set $G^{m+1}$ of $m$-simplices, depending on $\om$ (note it need not
be a Radon measure), and assumed to be defined on Borel subsets of $G^{m+1}$.
Assume the dependence of $\phi$ on $\om$ is {\em measurable}, i.e.\
the integral of a fixed test-\fn\ on $G^{m+1}$ (that is, nonnegative \cts\
\fn\ with compact support) w.r.t.\ $\phi$ is measurable in $\om$.
$\Delta$-right-invariance means that (\ref{eq:PhiInvr1}) is satisfied.

Suppose that one vertex of $\phi$ has a Radon $\mu$-integral, equal to the
(left Haar) measure $\la$ on $G$. Then every vertex has the same Radon
$\mu$-integral $\la$.
\end{Main}

\begin{Prf}
Since two vertices of an $m$-chain $\phi$ are vertices of some $1$-chain
``edge'' of $\phi$ (i.e.\ projection of the measure $\phi$ on some $G^2$),
and since measurability of the chain implies measurability of every ``edge'',
we have to prove the assertion only for $m=1$.

Let $\phi$ be our $1$-chain. Thus, for $\om\in\OM$ $\phi(\om)$ is a measure
on $G^2$. Let $f:G\to\bR$ be a test-\fn. Its integrals on the two vertices
are the integrals w.r.t.\ $\phi$ of $(x,y)\mapsto f(x)$ and
$(x,y)\mapsto f(y)$, which both depend on $\om$. We are told that, say,
$$\forall f\quad
\int_{\OM}\LP\int_{G^2}f(x)\,d\phi\RP\,d\mu(\om)=\LA\la,f\RA$$
and have to prove that
$$\forall f\quad
\int_{\OM}\LP\int_{G^2}f(y)\,d\phi\RP\,d\mu(\om)=\LA\la,f\RA.$$

Let $U$ be a \nbd\ of $e\in G$, with $U^{-1}=U$.  
Let $\LP h_n(x)\RP_n$ be a countable partition of unity on $G$
(thus, $h_n\ge0$, $\sum h_n(x)=1$), where $h_n$ are test-\fns\ with
$x,y\in\supp h_n\Rightarrow x^{-1}y\in U$.

The idea is to ``slice'' $(x,y)\mapsto f(x)$ (or, if one wishes, to slice
$\phi(\om)$) into diagonal slices using
the partition of unity, and then to use the shift-invariance to move
each slice diagonally so that the sum will approximate $(x,y)\mapsto f(y)$.

We have:
$$
\LA\la,f\RA=\int_{\OM}\LP\int_{G^2}f(x)\,d\phi\RP\,d\mu(\om)=
\sum_n\int_{\OM}\LP\int_{G^2}f(x)h_n(xy^{-1})\,d\phi\RP\,d\mu(\om)
$$

Denote by $\la_n$ the functional on the test-\fns\ $f$:
$$
\LA\la_n,f\RA:=
\int_{\OM}\LP\int_{G^2}f(x)h_n(xy^{-1})\,d\phi\RP\,d\mu(\om)
$$

Recall that we have the $\Delta$-right-invariance (\ref{eq:PhiInvr1}): 
$$\forall a\in G\quad \phi(\om)=\Delta(a)\,\phi(a\om)*\delta_{(a,a)}$$
which means that for a \fn\ $F$ on $G^2$
$$\int F(x_0,x_1)\,d(\phi(\om))(x_0,x_1)=
\Delta(a)\,\int F(x_0a,x_1a)\,d(\phi(a\om))(x_0,x_1).$$

Thus for $a\in G$:
\BER{l}
\LA\la_n,f\RA=
\int_{\OM}\LP\int_{G^2}f(x)h_n(xy^{-1})\,d(\phi(\om))\RP\,d\mu(\om)=
\Delta(a)\,
\int_{\OM}\LP\int_{G^2}f(xa)h_n(xy^{-1})\,d(\phi(a\om))\RP\,d\mu(\om)=\\=
\Delta(a)\,
\int_{\OM}\LP\int_{G^2}f(xa)h_n(xy^{-1})\,d(\phi(\om))\RP\,d\mu(\om)=
\Delta(a)\LA\la_n,(x\mapsto f(xa))\RA
\EER
Thus $\la_n$ is $\Delta$-right-invariant.
Since it is positive and finite on test-\fns\ (being dominated by
$\la$), it is a left Haar measure.
Thus, for $x_n\in G$, to be chosen later, we have
\BER{l}
\LA\la,f\RA= 
\sum_n\int_{\OM}\bLP\int_{G^2}f(x_nx)h_n(xy^{-1})\,d\phi\bRP\,d\mu(\om)
=\\=
\int_{\OM}\LP\int_{G^2}\LQ f(y)+\sum_n(f(x_nx)-f(y))h_n(xy^{-1})\RQ
\,d\phi\RP\,d\mu(\om)
\EER
For a fixed test-\fn\ $f$, we have to estimate the error:
\BER{l}
\cE=\int_{\OM}\bLP\int_{G^2}f(x)\,d\phi\bRP\,d\mu(\om)-
\int_{\OM}\bLP\int_{G^2}f(y)\,d\phi\bRP\,d\mu(\om)
=\\=
\int_{\OM}\LP\int_{G^2}\LQ\sum_n(f(x_nx)-f(y))h_n(xy^{-1})\RQ
\,d\phi\RP\,d\mu(\om)
\EER
To this end ($f$ is fixed, fix $\eps>0$ and fix a relatively compact
$e$-\nbd\ $U_0$), choose the above $U=U^{-1}\st U_0$
so that $xy^{-1}\in U\Rightarrow |f(x)-f(y)|<\eps$ and choose
$x_n^{-1}\in\supp h_n$.
If $xy^{-1}\in\supp h_n$ then $x_nxy^{-1}\in U$, hence
$|f(x_nx)-f(y)|<\eps$. Thus, if $k$ is a \cts\ nonnegative \fn\ with compact
support s.t.\ $k\ge1$ in $U_0\cdot(\supp f)$, then
$xy^{-1}\in\supp h_n\Rightarrow |f(x_nx)-f(y)|\le \eps k(x_nx)$,
hence, using the left Haar'ness of $\la_n$: 
\BER{l}
|\cE|\le
\int_{\OM}\LP\int_{G^2}\LQ\sum_n|f(x_nx)-f(y)|h_n(xy^{-1})\RQ
\,d\phi\RP\,d\mu(\om)
\le\\ \le
\eps\int_{\OM}\LP\int_{G^2}\LQ\sum_n k(x_nx)h_n(xy^{-1})\RQ
\,d\phi\RP\,d\mu(\om)
=\\=
\eps\int_{\OM}\LP\int_{G^2}\LQ\sum_n k(x)h_n(xy^{-1})\RQ
\,d\phi\RP\,d\mu(\om)\le
\eps\LA\la,k\RA
\EER
Since $\cE$ does not depend on $\eps$, one concludes that $\cE=0$
and we are done.
\end{Prf}

\NOT{
One may try to investigate the case of non-unimodular groups. (is there a
counterexample?)
}

\begin{Rmk}
Note that a statement like the CVE Thm.\ cannot hold for chains with
other kind of ``right-invariance'', i.e.\ w.r.t.\ a homomorphism from $G$ to
$\bR^{+\times}$ other than $\Delta$. Indeed, such invariance would be
inherited by the vertices and their expectations, and if a statement such
as CVE holds, these expectations must be left-invariant.
That follows from the fact that replacing a $1$-chain $\phi$ by
$\delta_{(0,a)}*\phi$ preserves its kind of right-invariance,
does not change one vertex but left-shifts the other.
\end{Rmk}

\subsection{Enhanced Functions and Infinitesimal Measures}
\label{SS:EnhInf}
For the rest of Section \ref{S:CT} we place ourself in the setting of a
locally compact group $G$ acting in a Borel and 
measure-preserving manner on a standard $\sigma$-finite\NOT{probability}
measure space $(\OM,\cB,\mu)$. Although we need not assume that
the measure $\mu$ is probability, we still speak of ``expectation'' instead
of $\mu$-integral, since our main interest lies in the probability case.
Denote a left Haar measure in $G$ by $\la$.

As already mentioned, unlike the discrete case, there is no
1-1 correspondence between measurable \fns\ on $\OM$ and
$\Delta$-right-invariant $0$-chains.
Given a (nonnegative) measurable \fn\ $f$ on $\OM$,
one can still correspond to it a measurable $\Delta$-right-invariant
$0$-chain, mapping each $\om\in\OM$ to $f(x\om)\,d\la(x)$. To
check that this $0$-chain is indeed $\Delta$-right-invariant,
that is, satisfies (\ref{eq:PhiInvr1}):
\begin{equation} \label{eq:fInvr}
\Delta(y)\,\LQ f(x\cdot y\cdot\om)\,d\la(x)\RQ*\delta_y=f(x\cdot\om)\,d\la(x).
\end{equation}
Its expectation is $\bE(f)\la$, as one finds using Fubini:
indeed, if $h:G\to\bR$ is a test-\fn,
\BER{l}
\int_\OM\LQ\int_G h(x)f(x\om)\,d\la(x)\RQ\,d\mu(\om)=\\=
\int_G\LQ\int_\OM h(x)f(x\om)\,d\mu(\om)\RQ\,d\la(x)=\\=
\int_G\LQ\int_\OM h(x)f(\om)\,d\mu(\om)\RQ\,d\la(x)=
\int_Gh(x)\bE(f)\,d\la(x),
\EER
and we sometimes identify this $0$-chain with the ``ordinary'' \fn\ $f$.

But let $\La$ be some other $\Delta$-right-invariant
measure (as usual, assumed complete, but not necessarily $\sigma$-finite)
on $G$, i.e. satisfying
\begin{equation} \label{eq:DInvMea}
\forall y\in G\quad \La=\Delta(y)\,\La*\delta_y
\end{equation}
(In the case of unimodular $G$ this is just right-invariance),
s.t.\ Borel sets are $\La$-measurable. An example is 
the counting measure, in the case of unimodular $G$.
Let $f:\OM\to\BAR{R^+}=[0,\I]$.
Suppose that for a.a.\ $\om$ $x\mapsto f(x\om)$ is $\La$-measurable.
One may consider the $0$-chain mapping each $\om\in\OM$ to the measure on
$G$ $f(x\om)\,d\La(x)$. This $0$-chain will be denoted by $f\frac{d\La}{d}$
(the rationale behind this notation will be seen below).
Suppose {\em $f$ is s.t.\ $f\frac{d\La}{d}$ is measurable}.
It will be $\Delta$-right-invariant {\em for the same reason as before},
(\ref{eq:fInvr}). 
And {\em even if $f$ has expectation $0$}, $f\frac{d\La}{d}$ may
have an expectation $a\la$, $a\ge0$. Here $a>0$ is
possible since {\em Fubini need not hold for $\mu(\om),\La(x)$}.
We shall write then
\begin{equation} \label{eq:EnhFn}
  \bE\LP f\frac{d\La}{d\la}\RP=a
\end{equation}
and view the expression $f\frac{d\La}{d\la}$ as an {\EM enhanced function},
$f$ being ``enhanced'' by the {\em ``infinite constant $\frac{d\La}{d\la}$''}.
(This mock-``Radon-Nikodym derivative'' should indeed be viewed as a
``constant'' since both $\la$ and $\La$ are $\Delta$-right-invariant.)
This notation is further justified by the fact that if $\La$ happens to
be Radon, hence some (other) left Haar, $\frac{d\La}{d\la}$ is just a number
multiplier in the above.

When the Haar measure $\la$ is fixed, we identify $f\frac{d\La}{d\la}$
with the $0$-chain $f\frac{d\La}{d}$ and speak of an enhanced function as
a special case of $0$-chain.

To be more precise, one may think of an enhanced function as a positively
homogeneous mapping from the half-line of left Haar measures on $G$ to the
cone of non-negative $\Delta$-right-invariant ($\om$-dependent) $0$-chains,
denoted by
$f\frac{d\La}{d\la}$ if it maps $\la\mapsto f\frac{d\La}{d}$. Such entities
can be added and multiplied by non-negative constants and by functions $f$ on
$\OM$ (i.e.\ by $(x;\om)\mapsto f(x\om)$).
The {\em ``infinite constants''}
such as $\frac{d\La}{d\la}$ are to be defined as positively homogeneous
maps from the half-line of left Haar measures on $G$ to the cone of
(not necessarily Radon) $\Delta$-right-invariant measures in $G$.

Note that (\ref{eq:EnhFn}) means just that for {\em arbitrary} test-\fn\
$h:G\to\bR^+$ (here we may take, in our case, any nonnegative Borel
\fn) we have:
\begin{equation} \label{eq:EnhFn1}
 \int_\OM\LP\int_G h(x)f(x\om)\,d\La(x)\RP\,d\mu(\om)=a\int_G h(x)\,d\la(x)
\end{equation}
Note that by invariance we are sure that if the enhanced $f$ has a Radon
expectation at all, it will be a multiple of $\la$ so the left-hand side
of (\ref{eq:EnhFn1}) is a constant multiple of $\int_Gh\,d\la$.

\begin{Rmk}
For example, for $\bR$-action (\ref{eq:EnhFn}) means just that for
$\om\in\OM$, if we consider the $\La$ integral on Time of $f$ over an
interval in the past or future, and take the expectation for $\om$, we get
$a$ times the $\la$-length of that interval.
\end{Rmk}

\begin{Rmk} \label{Rmk:OneFn}
Note that if $f$ has a $0$-chain ``enhancement'' which is {\em measurable}
and (\ref{eq:EnhFn1}) holds for some finite $a\ge0$ and {\em one} nonnegative
integrable Borel $h$ which is bounded below away from $0$ on an open set,
then the ``enhancement'' has the expectation $a\,d\la$, hence
(\ref{eq:EnhFn1}) holds for {\em every} nonnegative Borel \fn\, i.e.\
(\ref{eq:EnhFn}) holds. (by invariance (\ref{eq:EnhFn1}) holds
for all combinations of translates of $h$ which dominate any \cts\ with
compact support, hence the expectation is Radon.)
\end{Rmk}

\begin{Rmk}
In fact, the notation (\ref{eq:EnhFn}) is somewhat misleading -- it conceals
the fact that these notions depend on the particular action of the group $G$
on $\OM$.
\end{Rmk}

In order to clarify what was said, consider the following example alluded
to before:

\begin{Exm} \label{Exm:InfMeaRT}
Take $G=\bR$, $\OM=$ the circle $\bT=\bR/\bZ$ with Lebesgue probability
measure, $\bR$ acting by rotation.
Let $E=$ a finite union of closed intervals in $\bT$. Let $f$ be the
simple-minded {\em return function} to $E$. It is $0$ except for the
finite set of the upper extremities of the intervals (where the motion
``exits'' from $E$).
Of course, $f$ has zero expectation. But consider in $\bR$
the Haar $\la=dx$ and the invariant $\La=\cnt=\cnt_\bR$ -- the counting
measure.
One easily shows that for any $f:\bT\to\BAR{R^+}$ with countable support,
$\bE\LP f\frac{d\cnt}{dx}\RP=\sum_tf(t)$. Thus, ``enhanced'' by
$\frac{d\cnt}{dx}$ our return function has always the expectation
$1-\mu(E)$ which is an instance of a ``Kac'' theorem to be proved later
(Thm.\ \ref{Thm:CtsKac}) 
\end{Exm}

As seen in the last example, taking the expectation of an ``enhancement''
of functions $f$ on $\OM$ may amount to integrating $f$ on some measure
on $\OM$ ($\cnt_\bT$ in the example) that gives positive mass to
$\mu$-null sets, yet is ``induced'' by $\mu$ (and $\La$, $\la$).
Call it an {\EM infinitesimal measure}, since it may be thought of as
specifying the ``infinitesimal mass'' of a set $\st\OM$. Denote
it by $\frac{d\La}{d\la}\mu$. ($\mu$ is ``enhanced'' by ``multiplication''
by the ``constant'' $\frac{d\La}{d\la}$.)

In order to consider such ``infinitesimal measures'' in the general
setting, proceed as follows:

We use the generation of a measure via ``preintegrable'' \fns\, as explained
in \S\ref{SS:PRE}.
The set $\Pre$ of {\em preintegrable} \fns\ $f$ will consist
of the $[0,\I]$-valued \fns\ which have an enhancement which is a
measurable $0$-chain with a Radon expectation, thus the ``enhanced'' $f$
has a finite expectation, which will be the {\em integral} of $f$.

One checks easily that the axioms 1-\ref{it:comp} in \S\ref{SS:PRE} are
satisfied (recall $\La$ is assumed complete).
Thus the ``infinitesimal measure'' $\frac{d\La}{d\la}\mu$ is defined on
some $\sigma$-algebra of subsets of $\OM$. We have:

\begin{Prop} \label{Prop:InfMea}
Every Borel subset of $\OM$ is measurable w.r.t.\ the
``infinitesimal measure'' $\frac{d\La}{d\la}\mu$.
\end{Prop}  

\begin{Prf}
We refer to \S\ref{SS:CP}. One may assume $\OM$ is an invariant
Borel subset of a compact metrizable $G$-space $\bar{\OM}$, where the
action of $G$ on $\bar{\OM}$ is \cts\ in the two variables. It suffices to
prove the assertrion for $\bar{\OM}$ instead of $\OM$. Note that $\mu$ on
$\bar{\OM}$ is not neccessarily finite thus not neccessarily Radon.

One needs to prove that for every {\em closed} $K\st\bar{\OM}$, K is
$\frac{d\La}{d\la}\mu$-measurable, i.e.\ (see \S\ref{SS:PRE}) that
$\forall f\in\Pre,\: f\cdot 1_K\in\Pre$. One needs to know that if
the $0$-chain $f(x\om)\,d\La$ is measurable with Radon expectation,
so is $f(x\om)1_K(x\om)\,d\La$. In fact, this will hold for $1_{K'}(x,\om)$
for any closed $K'\st G\times\bar{\OM}$, instead of $1_K(x\om)$.
That follows from its holding for such closed $K'$ which are finite unions
of ``rectangles'' of the form
$\mbox{closed }\times\mbox{compact }\st G\times\bar{\OM}$,
which have a general closed $K'$ as a countable decreasing intersection
(note that we are always
checking test-\fns\ with {\em compact} support in $G$). For the latter
the measurability and having Radon expectation are evident.
\end{Prf}

Note that (see \S\ref{SS:PRE})
for a $\frac{d\La}{d\la}\mu$-measurable nonnegative $f$, in particular for
nonnegative Borel \fns\ $f$,
$f$ has a finite $\frac{d\La}{d\la}\mu$-integral iff it is in $\Pre$ 
i.e.\ iff it has an ``enhancement'' which is a measurable $0$-chain with
Radon expectation, its integral being equal to that expectation.

The following proposition is easily proved:

\begin{Prop} \label{Prop:InfMea0I}
Suppose $\La_1$ and $\La_2$ are $\Delta$-right-invariant measures on $G$
s.t.\ every set of finite $\La_1$ mass has zero $\La_2$ mass. Then every
$\frac{d\La_1}{d\la}\mu$ -integrable \fn\ has zero
$\frac{d\La_2}{d\la}\mu$ -integral.
\end{Prop}

\NOT{
\begin{Rmk}
Most of what was done in this \S\ does not need the assumption that the
original measure $\mu$ in $\OM$ is probability or has finite total mass.
\end{Rmk}
}

\begin{Exm} \label{Exm:InfMea}
of ``infinitesimal measures'':
\begin{enumerate}
\item
A situation which includes Exm.\ \ref{Exm:InfMeaRT}.
$G$ -- some unimodular Lie group with Haar $\la$, $\Gamma$ -- a discrete
subgroup s.t.\ $\la$ induces on the homogeneous space $\OM=G/\Gamma$
an (invariant) probability measure $\mu$.
$G$ acts on $\OM$ by left multiplication.
$\La$ -- some right-invariant measure in $G$.

Using the test-\fn\ $1_D$ where $D$ is a fundamental domain for $\Gamma$,
one concludes that $\frac{d\La}{d\la}\mu$ is just the measure on $\OM$
transported from $\La|_D$ by $(x\mapsto x\om):G\to\OM$.

\item \label{it:InfMea-a}
Let $\bR$ act on the torus $\bT^2=(\bR/\bZ)^2$ (with Lebesgue measure $\mu$)
by $x(s,t)=(s+ax,t+bx),\:\:x\in\bR$. When $a/b$ is {\em irrational},
this is ergodic. Let $\la=dx,\La=\cnt_\bR$. $\frac{d\cnt}{dx}\mu$
measures ``the area swept by a set in unit time'', i.e.\
$\frac{d\cnt}{dx}\mu(E)$ is the area swept by $E$ during some interval
of time, with multiplicity, divided by the length of the interval.
It follows from (\ref{eq:DifFormula}) below%
\NOT{One should be able to prove (although I have not elaborated the details)}
that for a smooth curve $E$ in $\bT^2$,
$$\frac{d\cnt}{dx}\mu(E)=\int_E|-b\,ds+a\,dt|$$
The latter differential form being the contraction of the area
$2$-differential in $\OM=\bT^2$ by the vector-field induced on $\OM$
by $\frac{\partial}{\partial x}$ on $\bR$.

\NOT{
(For the proof use the fact that, as mentioned above (Remark \ref{Rmk:OneFn}),
just one test-\fn\ needs to be checked, which may be taken with small
support near $0$).
}
\NOT{
\item
One may try to find formulas of this kind for more general situations,
such as $G$ an $n$-dimensional unimodular Lie group $G$ acting on a smooth
manifold $\OM$ with probability measure given by an $m$-form $\alpha$, $E$ a
submanifold and $\La$ some lower-dimensional invariant measure in $G$ given
by Haar measures on subgroups of $G$ or by integration of forms.
(see \cite{Federer} for pertinent theory).
}

\item \label{it:DifFormula}
Similarly, let an $n$-dimensional Lie group $G$ with left Haar
measure $\la$ act smoothly and measure-preservingly (on the left) on a smooth
$m$-dimensional manifold $\OM$ with (invariant) probability measure
$\mu$ given by an $m$-differential form, also donoted by $\mu$
(actually by the absolute value of this
differential form, so that orientations do not matter),%
\NOT{``twisted'' means: ``tensored with the twisted scalars'', the
latter being the $1$-dimensional vector bundle on the manifold associated
with the principal $\{\pm1\}$-bundle of orientations and with the action
of $\{\pm1\}$ on $\bR$ by usual multiplication of a real by $\pm1$.
Twisted $k$-differential forms can be viewed as scalar functions on pairs
(point, a {\em twisted} $k$-vector at the point) linear in the latter.
A twisted $k$-vector serves as an ``infinitesimal $k$-dimensional domain''
and a twisted $k$-differential form, or more generally a function on
(point, twisted $k$-vector) just homogeneous in the latter, can be
integrated on $k$-submanifolds or more general sets 
(see \cite{Federer}). If the manifold is oriented and an orientation
is fixed, ``twisted'' things are equivalent to ``non-twisted''. 
}
Let $\La$ be an $\Delta$-right-invariant\NOT{$k$-dimensional} measure in $G$
given on submanifolds by (the absolute value of) a $\Delta$-right-invariant
$k$-differential form, also denoted by $\La$.
The left Haar measure $\la$ is given by (the absolute value of)
a left-invariant $n$-differential form, to be denoted also by $\la$.
Then one has the following formula:
\NOT{
it seems that one can prove (at least making some assumptions of
transversality etc. -- I have not elaborated this, just came to the
formula from heuristic ``physicist-like'' considerations)
-- for the proof see \S\ref{SS:DIF}:
}

{
\newcommand{\lcntr}{{|\!\raisebox{-1pt}{\underline{\mbox{ }}}\vspace{2pt}}}
Let $\gamma$ be an $n$-vector in the tangent space $T_eG$ s.t.\
$|\LA\la,\gamma\RA|=1$.
Then $\frac{d\La}{d\la}\mu$ is given on submanifolds by (the absolute
value of) the $\ell$-differential form, $\ell=m+k-n$
\begin{equation} \label{eq:DifFormula}
\LP\La\lcntr\gamma\RP\om\lcntr\mu
=\LP (\om,u)\mapsto\LA\mu,\LP\LA\La,c\gamma\RA\RP\om\land u\RA\RP
\end{equation}
where $\lcntr$ denotes contraction, $c$ is the coproduct, $u$ is
a variable $\ell$-vector in $T_\om\OM$
and applying something in the tangent space $T_eG$ to $\om$ means applying
the derivative of $(x\mapsto x\om)$ at $e$ to the ``something''.

To prove (\ref{eq:DifFormula}), let $E$ be an $\ell$-dimensional submanifold
in $\OM$. Let $\phi:G\times\OM\to G\times\OM$ be given by
$\phi(x,\om)=(x,x\om)$, thus $\phi^{-1}(x,\om)=(x,x^{-1}\om)$.
By (\ref{eq:EnhFn}) and (\ref{eq:EnhFn1}), we compute
$\frac{d\La}{d\la}\mu(E)$ as follows:

Construct the $\ell+n=m+k$ -submanifold of $G\times\OM$
$\tilde{E}:=\phi^{-1}(G\times E)$, with the projection $\pi:\tilde{E}\to\OM$
given by $\pi(x,\om)=\om$ (with fibers $\pi^{-1}(\om)$ that we sometimes
identify with subsets of $G$). By Sard's Lemma (see \cite{Schwartz} Ch.\ I),
for a.a.\ $\om$ $\pi$ is a submersion throughout the fiber $\pi^{-1}(\om)$
so that the fiber is a $k$-manifold in $G$.
To compute $\frac{d\La}{d\la}\mu(E)$, choose an open $U\st G$ with
$0<\la(U)<\I$, for each $\om$ integrate $\La$ on $\pi^{-1}(\om)\cap U$,
then integrate on $d\mu(\om)$ and divide by $\la(U)$. This iterated
integration is given by integration on $\tilde{E}\cap(U\times\OM)$ of the
(absolute value of) the $(m+k)$-differential form $\nu$ in $G\times\OM$ given
by ($\xi_i,\xi'_i\in T_xG$, $v_i\in T_\om\OM$):
\begin{equation} \label{eq:DifFormula1}
\LA\nu,(\xi_1,0)\land\ldots\land(\xi_k,0)
\land(\xi'_1,v_1)\land\ldots\land(\xi'_m,v_m)\RA=
\LA\La,\xi_1\land\ldots\land\xi_k\RA\LA\mu,v_1\land\ldots\land v_m\RA
\end{equation}
(Note that at points $(x,\om)\in\tilde{E}$ where $\pi$ is not a submersion,
i.e.\ its derivative is not onto $T_\om\OM$, $\nu|{\tilde{E}}$ vanishes.)

To proceed, the measure $\frac{d\La}{d\la}\mu$ on $E$ will be given by an
$\ell$-differential form s.t.\ integration w.r.t.\ $\nu$ will be the result
of another iterated integration, here using the projection
$\pi':\tilde{E}\to E$ given by $\pi'(x,\om)=x\om$ and taking the integration
w.r.t.\ $\la$ on the fibers
(identified with $G$ by $x\mapsto\phi^{-1}(x,\om)=(x,x^{-1}\om)$)
and then integrating w.r.t.\ that $\ell$-differential form on $E$.
This means that for $v_1,\ldots,v_\ell\in T_\om E$, taking $x=e$ and writing
the above $\gamma$ as $\gamma=\xi_1\land\ldots\land\xi_n$, $\xi_i\in T_eG$
(note that the above identification $x\mapsto(x,x^{-1}\om)$ gives, on taking
derivatives, $\xi_i\mapsto(\xi_i,-\xi_i\om)$),
the following $\ell$-differential form will do:
\begin{equation} \label{eq:DifFormula2}
\LA\frac{d\La}{d\la}\mu,v_1\land\ldots\land v_\ell\RA=
\LA\nu,(0,v_1)\land\ldots\land(0,v_\ell)\land
(\xi_1,-\xi_1\om)\land\ldots\land(\xi_n,-\xi_n\om)\RA
\end{equation}
And (\ref{eq:DifFormula}) follows from (\ref{eq:DifFormula1}) and
(\ref{eq:DifFormula2}).

See \cite{Federer} for further pertinent theory about measures defined
by differential forms.
} 

\item \label{it:InfMeaR2T}
Yet one must be warned that ``infinitesimal measures'' may behave
strangely: Take, for instance $G=\bR^2$ acting on the circle
$\OM=\bT=\bR/\bZ$,
endowed with Lebesgue $\mu$, by rotation using one coordinate:
$(x,y)(t)=t+x\in\bT,\;(x,y)\in\bR^2$. Consider some cases for $\La$
(to compute the infinitesimal measures, take as test-\fn\ the characteristic
\fn\ of a fundamental domain of $\bZ^2$ in $\bR^2$):
\par\medskip
For $\La=\cnt_G$, $\frac{d\La}{dx\,dy}\mu=\I\cdot\cnt$, which gives mass
$\I$ to any set except $\es$.

For $\La$ an invariant measure which, reduced to any coset of some fixed
$1$-dimensional subspace $H\st G$, is Lebesgue on the coset:
if $H=\{(0,y):y\in\bR\}$ then
$\frac{d\La}{dx\,dy}\mu=\cnt_\OM$. If $H$ is otherwise, i.e.\ ``slanted'',
then if a set in $G$ that intersects an uncountable number of $H$-cosets
has $\La$-measure $\I$, then $\frac{d\La}{dx\,dy}\mu=\I\cdot\cnt$
If $\La$ is the sum of Lebesgue measures of the
intersections with all $H$-cosets, then $\frac{d\La}{dx\,dy}\mu=\I\cdot\mu$
where $\I\cdot\mu$ gives mass $0$ to Lebesgue-null subsets of $\OM$,
and mass $\I$ to sets with positive Lebesgue mass.

For $\La$ with mass of a set equal to the integral of the numbers of points
of intersection of the set with the $H$-cosets, w.r.t.\ some Haar on $G/H$
(on smooth curves it will be given by some\NOT{(twisted)} $1$-differential
$a\,dx+b\,dy$ -- see \cite{Federer}): If $H=\{(0,y):y\in\bR\}$ then
$\frac{d\La}{dx\,dy}\mu=\I\cdot\mu$. Otherwise $\frac{d\La}{dx\,dy}\mu$
is a multiple of $\cnt_\OM$.

\item \label{it:InfMeaSO3}
Further to the previous item, let $SO(3)$ act on the sphere $S^2$ by usual
rotations, $\la$ the normalized Haar on $SO(3)$ and $\mu$ normalized
(i.e.\ probabilty) invariant area on $S^2$. This example has in common with
the previous one the property that the stabilizer of any point is subgroup
of positive dimension.

Indeed, here $\frac{d\cnt}{d\la}\mu$ is $\I\cdot\cnt$
(take as test-\fn\ the constant $1$).

\item \label{it:PoisBrw}
``Infinitesimal measures'', restricted to important subsets $E\st\OM$
of $\mu$-mass $0$ provide interesting measure spaces, generalizing the
usual $\mu|_E$ when $\mu(E)>0$. 

One may think of $\OM$ $=$ the set of discrete subsets $\om$ of $\bR^n$
with $\mu$ given by
the Poisson distribution corresponding to a Haar measure $\la$ on $\bR^n$,
i.e.\ s.t.\ the expectation, for a Borel $A\st\bR^n$, of the number
$\#(\om\cup A)$ is $\la(A)$. Let $\bR^n$ act on $\OM$ by shift, and
let $E=\{\om:0\in\om\}$. Taking as a test-\fn\ a characteristic \fn\
of some bounded domain in $R^n$, one finds:
$$\frac{d\cnt}{d\la}\mu(E)=1$$
This is a case of Palm measure (see \S\ref{SS:Palm}).

More can be said on this in view of \S\ref{SS:CtAsso}.

Another example is $\OM=\cC(\bR)/\mbox{the constants}$ with $\mu$ $=$
Brownian motion and $\bR$ acting by shift, and $E=$ those motions that
return at $1$ to where they had been at $0$. One may guess that for this
$E$ taking as $\La$ some Hausdorff measure might be interesting,
but I have not thought on that.

\item \label{it:Stoc}
Suppose $\chi:\OM\to\bR$ is some stochastic variable. One may consider
$\bE\frac{d\cnt}{d\la}\{\chi=a\}$ $=$ the ``infinitesimal expectation that
$\chi=a$''. Curiously enough, in general this {\em cannot} be integrated on
$a$ to give $\bE\chi$, as can be seen considering $\chi$ a smooth stochastic
variable on $T^2$ acted by $\bR$ as in item \ref{it:InfMea-a}.
Moreover, this depends on the action of $G$.
Some such situations will appear in \S\ref{SS:CtAsso}.

\end{enumerate}
\end{Exm}

\subsection{Links with Palm measures}
\label{SS:Palm}
The notion of Palm measure is standard in the theory of stationary random
measures on locally compact groups, in particular stationary point
processes
(see \cite{Mecke}, \cite{Delasnerie}, \cite{NeveuP}, \cite{DaleyVereJones},
\cite{Neveu}).
In our context, Palm measures are measures on $\OM$ with enhanced \fns\
(or more general invariant $0$-chains) as densities.

For us, an invariant $0$-chain, which in the discrete case is an equivalent
way to give a function, in the \cts\ case played the role of a generalized
\fn. But a $0$-chain can be viewed also as a stationary measure-valued
stochastic variable, describing a stochastic random measure on $G$. To such
an object one associates a Palm measure on $\OM$. To describe the Palm
measure from our point of view, note that if $(\OM,\mu)$ is a Borel measure
space acted in a Borel manner by, say, an abelian locally compact group $G$,
a measurable \fn\ $f:\OM\to\BAR{\bR^+}$ defines a measure $f\,d\mu$ on $\OM$
with $f$ as density. Viewing $f$ as a $0$-chain
$\Phi:\om\mapsto f(x\om)d\la(x)$ (where $\la$ is a Haar measure on $G$),
one computes the integral of a non-negative Borel \fn\ $g$ w.r.t.\ $f\,d\mu$
by multiplying the $0$-chain $\Phi$ by $g(x\om)$ and taking the
``expectation'' of the resulting $0$-chain.
This means that if $h:G\to\BAR{\bR^+}$ is Borel, one has:
$$\LA f\,d\mu,g\RA\int h\,d\la=
\int\,d\mu(\om)\LP \int h(x) g(x\om)f(x\om)\,d\la(x)\RP$$
and this can be generalized to any Borel-measurable invariant $0$-chain
$\Phi$ to define the Borel measure $\Phi\,d\mu$ on $\OM$ with density $\Phi$,
called the Palm measure:
$$\LA \Phi\,d\mu,g\RA\int h\,d\la=
\int\,d\mu(\om)\LP h(x) g(x\om)\,\LP d\Phi(\om)\RP(x)\RP.$$

In particular, if $\La$ is an invariant measure on $G$ (not necessarily
$\sigma$-finite) and $E\st\OM$ is Borel, then the enhanced \fn\ 
$1_E\,\frac{d\La}{d\la}$ has a Palm measure $1_E\,\frac{d\La}{d\la}\,d\mu$
which is just the restriction of the ``infinitesimal measure'' 
$\frac{d\La}{d\la}\,d\mu$ to $E$. If $\cO_\om E$ is discrete for a.a.\ $\om$,
then this random discrete set defines a stationary point process with this
Palm measure. This is the situation in Ex.\ \ref{Exm:InfMea} \ref{it:PoisBrw}.

\subsection{Hypergraphs, Weighted Hypergraphs and the Continuous HG (CHG)
Theorem}
\label{SS:CHG}
We restrict ourselves to {\em unimodular} acting group $G$.

As mentioned above, in the \cts\ case {\EM hypergraphs} do not automatically
define $m$-chains. One may also consider ($\om$-dependent) nonnegative \fns\
on the set of $m$-simplices (``matrices''), now {\em not} the same as
$m$-chains (they do not have vertices). Such \fns\ will be referred to as
{\EM weighted hypergraphs} $F(x_0,\ldots,x_m;\om)\:\:x_i\in G,\om\in\OM$.
They will be assumed {\EM (right)-invariant}, in the sense that:
\begin{equation} \label{eq:WHInvr}
  F(x_0y,\ldots,x_my;\om)=F(x_0,\ldots,x_m;y\om)\quad x_i,y\in G, \om\in\OM
\end{equation} 

Such a weighted hypergraph may be converted into an $m$-chain if an
{\em $m$-simplex $(\La_0,\ldots,\La_m)$ of right-invariant measures
in $G$}
is given. Assuming enough measurability and ``Fubini'',
just multiply $F$, for each $\om$, by the {\em product measure} on $G^{m+1}$.

CVE leads to the following:

\begin{Main}
{\EM The Continuous Hypergraph (CHG) Theorem.}
Let a 2nd-countable {\em unimodular} locally compact group $G$ act
in a Borel and measure-preserving manner on a standard $\sigma$-finite
measure space $(\OM,\cB,\mu)$.
Let $\la$ be a Haar measure on $G$.

Let us be given a (right-)invariant weighted hypergraph.
This is a nonnegative function
$F(x_0,\ldots,x_m;\om)$
on the set $G^{m+1}$ of $m$-simplices, depending on $\om$.
(Right-)invariance means that (\ref{eq:WHInvr}) is satisfied.

Let us be given also an $m$-simplex $(\La_0,\ldots,\La_m)$ of
right-invariant complete measures in $G$ for which Borel sets are
measurable (they need not be a Radon measures).

Assume:
\begin{itemize}
\item
For a.a.\ $\om$, $F$ on $G^{m+1}$ is measurable w.r.t.\
to the product measure, and is $0$ outside a countable union of products in
$G^{m+1}$ of\NOT{Borel} sets with finite corresponding $\La_i$-measures,
so Fubini holds for the product of measures.

\item
When $F$ is multiplied by that product measure, the invariant
$m$-chain obtained is measurable.  
\end{itemize}

Then for all $i=0,\ldots,m$, the integral of the weighted hypergraph
$F$, over the set of simplices with $0$ as the $i$-th vertex, w.r.t.\ the
product of the measures $\La_{i'},\:\:i'\ne i$, this integral being a
function of $\om$, is $\frac{d\La_i}{d\la}\mu$-measurable.
If, for some $i=0,\ldots,m$, it has, when enhanced by $\frac{d\La_i}{d\la}$, 
a finite expectation (that is, $\mu$-integral),
then the same holds for any $i$, with the same
expectation (and since the only other possibility is all these expectations
being $\I$, the word ``finite'' may be deleted).
\end{Main}
The above $\om$-depending integral, or its above enhancement by
$\frac{d\La_i}{d\la}$, will be called:
{\EM the $i$'th vertex of the weighted
hypergraph w.r.t.\ the simplex of measures}.

\begin{Prf}
Note that {\em we do not have Fubini} for interchanging integrations w.r.t.\
$\La_i$ and $\mu$, but the fact that we have Fubini for the $\La_i$'s
defines the product measure unequivocally.

We formulate the proof for $m=1$, applying CVE to the $1$-chain. We have to
consider its two vertices, i.e.\ projections on $G$. These are:
\BER{lcl}
\mbox{source }&=&
 \om\mapsto\LQ\int_G F(x_0,x_1;\om)\,d\La_1(x_1)\RQ\,d\La_0(x_0)\\
\mbox{target }&=&
 \om\mapsto\LQ\int_G F(x_0,x_1;\om)\,d\La_0(x_0)\RQ\,d\La_1(x_1)
\EER
(Here Fubini for $\La_0$ and $\La_1$ was used implicitly.)

The source is a $\La_0$-enhancement of the function
$$\om\mapsto\int_G F(0,x_1;\om)\,d\La_1(x_1).$$
that follows from the fact that substituting $x_0\om$ for $\om$ in this
function gives
$$\int_G F(0,x_1;x_0\om)\,d\La_1(x_1)=\int_G F(x_0,x_1x_0;\om)\,d\La_1(x_1)=
\int_G F(x_0,x_1;\om)\,d\La_1(x_1)$$
(here we need actual right-invariance of $\La_1$, hence unimodularity),
similarly for the target. These are the enhanced \fns\ mentioned
in the theorem. If one of these has a finite expectation, then the $0$-chain
vertex has a Radon expectation and applying CVE we are done.
\end{Prf}

\beware
It might be helpful to try to
find conditions easier to check than measurability in
the CHG Thm.\beware

\subsection{The Question of Measurability}
\label{SS:MEAS}

To apply CHG,\NOT{(\S\ref{SS:CHG})}
one needs to know that the $m$-chain is measurable.

We refer to \S\ref{SS:PRD}.

In view of following sections (such as Section \ref{S:EDCT})
we wish not to refer to a particular (probability) measure $\mu$ in $\OM$.
This gives us the option to require Borel-measurability. Another significant
notion is sets or functions being {\EM universally measurable} -- measurable
w.r.t.\ any completion of a finite (equivalently, $\sigma$-finite)
Borel measure. It is known (see \cite{KechrisD} Thm.\ (21.10) or
\cite{Kuratowski} \S11 VII.\ or \cite{BourbakiTOP} Ch.\ IX \S6)%
\NOT{or\NOT{the proof of} Prop.\ \ref{Prop:UnivMeas} for $\La=\cnt$}
that any Suslin set in a standard space, i.e.\ any image of some (Borel
subset of a) standard space by a Borel mapping, is universally measurable.
(If the Borel mapping is countable-to-one, i.e.\ the preimage of every point
is at most countable, then the image is Borel 
-- see \cite{Lusin} Ch.\ III,IV).

Note that if a \fn\ $f:X\to Y$ between standard Borel spaces is given,
and a Lusin topology is given in each of the spaces, then $f$ is universally
measurable iff for any finite measure in $X$ s.t.\ Borel sets are
measurable, there is a compact $K\st X$ with $K^c$ of mass as small as we
please, s.t.\ $f|_K$ is \cts.

Consequently, the composition of universally measurable mappings between
standard Borel spaces is universally measurable.

\NOT{
Note that the composition $f_2\circ f_1$ of two universally measurable
mappings between Borel spaces is universally measurable. Indeed, the
preimage of a Borel set $E$ by $f_2$ is measurable w.r.t.\ any image by
$f_1$ of a completion of a ($\sigma$-)finite Borel measure, that image
being itself a completion of a Borel measure. Thus the preimage $f_2^{-1}(E)$
differs from a borel set by a null set, hence its preimage
$f_1^{-1}(f_2^{-1}(E))$ differs from a measurable set by a null set, hence
is measurable.
}

In many of our applications the measurability of the $m$-chain in an
application of CHG may be assured by constructing the (weighted) hypergraph
in two stages:

First, one corresponds to each $\om$ a
{\em closed subset} $F(\om)\st G$. In many cases this will be
the closure $\BAR{\cO_\om E}$ for some Borel $E\st\OM$. Another
alternative is the {\em essential closure} of $\cO_\om E$ -- the set of
points in $G$ no \nbd\ of which intersects $\cO_\om E$ in a Haar-null set.
In \S\ref{SS:HE} other alternatives are considered.

Second, the (weighted) hypergraph dependence on $\om$ is obtained by a rule
corresponding to each closed subset of $G$ a (weighted) hypergraph, with no
mention of $\OM$.

Now in the set of closed subsets of $G$ we always take the Effros Borel
structure (see \cite{KechrisD}) -- just identify each closed set $F$ with
the set of members $u$ of a fixed countable open base to the topology that
satisfy $u\cap F=\es$, that set being a member of $2^\bN$. This is a
standard Borel space (This Borel structure can be given in many other ways,
e.g.\ it is the Borel structure of the topological space $2^G$
(with members the {\em closed} subsets of $G$) with
subbasis consisting of sets of all closed set contained,
or all closed sets intersecting, an open set -- see \cite{Kuratowski}.
This topology is given by the Hausdorff metric for a restriction of
any metric of the Alexandrov (one-point) compactification of $G$.)

In the cases mentioned above one is sure that $\om\mapsto F(\om)$ is 
measurable. Indeed, for $F(\om)=\BAR{\cO_\om E}$, one has:

\begin{Prop} \label{Prop:EGDelta}
Let a 2nd countable locally compact group $G$ act in a Borel manner on a
standard Borel space $\OM$. Let $E\st\OM$ be Borel. Then:

\begin{itemize}
\item[(i)]
The mapping: $\om\mapsto\BAR{\cO_\om E}$ from $\OM$ to the Borel space of
closed subsets of $G$ is universally measurable.

\item[(ii)]
(see \S\ref{SS:CP})
Suppose $\cO_\om E$ is {\em countable} for every $\om\in\OM$. Then
there exists an embedding of the $G$-Borel space $\OM$ in a $G$-metrizable
compact space and a $E'\st\OM$ which is $G_\delta$ in the relative topology
s.t.\ $\forall\:\om\in\OM$, $\BAR{\cO_\om E}=\cO_\om E'$.

\item[(iii)]
Suppose $\cO_\om E$ is {\em countable} for every $\om\in\OM$. Then
the mapping: $\om\mapsto\BAR{\cO_\om E}$ from $\OM$ to the Borel space of
closed subsets of $G$ is Borel.
\end{itemize}
\end{Prop}

\begin{Prf}

\begin{itemize}
\item[(i)]
By the way the Borel structure in the set of closed sets is defined, it
suffices that for a fixed open $u\st G$, the set of $\om$ s.t.\ $u$
intersects $\BAR{\cO_\om E}$, that is, $u$ intersects $\cO_\om E$, is
universally measurable. But this set is Suslin, being a projection of a
Borel set $\st\OM\times G$.

\item[(ii)]
Let $\la$ be a right Haar measure in $G$.
Choose a sequence $h_n$ of $\I\times$the characteristic \fns\ of a decreasing
sequence $U_n$ of open \nbd{s} of $0$ in $G$ which forms a basis to the
\nbd{s} at $0$.

The \fns\ $(\om,x)\mapsto h_n(x)1_E(x\om)$ are Borel, for each $\om$
having at most countable number of $x$ where they do not vanish. By
\cite{Lusin} Ch.\ III,IV (cf.\ Rmk.\ \ref{Rmk:CntMeasFub} below) their
sums over $x$
$$f_n(\om):=\sum_{x\in G}h_n(x)1_E(x\om)$$
are Borel \fns\ of $\om$.

Now apply \S\ref{SS:CP}.
Since each $h_n$ is the increasing limit of non-negative $L^1$ -\fns,
there is an embedding of $\OM$ in a $G$-metrizable compact space s.t.\ all
the convolutions:
$$\tilde{f}_n(\om)=\int_G f_n(x\om)h_n(x)\,d\la(x)$$
are l.s.c.\ (lower semi-\cts) \fns. It is easy to see that these
$\tilde{f}_n$ are the same as $f_n$, but with $U_n^2$ instead of $U_n$.

Now, since $\tilde{f}_n$ are $\{0,\I\}$-valued l.s.c.\ \fns, the sets
$\{\tilde{f}_n=\I\}$
are open and their intersection $E'$ is $G_\delta$.
It is easy to see that this $E'$ will do. 

\item[(iii)]
Metrize $G$ by metrizing the Alexandrov (one-point) compactification of $G$.
Then the Borel structure in the closed sets $\st G$ is obtained from the
Hausdorff metric.
Choose a dense countable set $D\st G$. The mapping which maps each closed
$F\st G$ to the sequence $(\dist(x,F))_{x\in D}$ is 1-1 Borel, hence a 
Borel isomorphism with the image, thus it is enough to prove
$\om\mapsto\dist(x,\BAR{\cO_\om E})$ is Borel for each fixed $x\in D$.
But one may write the $E'$ in (ii) as
$E'=\cap_{n\ge1}W_n=\cap_{n\ge1}\BAR{W_n}$,
with $W_n$ open and $\BAR{W_{n+1}}\st W_n$ in some
relative topology from a $G$-metrizable compact space in which $\OM$ is
embedded. Then
$$\dist(x,\BAR{\cO_\om E})=\sup_n\inf_{y\in D,y\om\in W_n}\dist(x,y).$$
\end{itemize}
\end{Prf}

For the ``essential closure'' Borel-measurability follows from Fubini.

Thus measurability in $\om$ will follow if we insure that
the way the (weighted) hypergraph is constructed from the closed set
causes the $m$-chain made from the (weighted) hypergraph and the simplex
of measures to depend measurably on the closed set.

In most of the examples in the sequel, the weighted hypergraph is defined,
depending on the closed set, by Borel operations in the points of $G$ and
in closed sets (the given closed set, and also, say, a given closed
relation, such as a partial ordering). One uses the fact that 
{\em in a 2nd-countable locally compact space}
finite Boolean operations in closed sets
(and also countable intersection) are Borel (a not too hard exercise),
and so are $x\mapsto\{x\}$ for $x\in G$ and
$(x,F)\mapsto\{y\in G:(x,y)\in F\}$ for $x\in G$, $F\st G^2$ closed.
(Open sets or $F_\sigma$ sets can also be encoded in a Borel manner using
closed sets: encode an open set by its complement, and instead of a variable
$F_\sigma$ set take a variable sequence of closed sets. Yet, one should be
cautious: the relation $\cup_{n\ge1}F_n=F_0$ is co-Suslin, but {\em
not Borel}%
\footnote{%
Indeed, the relation $\{x_n\}=F$ for $F$ compact is not Borel in
$((x_n)_n,F)$. Otherwise its countable-to-one image: the family of countable
compact sets, would be Borel. But for any Borel family $\cF$ of countable
compact sets, the (countable ordinal) order of the first vanishing derivative
is bounded. Indeed (see \cite{Lusin} Ch.\ IV) there is a polish topology
in the set of pairs $\{(x,F):x\in F\in\cF\}$ finer than the product of
the topologies in $G\times2^G$. Moreover, we may assume this Polish topology
has a countable basis $\cU$ composed of clopens. Suppose the order $\gamma$
of the first vanishing derivatives was unbounded. Then one can find disjoint
$U_0,U_1\in\cU$ of diameter $<2^{-1}$ s.t.\
$\forall\gamma\exists F\in\cF,\:F^{(\gamma)}\cap\{x:(x,F)\in U_j\}\neq\es$
(not having that implies that for $\gamma$ big enough $F^{(\gamma)}$ is a
singleton $\forall F\in\cF$).
Then one finds disjoint $U_{ij}\in\cU,\:i,j\in\{0,1\}$ of
diameter $<2^{-2}$ and s.t.\ $U_{ij}\st U_i$ and
$\forall\gamma\exists F\in\cF,\:F^{(\gamma)}\cap\{x:(x,F)\in U_{ij}\}\neq\es$
and so on. The set $K=\cup U_i\cap\cup U_{ij}\cap\cup U_{ijk}\cap\ldots$
is a Cantor set in the Polish space, and all its elements must be pairs
with the same (uncountable) $F_\I$, otherwise, since the Polish topology
is finer than that of $G\times2^G$, there will be $U_{ij\ldots}$'s without
member pairs with common $F$.%
})

Next one applies a simplex of invariant measures, and one has to
know that their values on Borel sets, and the integral of Borel
\fns\ on them, depend measurably on the sets and \fns, in some sense.

A family of $E(\al)$ of Borel subsets of a standard Borel space $X$,
depending on a parameter $\al$ varying in a standard Borel space $\cA$,
will be called a {\EM Borel family} if the set
$$\{(x,\al):\om\in E(\al)\}\st X\times\cA$$
is Borel.

Similarly, a family $f(x;\al)$ of Borel \fns\ on a standard Borel space $X$,
depending on a parameter $\al$ varying in a standard Borel space $\cA$,
will be called a {\EM Borel family} if $(x,\al)\mapsto f(x;\al)$ is Borel.

\begin{Rmk} \label{Rmk:CntMeasFub}
The counting measure has the following property: for a Borel family
of sets (\fns), the counting measure of $E(\al)$
(the counting integral of $x\mapsto f(x;\al)$)
is {\em Borel measurable} in the parameter if
all the sets in the family are at most countable
(all the \fns\ in the family differ from zero only on a countable set
(which depends on the parameter)).

Indeed, By \cite{Lusin} Ch.\ III,IV if $X$ $\cA$ are standard
and $f:X\times\cA\to\BAR{\bR^+}$ is Borel s.t.\
for each $\al\in\cA$, $\om\in\OM\mapsto f(\om,\al)$ is different from $0$
only on a countable set, then one can write
$$f(x,\al)=\sum_{n\ge1}g_n(\al)\delta_{x,\psi_n(\al)}$$
for Borel $g_n$, $\psi_n$ with $g_n\ge0$. Thus
$$\int_X f(x,\al)\,d\cnt_X(x)=\sum_{n\ge1}g_n(\al)$$
is Borel in $\al$.
\end{Rmk}

\NOT{
\begin{Def}
\label{Def:UnivMeas}
Given a complete measure $\La$ on a standard Borel space $X$ with all Borel
sets measurable. We say that $\La$ is
{\EM universally measurable on Borel sets}, 
if for every Borel family $E(\al)$ of Borel subsets of a standard Borel
space $X$, the mapping $\al\mapsto\La(E(\al))$ is universally measurable.
equivalently, for any Borel family of $\BAR{\bR^+}$-valued Borel \fns\ on
$X$, $\int_X f(x;\al)\,d\La(x)$ is universally measurable in $\al$.
\end{Def}

By Fubini any $\sigma$-finite completion of a Borel measure is universally
measurable, the measure or integral of a Borel family being even 
Borel-measurable.
}

\begin{Def}
\label{Def:MeasFub}
Given a complete measure $\La$ on a standard Borel space $X$ with all Borel
sets measurable. We say that $\La$ has the {\EM Borel-} resp.\
{\EM universal- measurability Fubini property}
if for every Borel family $E(\al)$ of Borel
subsets of $X$ s.t.\ for each $\al\in\cA$, $E(\al)$ is $\La$-$\sigma$-finite,
the mapping $\al\mapsto\La(E(\al))$ is Borel-, resp.\ universally- measurable.
Equivalently, if for any Borel family of $\BAR{\bR^+}$-valued Borel \fns\ on
$X$ s.t.\ for each $\al\in\cA$, $x\in X\mapsto f(x;\al)$ is different
from $0$ only on a $\La$-$\sigma$-finite set (which depends on $\al$), the
integral $\int_X f(x;\al)\,d\La(x)$ is Borel- resp.\ universally-
measurable in $\al$.
\end{Def}

By Fubini any $\sigma$-finite completion of a Borel measure has the
Borel measurability Fubini property.

By Rmk.\ \ref{Rmk:CntMeasFub}, the counting measure has the Borel
measurability Fubini property.

\begin{Exm}
An example for a Borel measure which does not have the universal
measurability Fubini property: Let $f:[0,1]\to]0,\I[$ be a \fn\ which is not
universally measurable and let $\la$ be the Lebesgue measure on $[0,1]$.
Define a Borel measure on the square $[0,1]^2$ as follows:
a Borel set $E\st[0,1]^2$ whose projection on the first coordinate is
countable has measure $\sum_xf(x)\la\{y:(x,y)\in E\}$. If the projection is
uncountable the measure is $\I$. Then the mass of the member of the Borel
family with parameter $x$ -- $\{x\}\times[0,1]$, is $f(x)$, which is not
universally measurable.
\end{Exm}

As for other measures $\La$ on a 2nd-countable locally compact $X$
(such as those mentioned in \cite{Federer} \S2.10), one can say the following.

It follows from the fact that for a Borel family of Borel sets, the property
($n\in\bZ^+\cup\{\I\}$ fixed):
``the cardinality of the Borel set $\ge n$'' defines a Suslin set of the
parameters (it is the set of parameters s.t.\ the set composed of
``sequences of $n$ distinct points in the set'' is not empty),
and from Fubini, that the universal measurability Fubini 
property holds for measures whose value on a set $E$ is defined by taking
the counting measure of the intersection of $E$ with a variable closed
set $F$ and integrating on a fixed (say, invariant) $\sigma$-finite
measure on $F$.

This includes\NOT{besides the counting measure and $\sigma$-finite measures,}
$k$-dimensional measures in $\bR^n$, $k<n$ (intersect with $(n-k)$-dimensional
affine subspaces and take an invariant Radon measure on these) -- see
\cite{Federer}.

Many measures have the following property:
$K\mapsto\La(K)$ is Borel\NOT{(resp.\ universally measurable)}
on the set of compact $K\st X$
(one easily sees that this latter set is Borel in the space of closed sets).
We shall say then that $\La$ is {\EM Borel (measurable) on compacta}.%
\NOT{ resp.\ {\EM universally measurable on compacta}} 

Indeed, measures whose value on a set $E$ is defined by taking the counting
measure of the intersection of $E$ with a variable closed set $F$ and
integrating on a fixed (say, invariant) $\sigma$-finite measure on $F$
are Borel on compacta -- This follows from Borelness of the operations of
taking intersection of closed sets and taking the number of elements of a
closed set (that number being finite or $\I$).

Also, one has Borel measurability on compacta for measures, such as Hausdorff
measures, having the following property:
there exist a family $(U_{nm})_{n,m\ge1}$ of open sets in $X$
and numbers $\eta_{nm}\in\BAR{\bR^+}$ s.t.\ for compact $K\st X$:
$$\La(K)=\sup_n\inf_{\{m:K\st U_{nm}\}}\eta_{nm}$$
(For $r$-dimensional Hausdorff measure in $\bR^N$,
let $U_{nm}$ for fixed $n$ enumerate the finite unions of the members of
diameter $<1/n$ in a countable family of open convex sets, such as the family
of all finite intersections of rational half-spaces, that have a member
between any compact convex set and any of its open \nbd{s}, and take
as $\eta_{nm}$ the infimum of $\sum_{W\in\cF}(\diam W)^r$ over finite
coverings $\cF$ of $U_{nm}$, multiplied, if needed, by a normalizing
constant.)

\NOT{
For such measures, the universal measurability Fubini property is settled by

\begin{Prop} \label{Prop:UnivMeas}
Let $\La$ be a completion of a Borel measure on a 2nd-countable
locally compact space $X$. If $\La$ is Borel on compacta then it 
has the universal measurability Fubini property.%
\NOT{(Def.\ \ref{Def:MeasFub})}
\end{Prop}
}

Now, if a completion of a Borel measure $\La$ 
is Borel on compacta, then it has the universal
measurability Fubini property, moreover for a Borel family
$\al\mapsto E(\al):=\{x\in X: (x,\al)\in E\}$, $E\st X\times\cA$ Borel,
$\cA$ standard, s.t.\ $\forall\al$ $E(\al)$ is $\La$-$\sigma$-finite,
the set $\{\al\in\cA:\La(E(\al))>a\}$ is {\em Suslin}.

Indeed, endow $\cA$ with a Polish topology. Then $E$ is a Lusin space
hence $\exists$ a Polish space $Y$ and a 1-1 \cts\ onto $f:Y\to E$.
the compositions of $f$ with the projections on $X$ and on $\cA$ are
\cts. Thus the $\{\al\}\times E(\al)$ correspond to {\em closed} subsets
$F(\al)$ of $Y$, on which $\La$ induce $\sigma$-finite Borel measures.
Since $\sigma$-finite Borel measures on Polish spaces are regular
(for Borel sets), the masses of the $F(\al)$ are the supremum of the
masses of compacta $\st F(\al)$.

Therefore $\La(E(\al))>a$ iff $\exists$ a compact $K\st Y$, mapped
onto $\{\al\}$ by the projection on $\cA$ and mapped by the projection on
$X$ onto a (compact) set ($\st X$) with $\La>a$. Now, the compact sets
in a Polish space form, with Hausdorff metric, a Polish space, and mapping
compact sets to their image via a \cts\ \fn\ is \cts. Thus it follows from
$\La$ being Borel on compacta in $X$ that $\{\al:\La(E(\al))>a\}$ is
Suslin.

Similar considerations show that the {\em product} of two measures $\La_1$
and $\La_2$ on two 2nd-countable locally compact $X_1$ and $X_2$,
both completions of Borel measures,
has the universal measurability Fubini property,
if it is assumed that both measures are Borel on compacta, and the product
is defined so that only sets contained in a $\sigma$-finite Borel
``rectangle'' in $X_1\times X_2$ can have finite product mass.

Indeed, if both measures are Borel on compacta, then the product measure
is Borel on compact sets $K\st X_1\times X_2$ with projections on $X_1$ and
$X_2$ having finite mass
(The $\La_1\otimes\La_2$-mass of such $K$ is equal to an integral
$\int_{K_1(K)}f(K,x_1)\,d\La_1(x_1)$
with a Borel \fn\ $f(K,x_1)$ and with $K_1(K)$ a compact in $X_1$ of finite
$\La_1$-mass which depends in a Borel manner on $K$.
To prove such an integral is Borel in $K$, one may assume the Borel \fn\ $f$
is the characteristic \fn\ of a Borel set $F$ of $(K,x_1)$'s,
and since the collection of $F$'s satisfying what we want is {\em monotone}
(\cite{HalmosMEA} \S6), i.e.\ stable w.r.t.\ unions and intersections of
monotone sequences of Borel sets, and contains the Boolean algebra
consisting of the finite disjoint unions of ``rectangles'' with one side a
Borel set of $K$'s and the other side a difference of closed sets in $X_1$
-- here we use the fact that intersecting a fixed closed set in $X_1$
with $K_1(K)$ is a Borel operation in $K_1$ -- it contains all Borel sets.)
Then to prove the universal measurability Fubini property of the product
measure we proceed as above, restricting the compact sets in the Polish space
$Y$ to those having projections on $X_1$ and $X_2$ with finite mass.

\begin{Rmk}
Note that for a Borel action of a 2nd-countable locally compact
group $G$ on a standard Borel space $\OM$, $\spn E$ is, for any
Borel $E\st\OM$, a projection of a Borel set in a product, hence Suslin,
thus universally measurable.
(One cannot say more in general: for a projection $\mbox{pr}_1(E)$ of a Borel
set $E\st\bR^2$ on $\bR^1$ (which need not be Borel),
we have $\mbox{pr}_1(E)\times\bR=\spn E$ for
the action of $\bR$ on $\bR^2$ via shift on one coordinate.)

If $E$ has the property that $\cO_\om E$ is {\em countable} for all
$\om\in\OM$, then $\spn E$ is Borel since the image of a countable-to-one
Borel mapping among standard spaces is Borel (\cite{Lusin} Ch.\ III,IV --
see the beginning of this \S).

Even when $\spn E$ is universally measurable but not Borel, one can find,
for every invariant Borel probability measure $\mu$ in $\OM$, an invariant
Borel set $E'$, differing from $\spn E$ in a ($\mu$-)null set and
containing it, resp.\ is contained in it (in which case $E'=\spn(E\cap E')$,
with $E\cap E'$ differing from $E$ in a $\frac{d\La}{d\la}\mu$-null set
for any $\La$) as is shown by the following

\begin{Prop} \label{Prop:BorelAlmInv}
Let a 2nd-countable locally compact group $G$ act in a measure-preserving
Borel manner on a $\sigma$-finite standard measure space $(\OM,\cB,\mu)$.

\begin{itemize}
\item[(i)]
Let $f:\OM\to\BAR{\bR^+}$ be Borel and almost-invariant in the sense that
for each fixed $x\in G$ $f(x\om)=f(\om)$ a.e. Then $\exists$ an invariant
Borel \fn\ $\tilde{f}:\OM\to\BAR{\bR^+}$ s.t.\ $\tilde{f}(\om)=f(\om)$
a.e.\ and s.t.\ if $\om$ is s.t.\ $\cO_\om f$ is constant a.e.\ (w.r.t.\
Haar measure on $G$) then $\tilde{f}(\om)$ is equal to this constant.
Moreover, $\tilde{f}$ can be chosen s.t.\ there is an embedding of $\OM$
as an invariant Borel subset of a $G$-metrizable compact space s.t.\
$\tilde{f}$ is l.s.c.\ (lower semi-\cts) for the relative topology on $\OM$.

\item[(ii)]
Let $E\st\OM$ be Borel and almost-invariant (i.e.\ $\forall$ fixed $x\in G$
$E\triangle xE$ is null). Then $\exists$ an invariant Borel set 
$\tilde{E}$ which differs from $E$ only on a null set, and s.t.\ if
$\om$ is s.t.\ $\cO_\om E$ is null (resp. conull) w.r.t.\ Haar measure on
$G$, then $\om\notin\tilde{E}$ (resp. $\om\in\tilde{E}$).

\item[(iii)]
Every measurable invariant set $E$ differs by a null set from some Borel
invariant set $\tilde{E}$ contained in it (resp. which contains it).
\end{itemize}
\end{Prop}

\begin{Prf}

\begin{itemize}
\item[(i)]
We refer to \S\ref{SS:CP}.

Let $\la$ be a Haar measure on $G$ and let $h$ be some non-negative $L^1$
\fn\ on $G$. Let $f'$ be the Borel $\fn$ on $\OM$:
$$f'(\om):=\int_Gf(x\om)h(x)\,d\la(x).$$
The \fns\ $(\om,x)\mapsto f(x\om)$ and $(\om,x)\mapsto f(\om)$ on
$\OM\times G$ are Borel, and for every fixed $x$ are $\mu$-a.e.\ equal.
Therefore they are $\mu\otimes\la$- a.e.\ equal (here we use the fact that
$\mu$ is $\sigma$-finite!), hence for $\mu$-a.a.\ $\om$ $\cO_\om f$ is
equal $\la$-a.e. to the constant $f(\om)$, and for these $\om$'s
$\cO_\om f'$ is the constant $f(\om)$.

By \S\ref{SS:CP} there is an embedding of $\OM$ as an invariant Borel subset
of a $G$-metrizable compact space $K$ s.t.\ $f'$ can be extended to a
l.s.c.\ \fn\ on $K$. Let the \fn\ $\tilde{f}$ on $\OM$ be:
$$\tilde{f}(\om):=\sup_{x\in G}f'(x\om)$$
For each $\om$ s.t.\ $\cO_\om f$ is $\la$-a.e.\ the constant $f(\om)$, we
have $\tilde{f}(\om)=f(\om)$. Also $\tilde{f}$ is invariant and l.s.c.
Thus it satisfies our requirements.

\item[(ii)]
follows from (i).

\item[(iii)]
Choose a Borel set $E'\st E$ (resp. $E'\supset E$) which differs from $E$
by a null set. For every $x\in G$ $xE=E$ and $xE'$ and $E'$ differ from it
by null sets. Therefore $E'$ is almost-invariant. An invariant Borel
set $\tilde{E'}$ as in (ii) will do.
\end{itemize}
\end{Prf} 
\end{Rmk}

\subsection{An Assortment of Kac-like Theorems, Continuous Case}
\label{SS:CtAsso}
CHG gives us analogs of the Kac-like theorems of Section \ref{S:DS},
{\em mutatis mutandis}.

Recall that, when we deal with a (weighted) {\em graph} and a
$1$-simplex of invariant measures, {\em the source is the target measure
of the set of targets of arrows with source $0$}, analogously for the
target.

\subsubsection{A Continuous Kac Formula}

Let us start with the promised (see Example \ref{Exm:InfMeaRT}, cf.\ 
\cite{Birkhoff})

\begin{Thm} \label{Thm:CtsKac}
{\EM A Continuous Kac Theorem}
Let $\bR$ act measure-preservingly and in a Borel manner
on a standard probability measure space $(\OM,\mu)$. Let $E\st\OM$ be Borel.

Define:
\BER{lcl}
\Ex E&:=&\{\om\in\OM:
\forall\eps>0\:\:\om\mbox{ visits }E\mbox{ in the }]0,\eps[\mbox{-past}\}\\
\In E&:=&\{\om\in\OM:
\forall\eps>0\:\:\om\mbox{ visits }E\mbox{ in the }]0,\eps[\mbox{-future}\}\\
\EER
(with obvious meaning for ``visiting'' etc.)

Note that $E\cup\Ex E\cup\In E=\{\om\in\OM:0\in\BAR{\cO_\om E}\}$.

Define the {\EM return time} $\rho(\om)=\rho_E(\om):\OM\to[0,\I]$ by:
$$
\rho_E(\om):=\left\{\begin{array}{ll}
\inf\{t>0;T^t\om\in E\} &\mbox{if }\om\in E\cup\Ex E\\
0 &\mbox{otherwise}
\end{array}\right.
$$

Note that on $(E\cup\Ex E)\sm\In E$ we have $\rho>0$.

Then $\In E$, $\Ex E$, $E\cup\Ex E$ and $E\cup\In E$ differ only in null sets
and if $\tilde{E}$ is any set differing from them by a null set (such as
$E\cup\In E\cup\Ex E$) then we have:
$$\bE\LP\rho\frac{d\cnt}{dt}\RP = \mu((\spn E)\sm \tilde{E})$$
in other words,
\begin{equation} \label{eq:CtsKac}
  \bE\LP\rho\frac{d\cnt}{dt}+1_{\tilde{E}}\RP = \mu(\spn E)
\end{equation} 

Note that if a.e.\ $\cO_\om E$ is {\em closed}, one may take $\tilde{E}=E$;
if a.e.\ $\cO_\om E$ is {\em discrete}, one may take $\tilde{E}=\es$.
\end{Thm}
\begin{Prf}
(note the analogy with our proof of the usual Kac in \S\ref{SS:EqiChn})

First, we have Poincar\'e recurrence in the following sense: for a.a.\ 
$\om\in\spn E$ $E$ is visited in any $[t_0,\I[$-past or future.
Indeed, the sets
$$A_n=\{\om:[n,n+1[
 \mbox{ is the last such interval intersecting }\cO_\om E\}$$
are in the Boolean algebra generated by Suslin sets, hence are measurable,
are disjoint and have the same measure.

Now apply CHG to the $\om$-dependent graph consisting of the arrows
$(t_0,t_1)$ s.t.\ $t_0>t_1$ and $[t_1,t_0]\cap\BAR{\cO_\om E}=\{t_1\}$,
and for the $1$-simplex of invariant measures on $\bR$ $(dt,\cnt_\bR)$.
(Check that the conditions for applying CHG are satisfied.)
\par\medskip

The {\em source} equals a.e.\ $1_{\spn E\sm(E\cup\Ex E\cup\In E)}$.
\NOT{
(recall that the source is the target measure
of the set of targets of arrows with source $0$)
}

The {\em target} equals $\rho_E$.
\NOT{
(the source measure of the set of sources of arrows with target $0$)
}
\par\medskip
which, by CHG, proves the asserted formula for $\tilde{E}=E\cup\Ex E\cup\In E$.
This shows that $\rho_E$ has finite $\frac{d\cnt}{dt}\mu$-integral, hence,
by Prop. \ref{Prop:InfMea0I}, has zero $\mu$-integral,
hence $\mu((E\cup\Ex E)\sm\In E)=0$. Applying this for the inverse action
one deduces $\mu((E\cup\In E)\sm\Ex E)=0$, and we are done.
\end{Prf}

The second term in the left-hand side of (\ref{eq:CtsKac}) may be
interpreted as the result of points $\om$ in $\In E$
having ``infinitesimal return time''
(indeed, ``having a.e.\ return time $\frac{dt}{d\cnt}$'' !?).

\begin{Rmk}
Here one may appreciate the need for the
$\prod\La_i$-$\sigma$-finiteness assumption in the formulation of CHG.
Let us change the graph in the proof of \ref{Thm:CtsKac} by
requiring $t_0\ge t_1$ instead of $t_0>t_1$, thus adding
all the arrows $(t,t),\:t\in\BAR{\cO_\om E}$ to the graph. This does not
change the target, but the source is now a.e.\ $1_{\spn E}$ and the
formula that might be obtained is incorrect. Indeed, the $\sigma$-finiteness
is not satisfied (for $\cnt_\bR$).
\end{Rmk}

\begin{Rmk}
Let $T>0$, and let us integrate the test-\fn\ (nonegative Borel) $1_{[0,T[}$
on the $0$-chain corresponding to the enhanced \fn\ under the expectation in
(\ref{eq:CtsKac}), for $\tilde{E}=E\cup\In E\cup\Ex E$, thus
$\cO_\om 1_{\tilde{E}}=1_{\BAR{\cO_\om E}}$. The integral on
$1_{\BAR{\cO_\om E}}\,dt$ is complemented by the integral on
$\cO_\om\rho\,d\cnt_\bR$ and we have the following fact, equivalent, in fact,
to (\ref{eq:CtsKac}):
\begin{equation} \label{eq:CtsKacI}
   \forall T>0\quad
   {1\over T}\int_{\OM}\LQ\inf(\cO_\om E\cap\,]T,\I[\,)-
   \inf(\cO_\om E\cap\,]0,\I[\,)\RQ\,d\mu(\om)=1
\end{equation}
(Taking, e.g., $[0,\I]$ instead of $]0,\I]$ will change the integrand only
on $E\sm\In E$, a null set by Thm.\ \ref{Thm:CtsKac}.)

The integrand is $0$ if $\om$ will not enter $E$ in the $[0,t]$-future,
otherwise it is the time from its ``first'' future visit to $E$ to
its ``first'' one after time $T$.

One may define the {\EM arrival time} to $E$ 
$\xi(\om)=\xi_E(\om)$ by
      $$\xi_E(\om):=\inf\cO_\om E\cap\,]0,\I[$$
(One might take $[0,\I[$ instead of $]0,\I[$, changing $\xi$ only
on a null set.)

As in the discrete case (Rmk.\ \ref{Rmk:1Line}),   
If we knew that $\xi_E$ is $L^1$ in $\om$, we could prove (\ref{eq:CtsKacI})
in one line.

\end{Rmk}

For the rest of this section, unless stated otherwise, the setting is that
of Borel measure-preserving action on a standard Borel probability space
$(\OM,\cB,\mu)$, and $E$ borel.

\subsubsection{Some Assertions with No Discrete Analog}
\label{SS:NoDs}

Note first, that the mere defining property (\ref{eq:EnhFn1}) of the
expectation of an enhanced \fn, which can be written as:
\begin{equation} \label{eq:NoDs1}
 \bE \LP\int_G h(x)f(x\om)\,d\La(x)\RP=\bE f(\om)\frac{d\La}{d\la}
 \int_G h(x)\,d\la(x)
\end{equation}
is a case of CHG with weighted graph $F(x,y;\om)=h(xy^{-1})f(x\om)$
and $1$-simplex of measures $(\La,\la)$.

An application which involves two invariant measures, thus has no discrete
analog, comes from the weighted graph $f(x\om)h(xy^{-1})$ and $1$-simplex of
measures $(\cnt_G,\La)$, where $f\ge0$ is $\frac{d\cnt}{d\la}\mu$-integrable
(really it is enough that $f$ is Borel and $\cO_\om f$ is a.e.\ supported
on a countable set)
and $h$ is non-negative s.t.\ $h(x^{-1})$ is $\La$-integrable.
Note that the $\sigma$-finiteness requirement of CHG is satisfied -- the
set of relevant $y$'s for each fixed $\om$ is $\La$-$\sigma$-finite.
To check measurability of the obtained $1$-chain we may view $x$ and $yx^{-1}$
as two independent variables. This separates $f$ and $h$.
(Note that $yx^{-1}$ is restricted to a fixed $\La$-$\sigma$-finite set.)
It suffices to check test-\fns\ on $G^2$ where these variables are similarly
separated, i.e.\ of the form $g_1(x)g_2(yx^{-1})$.
For the variable $x$ we may use Rmk.\ \ref{Rmk:CntMeasFub}.
We conclude that the obtained $1$-chain is measurable,
and CHG gives the following:

\begin{equation} \label{eq:NoDs2}
\mbox{for }g(\om):=\int_G f(x\om)h(x)\,d\cnt_G,\quad
\bE g(\om)\frac{d\La}{d\la}=
\LQ\int_G h(x^{-1})\,d\La(x)\RQ\bE f(\om)\frac{d\cnt}{d\la}
\end{equation}

For example, referring to Exm.\ \ref{Exm:InfMea} item \ref{it:InfMea-a},
if $S\st\bR$ has finite measure w.r.t.\ some Hausdorff measure
$\La$ on $\bR$ then the union of the translates of a smooth curve by $S$
(with multiplicity) has finite $\frac{d\La}{dt}\mu$-mass.

The formulas in this \S\ have applications to the problem of recovering the
original measure from the infinitesimal measure (\S\ref{SS:REC}).

\subsubsection{Some Analogs to the Discrete Case}

One may formulate analogs to the facts in \S\ref{SS:DsAsso}, the analogy
being both in formulation and in proof.

Note that assertions like the \cts\ Kac Thm.\ \ref{Thm:CtsKac} are more
lucid in the case when for a.a.\ $\om\in\OM$ $\cO_\om E$ is {\em closed},
and in most of what follows it will be assumed that $\cO_\om E$
is {\em discrete} in $\bR$. Such is the case for Exm.\ \ref{Exm:InfMea}
item \ref{it:PoisBrw} (Poisson distribution), or for other stationary
processes defined by probability over the discrete subsets of $G$.
For $\bR$-action, we have by Poincar\'e's recurrence that when $\cO_\om E$
is a.e.\ discrete, it is unbounded above and below a.e.\ in $\spn E$. 

Note that if $\cO_\om E$ is a.e.\ discrete then $E$ is always $\mu$-null
(Fubini for $\mu$ and Haar in $G$). Of course, it may be non-null w.r.t.\
infinitesimal measures.  

\begin{Prop}
Let $\bR$ act in a measure-preserving Borel manner on a standard probability
space $(\OM,\mu)$.
Suppose $E\st\OM$ is s.t.\ $\forall\om\in\OM$ $\cO_\om E$ is {\em discrete}
in $\bR$.
Then the $\bZ$-action on $T_E:E\to E$ given by the
{\EM induced transformation} (see \S\ref{SS:Asso1D}):
$$T_E(\om):=T^{\min(\cO_\om E\cap]0,\I])}\om\quad\om\in E$$
is Borel, and, w.r.t.\ $\frac{d\cnt}{dt}\mu$ on $E$, defined a.e.\ and
measure preserving.
\end{Prop}

\begin{Prf}
First, let $\OM_1$ be the set
$$\{\om\in\OM:\cO_\om E\mbox{ is unbounded above and below}.\}$$
$\OM_1$ is invariant, it is Borel since the discrete set $\cO_\om E$ depends
in a Borel manner on $\om$ (Prop.\ \ref{Prop:EGDelta}), and by Poincar\'e
recurrence it is conull. Therefore, by the definition of
$\frac{d\cnt}{dt}\mu$, $E\cap\OM_1$ is $\frac{d\cnt}{dt}\mu$-conull in $E$,
and we may and do replace $\OM$ by $\OM_1$ and $E$ by $E\cap\OM_1$,
making $T_E$ 1-1, onto and Borel on $E$.

Let $f:E\to\bR^+$ be Borel $\frac{d\cnt}{dt}\mu$-integrable.

Weighted graph: the arrows $(t_0,t_1)$ s.t.\ $t_1<t_0$ and
$[t_1,t_0]\cap\cO_\om E=\{t_1,t_0\}$, with weight $f(T^{t_0}\om)$.

$1$-simplex of measures: $(\cnt_\bR,\cnt_\bR)$.

Check that CHG is applicable. Note that the weighted graph is a Borel \fn\
of $(\mbox{the discrete }\cO_\om E, t_0, t_1)$.

We have: source: $f(\om)$; target: $f(T_E(\om))$.

\NOT{
CHG is applicable: to check measurability of the $1$-chain: the closed
$\cO_\om E$ depends measurably on $\om$. So does the sum of
$\cO_\om f$ on the intersection of $\cO_\om E$ with some fixed closed.
The graph of the weighted hypergraph is Borel measurable in
$\cO_\om E$, the operation of intersection of closed sets is Borel, so
is the mapping $(\mbox{closed }K\st\bR^2,y\in\bR)\mapsto\{x:(x,y)\in K\}$.
Now take as test-\fns\ compact sets $\st\bR^2$. 
}

CHG implies $f(T_E(\om))$ is $\frac{d\cnt}{dt}\mu$-integrable with the same
integral as $f(\om)$. Consequently, $T_E$ is measure-preserving.
\end{Prf}

In order to show what happens, in the \cts\ case, to assertions such as
Propositions \ref{Prop:KacDist}, \ref{Prop:InvDist} and \ref{Prop:SsDist},
let us state:

\begin{Prop}
For measure-preserving Borel $\bR$-action on a standard $(\OM,\mu)$,
let $E\st\OM$ be Borel and assume $\spn E=\OM$ and
$\cO_\om E$ always discrete. Define the return time $\rho=\rho_E$
and the arrival time $\xi=\xi_E$ as usual ($\rho_E$ is $0$ outside $E$).
Let the superscript $(-)$ refer to the inverse action. Then
\begin{itemize}
\item
$\xi_E$ and $\xi_E^{(-)}$ have the same ordinary (i.e.\ $\mu$-) distribution.
\item
$\rho_E$ and $\rho_E^{(-)}$ have the same $\frac{d\cnt}{dt}\mu$-distribution.
\item
For any $a\ge0$, the sets $\{\xi_E=a\}$ and $\{\xi_E^{(-)}=a\}$ have the
same $\frac{d\cnt}{dt}\mu$-mass.
\item
For any $a\ge0$, the $\frac{d\cnt}{dt}\mu$-masses of $\{\xi_E=a\}$ and
$\{\rho_E>a\}$ are equal.
\end{itemize}
Now let $s:\bR^+\to\bR^+$ be Borel and let $S(x):=\int_0^xs(t)\,dt$. We have:
\begin{itemize}
\item
$$\bE s(\xi_E(\om))=\bE S(\rho_E(\om))\frac{d\cnt}{dt}$$
($s\equiv1$ gives Kac's.)
\end{itemize}
(Note that the last assertion is {\em not} a direct consequence of the
preceding one -- see Exm.\ \ref{Exm:InfMea} item \ref{it:Stoc}.)
\end{Prop}

\begin{Prf}
Exercise in using CHG, in analogy with \S\ref{SS:Asso1D}. 
\end{Prf}

Thus, $\xi$ is in $L^p$ iff $\rho$ is in $L^{p+1}$ for $\frac{d\cnt}{dt}\mu$.

\beware In the two last assertions of the above Cor.\ one could,
in retrospect, integrate the infinitesimal measure of $\{\chi=a\}$
as in Exm.\ \ref{Exm:InfMea} item \ref{it:Stoc}. One may try to formulate
general conditions for the possibility of such integration for a stochastic
variable.\beware

As an exercise, one may formulate and prove assertions about
{\em return and arrival to two sets}, in analogy with Prop.\
\ref{Prop:TwoSets}.

Now, for $\bR$-action with discrete $\cO_\om E$, an analog of
Thm.\ \ref{Thm:JointDist} about repeated return and arrival
can be formulated, referring to ordinary ($\mu$-) and
$\frac{d\cnt}{dt}\mu$- distributions. We mention just the analog of
(\ref{eq:KacDec}) (see Rmk.\ \ref{Rmk:KacDec}) which will give us an
opportunity to apply CHG to a weighted $2$-hypergraph,
composed of $2$-simplices.

Assume measure preserving Borel $\bR$-action on a standard $(\OM,\mu)$,
$E\st\OM$ Borel with $\spn E=\OM$ and $\cO_\om E$ always discrete.

Let $s,s':\bR^+\to\BAR{\bR^+}$ be Borel,
and consider their {\em convolution}: $(s*s')(x)=\int_0^xs(t)s'(x-t)\,dt$.

The weighted $2$-hypergraph:
$$\{(t_0,t_1,t_2):t_1<t_0<t_2,\;
[t_1,t_2]\cap\BAR{\cO_\om E}=\{t_1,t_2\}\}$$
with weight $s(t_2-t_0)s'(t_0-t_1)$.

The $2$-simplex of measures: $(dt,\cnt_\bR,\cnt_\bR)$.

Check that CHG is applicable. The three vertices are:
\par\medskip
the $0$-vertex: $s(\xi_E(\om))s'(\xi_E^{(-)}(\om))$\par
the $1$-vertex: $(s*s')(\rho_E(\om))$\par
the $2$-vertex: $(s*s')(\rho_E^{(-)}(\om))$\par
\par\medskip
and one gets the following \cts\ variant of (\ref{eq:KacDec}) -- a
``decomposition'' of Kac's formula: 
\begin{equation} 
\bE \LQ s(\xi(\om))s'(\xi^{(-)}(\om))\RQ=
\bE \LQ(s*s')(\rho(\om))\frac{d\cnt}{dt}\RQ 
\end{equation}
Kac's obtains from $s=s'=1_{]0,\I[}$; then $(s*s')(x)=x$.
\par\bigskip

Let us conclude with further analogs of \S\ref{SS:Asso1D}.
Retain the setting of $\bR$-action and $E\st\OM$ Borel with $\spn E=\OM$
and $\cO_\om E$ always discrete.
Let $s:\bR^+\to\BAR{\bR^+}$ be Borel and let $S(x):=\int_0^xs(t)\,dt$.

Consider our usual Kac graph: 
$$\{(t_0,t_1):t_0>t_1,\;[t_1,t_0]\cap\BAR{\cO_\om E}=\{t_1\}\}$$
and the $1$ simplex of measures $(dt,\cnt_\bR)$.

But now we turn the graph into a weighted graph, giving weights to its
$1$-simplices, in two different manners:

Let $f:\OM\to\bR^+$ be Borel and consider the
weight $f(T^{t_0}\om)s(t_0-t_1)$, CHG gives:
\begin{equation} \label{eq:Ctsfs}
\int_{\OM}f(\om)s(\xi^{(-)}\om)\,d\mu(\om)=
\int_E\LP\int_0^{\rho(\om)}f(T^t\om)s(t)\,dt\RP\,\frac{d\cnt}{dt}\,d\mu(\om)
\end{equation}

Now let $f:E\to\bR^+$ be Borel and consider the
weight $f(T^{t_1}\om)s(t_0-t_1)$. CHG gives:
\begin{equation} \label{eq:CtsfSs}
\int_{\OM}f\LP T^{-\xi^{(-)}\om}(\om)\RP s(\xi^{(-)}\om)\,d\mu(\om)=
\int_E f(\om)S(\rho(\om))\,\frac{d\cnt}{dt}\,d\mu(\om)
\end{equation}

Note that, in fact, (\ref{eq:Ctsfs}) for $s\equiv 1$, with suitable $f$'s,
encompasses both (\ref{eq:Ctsfs}) and (\ref{eq:CtsfSs}).

\begin{Rmk} \label{Rmk:CtsFlow} (compare with Rmk.\ \ref{Rmk:DsFlow}).
Ambrose and Kakutani \cite{AmbroseKakutani} (see also \cite{Jacobs},
\cite{Nadkarni})
have shown, for measurable $\bR$-action, that, taking apart an invariant
subset where the action is trivial on the Boolean Algebra of measurable
sets modulo null sets, there always exists an $E$ s.t.\ $\cO_\om E$ is
discrete a.e.\ and $\spn E$ is conull. Then $\OM$ has the structure of the
``flow under a \fn'' construction (see \cite{Petersen}, p.~11),
while (\ref{eq:Ctsfs}) for $s\equiv1$ 
shows that $\mu$ is indeed the measure making
$(\om,\mu,(T^t)_{t\in\bR})$ the result of the ``flow under a \fn''
construction starting from the discrete system
$(E,\frac{d\cnt}{dt}\mu|_E,T_E)$ and the function $\rho_E:E\to[0,\I[$.
Note that by Kac's $\rho_E$ has integral $1$ on
$(E,\frac{d\cnt}{dt}\mu|_E)$.

Thus, in this case $\OM$ can be recovered from $E$ and $\mu$ can be recovered
from the infinitesimal measure $\frac{d\cnt}{dt}\mu$. More on this in
\S\ref{SS:REC}.
\end{Rmk}

One can play with CHG to find analogs, for \cts\ preordered Time,
(take $G=\bR^n$), to what is said in \S\ref{SS:DsPOG}. In particular,
if $\cO_\om E$ is always closed, one can define the return and arrival
duration and epoch, and one has formulas analogous to \S\ref{SS:DsPOG},

As an example of such statement, one has the following, which for
$n=1$ and the usual ordering gives \cts\ Kac's: 

\begin{Prop} \label{Prop:CtsEpoDur}
Let $\bR^n$ act in a Borel and measure-presrving manner on a standard
probability space $(\OM,\cB,\mu)$.

Fix a preordering $\ge$ in $\bR^n$ with $F_\sigma$ graph.

\NOT{ 
Given a measure-preserving Borel $\bR^n$-action on a standard
$(\OM,\mu)$.
}

Let $E\st\OM$ be Borel with $\cO_\om E$ always {\em closed}.
(W.l.o.g.\ one may think of $\OM=$ the Borel space of closed subsets
$\st\bR^n$ with some shift-invariant probability measure and
$E=\{\om:0\in\om\}$.) 

Let $\La$ be some invariant measure in $\bR^n$ which is a completion of a
Borel measure and has the universal measurability Fubini property
(Def.\ \ref{Def:MeasFub}).

Denote the Lebesgue measure in $\bR^n$ by $dx$.

Then:
\begin{itemize}
\item
The expectation of the Lebesgue mass of the arrival duration:
$$\{x\in\bR^n: [0,x]\cap\cO_\om E=\es\}$$
is equal to the expectation of the Lebesgue mass of the arrival duration
w.r.t.\ the inverse action.

Also, for any $z\in\bR^n$, the $\frac{d\cnt}{dx}\mu$-mass of the set of
$\om$ s.t.\ $z$ is in the arrival duration is equal to the same for
the inverse action.

\item
Suppose that for all $\om$ the set
$$\{y\in\cO_\om E:\exists x\in\bR^n\:y<x,\:[y,x]\cap\cO_\om E=\{y\}\}$$
($x>y$ means $x\ge y\:\&\:x\ne y$)\par
is $\La$-$\sigma$-finite.

Then the $\frac{d\La}{dx}\mu$-integral over $E$ of the Lebesgue mass of the
return duration:
$$\{x\in\bR^n: x>0,\:[0,x]\cap\cO_\om E=\{0\}\}$$
is equal to the expectation of the
$\La$-mass of the arrival epoch w.r.t.\ the inverse action:
$$\{x\in\bR^n: x>0,\:[-x,0]\cap\cO_\om E=\{-x\}\}$$
\end{itemize} 
\end{Prop}
\qed

Here one may think of $E$ s.t.\ $\cO_\om E$ is discrete, with the
$\La=\cnt$. For a non-discrete case take
Exm.\ \ref{Exm:InfMea} item \ref{it:PoisBrw} -- the case of Brownian
motion in $R^n$ or the case of Poisson's distribution, but $E=$
the $\om$'s s.t.\ $0$ is distanced $<r_0$ from a point in the
discrete set $\om$. Here invariant measures different from $\cnt$
come into play. As another example check the above for the \cts\
(and simpler) analog of Exm.\ \ref{Exm:Tor} (triangle in the torus
$\bT^2$ acted by $\bR^2$).

\subsubsection{The Nearest Point}
\label{SS:Nearest}
This construction gives some multi-dimensional \cts\ analog to the
one-dimensional ``future until the first return'' engaged in Kac's and in
the ``flow under a \fn'' construction. It will be applied in
\S\ref{SS:CntCls}.

Let $G$ be a {\em (connected) unimodular Lie group}
acting in a Borel and measure preserving manner on a standard probability
$(\OM,\mu)$.
Choose a positive definite symmetric bilinear form on the tangent space
$T_eG$ and by right-shifting it to each tangent space $T_xG$ via the
derivative at $e$ of $y\mapsto yx$, construct a Riemannian metric,
turning $G$ into a Riemannian manifold with right-invariant metric. Note
that since the exponential map at $e$ defined by geodesics (see \cite{Hicks})
is the same as the Lie-group-theoretic exponential map, the
exponential map is everywhere defined, hence $G$ is complete as a Riemannian
manifold.

Let $E\st\OM$ and assume for simplicity that $\cO_\om E$ is always
{\em closed}.

What we have in mind is to take the graph 
\begin{equation} \label{eq:Nearest}
\{(x,y)\in G^2: y \mbox{ is the unique nearest point to }x
\mbox{ in }\cO_\om E.\}
\end{equation}

To this end, the following proposition is helpful:

\begin{Prop} \label{Prop:Geodesic}
Let $X$ be a complete Riemannian manifold. Let $Y\st X$ be closed. Then 
if $x_0\in X$ is a point where $x\mapsto\dist(x,Y)$ (the distance from $x$ to
$Y$) is differentiable giving gradient $-v\in T_xX$ (that is, it has a
differential equal to the linear functional $w\mapsto\LA -v,w\RA$ where
$\LA,\RA$ is the inner product), then there is a unique point in $Y$
nearest to $x$, that point being
$$\exp_{x_0}(\dist(x_0,Y)v)$$
(that is, the point with parameter $1$ on the
geodesic emanating from $x$ with tangent vector $\dist(x_0,Y)v$
-- see \cite{Hicks}).

Note that since $x\mapsto\dist(x,Y)$ is Lipschitz, it is differentiable
a.e.\ by Rademacher's Theorem (see \cite{Federer} \S3.1.6.).
(Note that in the \nbd\ of each point the Riemannian metric is
Lipschitz-equivalent to the Euclidean metric of a coordinate chart.)
Thus for a.a.\ $x$ there is a unique nearest point.
\end{Prop}

\begin{Prf}
First, in a complete Riemannian manifold any bounded closed set is compact,
therefore there are nearest points $y$, and it suffices to prove that
if there is a differential giving gradient $-v$ then each nearest point
$y$ equals $\exp_{x_0}(\dist(x_0,Y)v)$.

Let $t\mapsto\exp_{x_0}tv'$ be the shortest geodesic from $x_0$ to $y$
(this exists in a complete Riemannian manifold), with $\exp_{x_0}v'=y$.
(thus $\|v'\|=\dist(x_0,Y)$). Clearly,
$$\dist\LP\exp_{x_0}tv',y\RP=(1-t)\|v'\|\quad 0\le t\le 1,$$
and by computing the derivative at $x_0$ in this direction we find
$$\LA -v,v'\RA=-\|v'\|.$$ 
But $\|v\|\le1$, since the Lipschitz constant of $x\mapsto\dist(x,Y)$
is $\le1$. This implies $v'=\|v'\|v$, and $y=\exp_{x_0}(\dist(x_0,Y)v)$.

\NOT{
Now suppose there is a unique nearest point $y_0$ and prove there is a
G\^ateau derivative to $x\mapsto\dist(x,Y)$ at $x_0$. 

Let $\eps_n\downarrow0$, and let $B_n=B_n(y_0,\eps_n)$ be the {\em open}
ball. By the uniqueness of $y_0$,
$$\dist(x_0,Y\sm B_n)>\dist(x_0,Y)=\dist(x_0,y_0),$$
hence $\dist(x,Y\sm B_n)>\dist(x,Y)$ for $x$ in a \nbd\ $U_n$ of $x_0$.
Thus, for $x\in U_n$,
$$\dist(x,Y)=\dist(x,Y\cap B_n)\le\dist(x,y_0),$$
while for $x=x_0$, $\dist(x_0,Y)=\dist(x_0,y_0)$. 

Thus
$$\dist(x,Y)-\dist(x_0,Y)\le\dist(x,y_0)-\dist(x_0,y_0)$$
therefore (where a coordinate chart near $x_0$ is used),
\begin{equation} \label{eq:GeodesicLe}
\limsup_{t\downarrow0}\frac{1}{t}\LQ\dist(x_0+tv,Y)-\dist(x_0,Y)\RQ\le
\LA\frac{\partial}{\partial x}\dist(x,y)|_{x_0,y_0},v\RA
\end{equation}
the function $(x,y)\mapsto\dist(x,y)$ being $\cC^\I$ ??????

On the other hand, every limit of quotients of the form
$$\lim_k\frac{1}{t_k}\LQ\dist(x_0+t_kv,Y)-\dist(x_0,Y)\RQ\quad
t_k\downarrow0$$
has the form
$$\lim_k\frac{1}{t_k}\LQ\dist(x_0+t_kv,y_k)-\dist(x_0,Y)\RQ
\ge\limsup_k\frac{1}{t_k}\LQ\dist(x_0+t_kv,y_k)-\dist(x_0,y_k)\RQ$$
with $y_k\in Y\cap B_n$, while the last $\limsup$ is bounded by the
bounds of the partial derivative w.r.t.\ x at $(x_0,y)$ of 
$(x,y)\mapsto\dist(x,y)$ for $y\in B_n$. Taking $n\to\I$, we find that
\begin{equation} \label{eq:GeodesicGe}
\liminf_{t\downarrow0}\frac{1}{t}\LQ\dist(x_0+tv,Y)-\dist(x_0,Y)\RQ\ge
\LA\frac{\partial}{\partial x}\dist(x,y)|_{x_0,y_0},v\RA
\end{equation}
And (\ref{eq:GeodesicLe}) and (\ref{eq:GeodesicGe}) imply that we have
our G\^ateau derivative.

(The last part of the proof is, without naming it, working with
subdifferentials of the not necessarily differentiable \fns.)
}

\end{Prf}

\begin{Prop} \label{Prop:MeasNearest}
In the above situation, the set of $(x,Y)$, $x\in X$, $Y\st X$ closed,
satisfying: there is a unique nearest point to $x$ in $Y$, is closed in
the product $X\times2^X$
($2^X$ denotes the Borel space of closed subsets of $X$ -- 
see \S\ref{SS:MEAS}),
and the mapping which maps each $(x,Y)$ in this set to the unique nearest
point is Borel.
\end{Prop}

\begin{Prf}
Let $D$ be a fixed dense set in $X$. There is a unique nearest point iff
for each $n\in\bN$ there are rational $q_1<q_2$ and $x\in D$ with
$q_2-q_1<1/n$, the closed ball $B(x,q_1)$ not intersecting $Y$ and the
closed ball $B(x,q_2)$ intersecting $Y$ in a non-empty set $\st B(x,1/n)$.
$y$ is the unique nearest point if one may require also
$y\in B(x,1/n)\:\forall n$.
Since a closed ball is a continuous \fn\ of its center and radius
(Hausdorff topology in the space of compact sets), and intersection of
closed sets is a Borel \fn, we are done.
\end{Prf}

Return now to our connected unimodular Lie group $G$ acting on $\OM$
(choose a Haar measure $\la$ on $G$) and to our $E\st\OM$,
and assume $\cO_\om E$ is always {\em countable closed}.
For $\om\in\OM$, denote by $\pi(\om)=\pi_E(\om)\in G$ the
unique nearest point to $0$ in $\cO_\om E$, if it exists
(If there is no unique nearest point, $\pi_E(\om)$ is undefined);
for $\om\in E$, denote by $P(\om)=P_E(\om)$ the set of $x\in G$ s.t.\ $0$ is
the unique nearest point to $x$ in $\cO_\om E$.

Taking into account Prop.\ \ref{Prop:EGDelta}
and Prop.\ \ref{Prop:MeasNearest} above, we have
that the set $\OM'$ where $\pi_E$ is defined is Borel and
$\pi_E:\OM'\to G$ is Borel, moreover $\cO_\om\OM'$ is conull
(w.r.t.\ Haar) $\forall\om\in\spn E$. Also 
$\{(\om,x):x\in P_E(\om)\st\OM\times G\}$ is Borel.

So consider the $\om$-dependent graph (\ref{eq:Nearest}).
Since the Riemannian metric is right-invariant, so is the graph,
and it will remain so if we take a weight
$f(x\om)s(xy^{-1})$ or $f(y\om)s(xy^{-1})$, $s,f$ Borel.
Take the $1$-simplex of measures $(\la,\cnt_G)$.
The $\sigma$-finiteness requirements of CHG are satisfied.
By the above, the relation: ``$(x,y)$ belongs to the graph at $\om$''
is Borel in $(x,y,\om)$. Also $\cO_\om E$ is always countable.
Therefore (see Rmk.\ \ref{Rmk:CntMeasFub}) the obtained $1$-chain is
Borel-measurable.
By Prop.\ \ref{Prop:Geodesic}, 
the source of the unweighted graph is a.e.\ $1_{\spn E}$.
CHG gives now the following statements, in the spirit of Kac's 
(compare (\ref{eq:Ctsfs}) and (\ref{eq:CtsfSs})):

\begin{Prop} \label{Prop:Nearest}
Let $G$ be a connected unimodular Lie group acting in a Borel and
measure-preserving manner on a standard probability space $(\OM,\mu)$.
Let $E\st\OM$ be Borel s.t.\ $\cO_\om E$ is always 
non-empty countable closed. (Thus $\spn E=\OM$.) 
Endow $G$ with a right-invariant Riemannian
structure and let $\la$ be a Haar measure in $G$.

For $\om\in\OM$, let $\pi_E(\om)\in G$ be the unique nearest point
(w.r.t.\ the Riemannian metric) to $0$ in $\cO_\om E$, if it exists.
If there is no unique nearest point, $\pi_E(\om)$ is undefined.

For $\om\in E$, let $P(\om)$ be the set of $x\in G$ s.t.\ $0$ is the unique
nearest point (w.r.t.\ the Riemannian metric) to $x$ in $\cO_\om E$.

Then ((i) is, of course, a special case of (ii) which is a special case of
(iii) or (iv), these being special cases of (iii) for $s\equiv1$):
\begin{itemize}
\item[(i)]
(consider the above unweighted graph (\ref{eq:Nearest}). The simplex of
measures is always $(\la,\cnt_G)$)
$$\int_E\la(P(\om))\,\frac{d\cnt}{d\la}\,d\mu(\om)=1$$
\item[(ii)]
(consider the above graph (\ref{eq:Nearest} with weight $s(xy^{-1})$)
\begin{equation}
\int_{\OM} s\LP\LP\pi_E\om\RP^{-1}\RP\,d\mu(\om)=
\int_E\LP\int_{P(\om)}s(x)\,d\la(x)\RP\,
\frac{d\cnt}{d\la}\,d\mu(\om)
\end{equation}
\item[(iii)]
Let $f:\OM\to\BAR{\bR^+}$ and $s:G\to\BAR{\bR^+}$ be Borel. Then
(consider the above graph (\ref{eq:Nearest}) with weight $f(x\om)s(xy^{-1})$):
\begin{equation} 
\int_{\OM}f(\om) s\LP\LP\pi_E\om\RP^{-1}\RP\,d\mu(\om)=
\int_E\LP\int_{P(\om)}f(x\om)s(x)\,d\la(x)\RP\,
\frac{d\cnt}{d\la}\,d\mu(\om)
\end{equation}
\item[(iv)]
Now let $f:E\to\BAR{\bR^+}$ and $s:G\to\BAR{\bR^+}$ be Borel. Then
(consider the above graph (\ref{eq:Nearest}) with weight $f(y\om)s(xy^{-1})$):
\begin{equation} 
\int_{\OM}f\LP\pi_E\om\RP s\LP\LP\pi_E\om\RP^{-1}\RP\,d\mu(\om)=
\int_E f(\om)\LP\int_{P(\om)}s(x)\,d\la(x)\RP\,\frac{d\cnt}{d\la}\,d\mu(\om)
\end{equation}
\end{itemize}
\end{Prop}

\newpage
\section{Equidecomposability and the Totality of Invariant Measures}
\label{S:ED}
\subsection{Introduction}
In this section the acting group $G$ is assumed {\em discrete}.

As we have seen, Kac's theorem is a special case of equidecomposable
\fns\ which trivially have the same integral with respect to the
$G$-invariant measure.

A significant fact, however, is that the property of \fns\ to be
equidecomposable (such as $\rho_E(\om)$ and $1_{\spn'E}$ in the formulation
of Kac's in \S\ref{SS:EqiChn}) has nothing to do with the particular
invariant measure. This suggests an investigation of the relationship
between equidecomposability and {\em the totality of invariant measures} in
suitable frameworks.

The results which will be presented try to state reverse implications:
if two \fns\ have the same integral (or one has always a greater
integral) w.r.t.\ a comprehensive set of
invariant measures, then they are close to being equidecomposable in a
suitable sense (or one may find a \fn\ equidecomposable with one
\fn\ and dominated by the other, etc.).

As an example to such ``reverse implication'', one deduces directly from
the Hahn-Banach Thm.\ that for a group acting by homeomorphisms on a
compact space, two \cts\ \fn\ have the same integral w.r.t.\ all
invariant regular probability measures iff their difference can be
uniformly approximated by sums of \fns\ of the form $f(\om)-f(x\om)$,
$f(\om)$ \cts, $x\in G$. Such \fns\ are ``close to being equidecomposable''
via signed \cts\ \fns. Our interest will lie, however, with
equidecomposability via {\em nonnegative} \fns, where the matters are
a bit less simple. Yet, our main tools will still be theorems akin to convex
separation, which are, of course, a part of the Hahn-Banach philosophy.

\begin{Rmk}
Equidecomposability of \fns\ has been studied extensively by Friedrich
Wehrung (see, for example, \cite{WehrungInj1}, \cite{WehrungInj2},
\cite{WehrungHB}),
using his algebraic methods concerning positively ordered monoids
and in connection with the Banach-Tarski paradox (see \cite{WagonBT}).
He calls it {\em \cts}
equidecomposability, ``\cts'' referring to the \fns\ being
$\BAR{\bR^+}$-valued, while the Banach-Tarski paradox about decomposability
of sets refers to $\BAR{\bZ^+}$-valued \fns.
It seems that our methods and results, being more functional analytic, are
somewhat different.
\end{Rmk}

\subsection{Upper Semi-Continuous and Baire Lower Semi-Continuous Functions}
\label{SS:ULSC}

For the rest of Section \ref{S:ED} our groups will be discrete
(not necessarily countable).
Thus {\em amenable} will mean: {\em amenable as a discrete group}.

It seems better, when one wishes to speak about the totality of invariant
measures, to deal with a compact or at least locally compact space.
When a group $G$ acts on it by homeomorphisms, we have a compact $G$-space.
Note that $G$-measurable spaces can often be related to compact $G$-spaces
(see \S\ref{SS:CP}, which deals with the continuous 2nd-countable case,
containing the case of countable discrete $G$).
By Stone's duality (\cite{HalmosBOO}) $G$-Boolean
algebras are equivalent to compact totally disconnected $G$-spaces.

\begin{Rmk} 
Moreover, if a Boolean $\sigma$-field $\cB$ of subsets of some set $\OM$
is given, with a given fixed Boolean $\sigma$-ideal $\cJ$ of
{\em null sets}, then many measure theoretic concepts in $\OM$ correspond
canonically to topological concepts in the Stone space $\cS$ of $\cB/\cJ$
(see \cite{Ellis}). In particular, a measure in $\cB$ zero on $\cJ$
induces a measure in $\cS$; bounded measurable \fns\ from $\OM$ to some
metrizable compact space $K$, modulo null sets,
are in 1-1 correspondence with \cts\ \fns\ on $\cS\to K$, (taking
$K=[0,\I]$ deals with unbounded \fns), where convergence a.e.\ corresponds
to convergence on $\cS$ modulo meager sets.  
\end{Rmk}

Let us fix some notations.

By a {\EM p.m.\ (probability measure)} we will mean a regular Borel
probability measure on a compact or locally compact $\OM$.

We shall deal with {\EM u.s.c.\ (upper semi-\cts)} \fns\ and
{\EM l.s.c.\ (lower semi-\cts)} \fns. Unless stated otherwise,
u.s.c.\ \fns\ will be assumed $[-\I,\I[$-valued and l.s.c.\ \fns\
$]-\I,\I]$-valued.

Special attention will be given to nonnegative {\EM Baire l.s.c.\
(b.l.s.c.)} \fns\ in the sense of \cite{HalmosMEA} Ch.\ X.
These are nonnegative l.s.c.\ \fns\ $f$ s.t.\ the open\NOT{$\{f>0\}$ is
$\sigma$-compact, and $f$ is Baire-measurable, which is equivalent to}
$\{f>a\}$ is $\sigma$-compact for all $a\in\bR^+$.

\begin{Rmk} \label{Rmk:USC-LSC}
Note the following well-known facts, where all \fns\ are real on a compact
space $\OM$: ($f>g$ means $\forall\om\in\OM\:\:f(\om)>g(\om)$)
\begin{enumerate}
\item \label{it:USC-LSCInf}
any u.s.c.\ \fn\ is the pointwise infimum of the
\cts\ \fns\ dominating it.
\item \label{it:USC-LSCInt}
Integration w.r.t.\ a (positive)
regular finite measure commutes with taking pointwise infima of a
directed downward family of u.s.c.\ \fns.
\item \label{it:USC-LSCCof}
Any family, directed downwards w.r.t.\ $\le$, of \cts\ \fns\ whose
infimum is $f$ (necessarily u.s.c.) is cofinal downwards w.r.t.\ the
\cts\ \fns\ strictly $>f$.
\item \label{it:USC-LSCSum}
If $(f_i)_i$ is a {\em finite} family of u.s.c.\ \fns,
then any \cts\ $h>\sum_if_i$ can be written as $h=\sum_ih_i$,
$h_i$ \cts\ and $h_i>f_i$. (this follows from item \ref{it:USC-LSCCof},
for the family of all $h$'s expressible as such sums.)
\item
Dual facts hold for l.s.c.\ \fns.
\item
The sum of a (not necessarily countable) family of {\em nonnegative} l.s.c.\
\fns\ is l.s.c. 
\item \label{it:USC-LSCCont}
If a \cts\ \fn\ $f$ is the sum of a (not necessarily countable)
family of nonnegative l.s.c.\ \fns\ $(f_\al)$, then each $f_\al$ is 
\cts\ (this follows from $f_\al$ being also u.s.c., since
$f_\al=f-\sum_{\be\ne\al}f_\be$). 

Moreover, then the summation is uniform on compacta (Dini's Thm.)

\item
Similar facts hold for \fns\ on a locally compact $\OM$, where
the u.s.c.\ (in particular, \cts) \fns\ are required to
have compact support, and for such \fn\ $g$ one requires, instead of $f>g$,
that $f>g$ {\em on $\supp g$}.

\item \label{it:BLSCCount}
Any b.l.s.c.\ \fn\ is the supremum of a {\em countable} family of \cts\ \fns\
with compact support
(this family can be chosen monotone nondecreasing.)
\end{enumerate}
\end{Rmk}

\begin{Rmk}
Note that if $\bZ$ is acting on a compact $\OM$ and $E\st\OM$ is clopen,
then $\rho_E$ (see \S\ref{SS:EqiChn}) is {\em b.l.s.c.},
this following from the fact that its
$0$-chain is a vertex of a $1$-chain \cts\ in $\om\in\OM$. Note that
by \S\ref{SS:CP}, for any standard Borel $\OM$ on which $\bZ$ acts and
for any Borel $E\st\OM$, $\OM$ can be embedded as a Borel subset in a
totally disconnected compact $\bar{\OM}$ with $E$ an intersection of
$\OM$ with a clopen $\bar{E}$.  
\end{Rmk}

\subsection{Two Theorems}
\label{SS:2Thms}

One has the following two theorem, one dealing with u.s.c.\ \fns\ and the
other with b.l.s.c.\ \fns. They will be presented here with some
corollaries. The proofs are given in \S\ref{SS:Proofs}.

\begin{Thm}\label{Thm:USC}
Let $G$ be a discrete group (possibly uncountable).
Let $\OM$ be a $G$-compact space.

Let $f:\OM\to\bR^+$ and $g:\OM\to\bR^+$ be nonnegative u.s.c.\
(as mentioned above they are assumed finite).

Suppose any $G$-invariant p.m.\ on $\OM$ gives to $f$ a greater or equal
integral than to $g$.

Then every \cts\ \fn\ strictly greater than $f$
is strictly greater than some u.s.c.\ \fn\ finitely decomposable with
$g$ via u.s.c.\ \fns.

Equivalently (see Remark \ref{Rmk:USC-LSC} item
\ref{it:USC-LSCSum}):

Every \cts\ \fn\ strictly greater than $f$
is finitely equidecomposable via nonnegative \cts\ \fns\ with some
(\cts) \fn\ strictly greater than $g$.
\end{Thm}

The phrase: $f'$ and $f''$ being ``finitely equidecomposable via \cts\
\fns'' is self-explanatory: it means that $\exists$ a finite $I\subset
G$ and nonnegative \cts\ \fns\ $f_x,\,x\in I$ s.t.\ $f'=\sum_{x\in
I}f_x$ while $f''(\om)=\sum_{x\in I}f_x(x\om)\;\om\in\OM$. Alternatively,
this may be expressed as: the $0$-chains of $f'$ and $f''$ are the two
vertices of a \cts\ $1$-chain (i.e.\ the coefficient of every $1$-simplex
depends \cts{ly} on $\om$) which is supported in a fixed finite union of 
sets of the form $\{(x_1y,x_2y):y\in G\}$.
Similar expressions will have similar meanings.

\begin{Rmk}
One cannot replace in Thm.\ \ref{Thm:USC}
``\cts\ \fn\ strictly greater than $f$'' by $f$ itself.
A simple example is given by $f\equiv0$ and $g$ being $0$ except at a single
point $x_0$ where it is $1$, where $x_0$ has an infinite $G$-orbit.
\end{Rmk}

Now, take as one of the \fns\ in Thm.\ \ref{Thm:USC} a constant $a$,
to obtain:

\begin{Cor} \label{Cor:USC1}
Let $G$ be a discrete group (possibly uncountable).
Let $\OM$ be a $G$-compact space.

Let $f$ be nonnegative u.s.c.\ on $\OM$. If a constant $a$ is greater
than the supremum of the integrals of $f$ w.r.t.\ all $G$-invariant
p.m.\ then $\exists$ a \fn\ finitely equidecomposable with $f$ via
\cts\ (resp.\ u.s.c.) \fns\ and dominated by $a$.

In other words, for a \cts\ (resp. u.s.c.) \fn\ $f$, the infimum
of the maxima of all \fns\ finitely equidecomposable with it via
\cts\ (resp.\ u.s.c.) \fns\ is equal to the supremum of the
integrals of $f$ w.r.t.\ $G$-invariant p.m.'s.
\end{Cor}

Clearly, in Cor.\ \ref{Cor:USC1} one may replace ``\cts'' by ``member of a
fixed dense subalgebra $\cA$ of $\cC(\OM)$ s.t.\ for any nonvanishing
$f\in\cA$ one has $1/f\in\cA$''. As $\cA$ one may take the set of smooth
\fns\ (for $\OM$ a compact manifold) of the set of \cts\ simple (i.e.\
with finite range) \fns\ (for $\OM$ compact totally disconnected).
Since (see \cite{HalmosBOO}) such $\OM$ is equivalent, by Stone's duality,
to a $G$-Boolean algebra $\cB$, where every finitely additive p.m.\ on
$\cB$ can be uniquely extended to a $\sigma$-additive regular Borel
p.m.\ om $\OM$, this leads to the following statement dealing with
equidecomposition of ``sets'', in the spirit of the Banach-Tarski paradox
(see \cite{WagonBT}):

\begin{Cor} \label{Cor:USC-Boo}
Let $G$ be a discrete group (possibly uncountable).
Let $\cB$ be a $G$-Boolean algebra.

Let $E\in\cB$. If a constant $a$ is greater
than the supremum of the masses of $E$ w.r.t.\ all $G$-invariant
finitely additive p.m.\ on $\cB$ then $\exists$ $p\in\bN$ and a finite
decomposition of $p$ times $E$ into members of $\cB$, with $G$-translates
of the pieces covering every part of $1$ no more than $\lfloor ap\rfloor$
times.
\end{Cor}

\beware
In this, ``Banach-Tarski paradox-like'' context, one may inquire, for
instance:

Can $p$ be given ``in advance''? That is, suppose $q/p$ is greater
than the above supremum.
Are we sure $p$ times $E$ can be decomposed to
pieces whose translates cover every part no more than $q$ times? 
Note that $\cB$ is here arbitrary.

What happens if $q/p$ is {\em equal} to the above supremum? (or if in
Cor.\ \ref{Cor:USC1} ``greater than the supremum'' is replaced by
``greater or equal''?)

Compare with Examples \ref{Exm:MaxAv-Boo1}, \ref{Exm:MaxAv-Boo2}
and \ref{Exm:MaxAv-Boo3}.
\beware

\par\medskip
Now turn to the theorem dealing with b.l.s.c.\ \fns.

\begin{Thm}\label{Thm:LSC}
Let $G$ be a discrete group (possibly uncountable).
Let $\OM$ be a $G$-locally compact space.

Let $f:\OM\to\BAR{\bR^+}$ and $g:\OM\to\BAR{\bR^+}$ be nonnegative b.l.s.c. 
(as mentioned above, they may take the value $+\I$.)

\begin{itemize}
\item[(i)]
suppose:

\begin{itemize}

\item
An orbit of $G$ intersects the set $\{g>0\}$ iff it intersects the set $\{f>0\}$.

\item
If $U$ is the (open) union of the orbits that intersect $\{g>0\}$ and $\{f>0\}$,
then any (non-negative) $G$-invariant Radon measure%
\footnote{Recall that a Radon measure on a locally compact space is a
measure, with every open set measurable, finite on compact sets,
and regular, in the sense that the mass of every measurable set is the
infimum of the mass of opens containing it, and the mass of every open set
is the supremum of the mass of compacts contained in it. If the locally
compact space is 2nd-countable regularity for Borel sets is automatic,
and the mass of every Borel set is the supremum of the mass of compacts
contained in it as well as the infimum of the mass of opens containing it.}%
on $U$ (possibly with infinite mass)
gives the same integral to $f$ and $g$.

\end{itemize}

Then $f$ and $g$ are countably equidecomposable via nonnegative
b.l.s.c.\ \fns.

{\EM note:} if $f$ is \cts, all the \fns\ in the decomposition are
\cts\ and the decomposition of $f$ is uniformly convergent
on compacta -- see Remark \ref{Rmk:USC-LSC} item \ref{it:USC-LSCCont}. 

\item[(ii)]
suppose:

\begin{itemize}

\item
Any orbit of $G$ that intersects the set $\{g>0\}$ intersects the set $\{f>0\}$.

\item
If $U$ is the (open) union of the orbits that intersect $\{g>0\}$, then any
(non-negative) $G$-invariant Radon measure on $U$ (possibly with
infinite mass) gives to $f$ an integral greater or equal than to $g$.
\end{itemize}

Then $f$ is greater or equal than some \fn\
countably equidecomposable with $g$ via nonnegative b.l.s.c.\ \fns.
\end{itemize}
\end{Thm}

\begin{Exm}
From Thm.\ \ref{Thm:LSC} one obtains that if $\OM$ is a $G$-metrizable
compact space, and $x_0\in\OM$ is a point with infinite orbit,
then $1$ and $1-1_{\{x_0\}}$ are countably equidecomposable
(via l.s.c.\ \fns) while $1-1=0$ and $1-(1-1_{\{x_0\}})=1_{\{x_0\}}$
are not equidecomposable. Thus the relation of countable equidecomposability
does not carry over to differences.
\end{Exm}

From Stone's duality point of view, the setting of Thm.\ \ref{Thm:LSC} 
corresponds to $\BAR{R^+}=[0,\I]$-valued finitely additive measures $\mu$
on a $G$-Boolean algebra $\cB$ with dual compact totally disconnected
Stone space $\OM$. The union of clopens in $\OM$ which, as members of $\cB$,
have finite $\mu$-mass is an open, hence locally compact, $U\st\OM$ on which
$\mu$ defines a unique Radon measure. Thus if $E\in\cB$ then
$G$-invariant $\BAR{R^+}$-valued measures $\mu$ with $0<\mu(E)<\I$
correspond to $G$-invariant Radon measures $\mu$ on the open union of
orbits in $\OM$ that intersect the clopen $E$, s.t.\
$0<\int 1_E\,d\mu<\I$. The following corollary of Thm.\ \ref{Thm:LSC} may
be thought of as a ``continuous'' ``topological'' analog of Tarski's
theorem which states, for $\cB$ the field of subsets of some set, that
such $\mu$ exists iff $2E$ is not equidecomposable (as sets) with $E$
(see \cite{Paterson}, \cite{TarskiCA})

\begin{Cor} \label{Cor:LSC-Tar}
Let $G$ be a discrete group (possibly uncountable).
Let $\OM$ be a $G$-locally compact space.

Let $f:\OM\to\bR$ be nonnegative b.l.s.c. 
(as mentioned above, it may take the value $+\I$.)

Let $U$ be the open union of the orbits which contain a point where $f$
does not vanish.

If there is no $G$-invariant Radon measure on $U$ that gives to $f$ the
integral $1$, then for every $0<a<\I$ $af$ and $f$ are countably
equidecomposable via b.l.s.c.\ \fns.

(If $f$ is continuous the \fns\ in
the decomposition are continuous and the decompositions converge
uniformly on compacta.)
\end{Cor}

\begin{Prf}
Just apply Thm.\ \ref{Thm:LSC} (i) to $f$ and $af$.
\end{Prf}

Thm.\ \ref{Thm:LSC} can be applied to Stone-\v{C}ech compactifications of
discrete sets:

Let a group $G$ act (on the left) on a set $S$.
Call a set $T\st S$ {\EM syndetic} if $\exists$ a finite $F\st G$ s.t.\ 
$F\cdot T=\{xt:x\in F,\:t\in T\}=S$. Call a \fn\ $f:S\to\BAR{R^+}$
{\EM strictly syndetically supported} if for some $\eps>0$ the set
$\{f\ge\eps\}$ is syndetic.

Now consider the Stone-\v{C}ech compactification $\be S$, which is a
$G$-compact space. a \fn\ $f:S\to\BAR{R^+}$ extends to a \cts\ \fn, to be
denoted also by $f$, from $\be S$ to the compact $\BAR{R^+}$, which is
thus a b.l.s.c.\ \fn\ $f:\be S\to\BAR{R^+}$. If $f$ is strictly syndetically
supported, there is no $G$-invariant filter $\cF$ on $S$ s.t.\
$\lim_\cF f=0$, since such filter must contain the complements of all
sets $x\cdot\{f\ge\eps\}$, hence must contain $\es$ since some $\{f\ge\eps\}$
is syndetic. Thus there is no non-empty $G$-invariant closed set in $\be S$
on which $f$ vanishes, hence the union of $G$-orbits intersecting $\{f>0\}$
in $\be S$ is the whole $\be S$. From Thm.\ \ref{Thm:LSC} one can now deduce:

\begin{Cor}
Let a group $G$ act on a set $S$. Let $f,g:S\to\BAR{R^+}$ be strictly
syndetically supported.

If $f$ and $g$ have the same integral w.r.t.\ all
$G$-invariant finitely additive probability measures on $S$
(where the integral of an unbounded \fn\ is defined as the supremum
of the integrals of bounded \fns\ majorized by it), then $f$ and $g$
are countably equidecomposable via non-negative \fns\ on $S$,
the decompositions converging uniformly on any set where the relevant sum
is bounded.

If $f$ has greater or equal integral than $g$ w.r.t.\ any $G$-invariant
finitely additive probability measure on $S$, then $f$ is greater or equal
than some \fn\ countably equidecomposable with $g$ via non-negative \fns\
on $S$, with the above uniform convergence property.
\end{Cor}

\begin{Prf}
Consider $f$ and $g$ extended to $\be S$ as \cts\ to $\BAR{R^+}$.
By the remarks preceding the statement of the corollary, the set $U$ in
Thm.\ \ref{Thm:LSC} for $f$ or $g$ on $\be S$ is the whole $\be S$. Hence
Radon measures on $U=\be S$ are finite, and they correspond to the finitely
additive finite measures on $S$. Thus, by Thm.\ \ref{Thm:LSC}, $f$ and $g$
in our case are countably equidecomposable via b.l.s.c.\ \fns\ on $\be S$.
The uniform convergence follows from Dini's Thm.
\end{Prf}

Note that some requirement, such as $f$ and $g$ being strictly syndetically
supported, is needed, as is shown by the example of $f\equiv 0$ and $g$
s.t.\ for every $\eps>0$ $\{g>\eps\}$ is a set with infinite number of
disjoint translates (e.g.\ $g$ tending to $0$ at infinity).

For the case that $G$ is amenable and 
the invariant probability measure is unique (e.g.\ a $\bZ$ action by
irrational rotation on the circle) one has

\begin{Cor}
Suppose $G$ is amenable.

Let us be given a uniquely ergodic $G$-compact space $\OM$, i.e.\ a \cts\
action of $G$ on $\OM$ s.t.\ on $\OM$ $\exists$ a unique invariant probability
measure $\mu$.
(equivalently, $\exists$ a unique {\em ergodic} invariant
probability measure). Suppose, moreover, that the support of $\mu$ is $\OM$

Then two b.l.s.c.\ non-negative \fns\ that have the same $\mu$-integral
are countably equidecomposable via b.l.s.c.\ \fns.
\end{Cor}

\begin{Prf}
Just apply Thm.\ \ref{Thm:LSC}, taking into account the fact that
in our case every orbit is dense (otherwise its closure is an invariant
proper subset, on which $\exists$ some invariant probability measure
different from $\mu$).
\end{Prf}

\begin{Rmk}
When a standard $G$-space $\OM$ is embedded as an (invariant) Borel subset
in a metrizable compact $G$-space $\BAR{\OM}$, the totality of $G$-invariant
measures in $\OM$ is a subset of the totality of $G$-invariant measures in
$\BAR{\OM}$, namely, those giving mass $0$ to $\BAR{\OM}\sm\OM$.
For example, when $\BAR{\OM}$ is constructed as a Stone space of an
(invariant) countable Boolean Algebra $\cB$ forming a basis to the Borel
sets, any {\em finitely} additive (say, invariant) measure on $\cB$
corresponds to a ($\sigma$-additive) measure in $\BAR{\OM}$, the latter
being concentrated in $\OM$ if and only if the original measure was
$\sigma$-additive.
One feels that considering invariant measures in the compact $\BAR{\OM}$
``includes measures that have flown away''.
On the other hand (see \S\ref{SS:CP}) for any countable family of Borel
non-negative \fns\ (resp.\ bounded Borel non-negative \fns) on $\OM$ there
is a $\BAR{\OM}$ such that these \fns\ can be extended to l.s.c.\
(resp.\ \cts) \fns\ on $\BAR{\OM}$.
One may say that when the theorems of this \S\ are applied to an $\BAR{\OM}$
``also some flown-away measures are taken into account'', while the resulting
equidecomposability ``holds also for limit values''. They do not address the
case when only ($\sigma$-finite) measures in $\OM$ ``proper'' are considered.
\end{Rmk}

\subsection{Proofs} \label{SS:Proofs}

\begin{Prff} {\EM Proof of Thm.\ \ref{Thm:USC}:}

Introduce the following notation: If $\cA\subset\cC(\OM)^+$, then
$\tilde{\cA}$ is the set of all $h\in\cC(\OM)$
finitely equidecomposable via nonnegative \cts\ \fns\ with some
$f\in\cA$. For $A=\{f\}$ write $\tilde{f}:=\tilde{A}$.

\begin{Lem} \label{Lem:USC}
Let $\OM$ be a compact $G$-space. Let $\mu$ be a positive measure on
$\OM$. Then the functional $\tilde{\mu}$ on $\cC(\OM)$, corresponding to
every $f\in\cC(\OM)^+$ the infimum of $\mu$ on $\tilde{f}$
is additive, so it extends to a bounded positive functional on $\cC(\OM)$,
that is, to an (evidently invariant) finite positive regular measure
$\tilde{\mu}$.
\end{Lem}

\NOT{
\begin{Prff} {\EM Proof of Lemma \ref{Lem:USC}:}
Clearly $\tilde{\mu}(f+g)\le\tilde{\mu}(f)+\tilde{\mu}(g)$. The opposite
inequality follows from Remark \ref{Rmk:USC-LSC} item \ref{it:USC-LSCSum}
which implies that any member of $\tilde{}(f+g)$ 	
is $>$ $-\eps+$ some member of the sum of the sets $\tilde{f}+\tilde{g}$.
\end{Prff}
}

\begin{Prff} {\EM Proof of Lemma \ref{Lem:USC}:}
This follows from the fact that $\tilde{}(f+g)=
\tilde{f}+\tilde{g}$. The latter is a consequence of $\cC(\OM)$ being
an Abelian group lattice, hence has the {\em decomposition property}
(or {\em refinement property} -- see \cite{WehrungInj1}): 
If $f_i,g_j\ge0\:\:i=1,\ldots,n\:j=1,\ldots,m$ and $\sum_if_i=\sum_jg_j$
then $\exists h_{ij}\ge0$ s.t.\ $f_i=\sum_jh_{ij},\:g_j=\sum_ih_{ij}$ --
see \cite{BourbakiALG} Ch.\ VI \S1 Thm.\ 1. 
\end{Prff}

{\EM Continuation of the Proof of Thm.\ \ref{Thm:USC}:}
For our u.s.c. \fn\ $f$, denote by $\cH_f$ the set of \cts\ \fns\
$>f$. Define $\cH_g$ analogously.

A rephrasing of the theorem is (see Rmk.\ \ref{Rmk:USC-LSC}): if 
$\LA \nu,f\RA\ge\LA\nu,g\RA$ for any $G$-invariant p.m.\ $\nu$ on $\OM$, then
\NOT{$\tilde{\cH_f}\subset \tilde{\cH_g}$.}
$\cH_f\subset\tilde{\cH_g}$.
Recall that $\tilde{\cH_g}$ is the set of \cts\ \fns\ $>$ some \fn\
finitely equidecomposable with $g$ via u.s.c.\ \fns. 

It is easily proved that $\tilde{\cH_g}$ is convex. It is also open in
$\cC(\OM)$.
Indeed, if $h\in \cH_g$, then since $g-h$ is u.s.c.\ $<0$, $\max(g-h)<0$,
thus for some $\eps>0$\quad$h-\eps$ is still $>g$,
hence for any $h'\in\tilde{h}$,\quad$h'-\eps\in\tilde{\cH_g}$.

Thus one can apply separation theorems for open convex sets.
Suppose $f_0\in \cH_f$ but $f_0\notin\tilde{\cH_g}$. Then $\exists$ a bounded
functional $\mu\ne0$ on $\cC(\OM)$ with $\mu>\mu(f_0)$ on $\tilde{\cH_g}$.
Since $\tilde{\cH_g}+\cC(\OM)^+\subset\tilde{\cH_g}$, 
$\mu\ge0$ and one may assume $\mu$ a p.m. Consider $\tilde{\mu}$
(Lemma \ref{Lem:USC}). It is an invariant positive measure dominated by
$\mu$. By definition $\tilde{\mu}\ge\mu(f_0)$ on $\cH_g$, hence
$\tilde{\mu}(g)\ge\mu(f_0)$ (Remark \ref{Rmk:USC-LSC} items
\ref{it:USC-LSCInf} and \ref{it:USC-LSCInt}).
But we have $f_0\in \cH_f$, implying,
as we saw above, $f_0-\eps>f$ for some $\eps>0$, hence
$\tilde{\mu}(g)\ge\mu(f_0)>\mu(f)\ge\tilde{\mu}(f)$. This contradicts
the assumption of the theorem, that 
$\LA \nu,f\RA\ge\LA\nu,g\RA$ for any $G$-invariant p.m.\ $\nu$ on $\OM$,
{\em a fortiori} for any invariant finite positive $\nu$. 

\end{Prff}

\begin{Prff} {\EM Proof of Thm.\ \ref{Thm:LSC}}
We start with an analog of Thm.\ \ref{Thm:USC}.
Denote by $\cC_{00}(\OM)$ the vector lattice of real \cts\ \fns\ on $\OM$
with compact support. Endow $\cC_{00}(\OM)$ with the (Hausdorff)
locally convex topology which is the direct limit of the spaces,
for compact $K\subset\OM$,
$\cC(K)\cap\cC_{00}$ (normed by supremum norm on $K$). This is the strongest
locally convex topology making all the inclusion maps from these spaces \cts\
(see \cite{BourbakiEVT}). 
The positive cone of the dual space $\cC_{00}(\OM)^*$ is identified with the
set of Radon measures on $\OM$ (see \cite{BourbakiINT}). 

\begin{Lem} \label{Lem:LSC}
Let $G$ be a discrete group (possibly uncountable).
Let $\OM$ be a $G$-locally compact space.

Let $f:\OM\to\bR$ and $g:\OM\to\bR$ be nonnegative l.s.c.\
(as mentioned above they may take the value $+\I$).

Suppose

\begin{itemize}

\item
Any orbit of $G$ intersects the set $\{g>0\}$ and the set $\{f>0\}$.

\item
\NOT{If $U$ is the (open) union of the orbits that intersect $\{g>0\}$, then}%
Any (non-negative) $G$-invariant Radon measure on $\OM$ (possibly with
infinite mass) gives to $f$ an integral greater or equal than to $g$.
\end{itemize}

Then every \cts\ \fn\ with compact support which is strictly
dominated on its support by $g$
\NOT{
is strictly less than some u.s.c.\ \fn\ finitely decomposable with
$g$ via u.s.c.\ \fns, equivalently (see Remark \ref{Rmk:USC-LSC} item
\ref{it:USC-LSCSum}) every \cts\ \fn\ strictly greater than $f$
}%
is finitely equidecomposable via nonnegative \cts\ \fns\ with compact support
with some \fn\ strictly dominated on its support by $f$.
\end{Lem}

Let us first show how the theorem follows from Lemma \ref{Lem:LSC}:

\begin{itemize}
\item[(i)]
We may and do assume $U=\OM$ (note $U$ is open in a locally compact space
hence is itself locally compact). 

Denote by $\sim$ the relation: ``finitely equidecomposable via nonnegative
\cts\ \fns\ with compact support''.

Write $f=\sup f_n$ $g=\sup g_n$ where $f_n, n=1,2,\ldots$ and
$g_n, n=1,2,\ldots$ are nondecreasing sequences of \cts\ \fns\
with compact support, with $f_n<f$ on $\supp f_n$, $g_n<g$ on $\supp g_n$
(see Remark \ref{Rmk:USC-LSC} item \ref{it:BLSCCount}). 

By Lemma \ref{Lem:LSC}, $f_1\sim g_1'$, $g_1'$ being strictly dominated on
its support by $g$.
Replace $g_n,n\ge1$ by $\max(g_n,g_1')$. Again by Lemma \ref{Lem:LSC},
$g_1-g_1'\sim f_1'-f_1$, 
$f_1'$ being strictly dominated on its support by $f$.
Replace $f_n,n\ge2$ by $\max(f_n,f_1')$. Again, $f_2-f_1'\sim g_2'-g_1$,
$g_2'$ being strictly dominated on its support by $g$. Replace $g_n,n\ge2$ by
$\max(g_n,g_2')$. Now $g_2-g_2'\sim f_2'-f_2$, and we
continue, addressing alternatively $f$ and $g$, {\em ad infinitum}. 
We have
\BER{l}
f=f_1+(f_1'-f_1)+(f_2-f_1')+(f_2'-f_2)+\ldots\\
g=g_1'+(g_1-g_1')+(g_2'-g_1)+(g_2-g_2')+\ldots
\EER
and the terms are mutually
finitely equidecomposable via nonnegative
\cts\ \fns\ with compact support.
This means that they can be written as sums $\sum k_x$ and $\sum xk_x$,
$k_x$ \cts\ with compact support. Collecting all $k_x$ with the same $x$
for all the terms, we have our conclusion that
$f$ and $g$ are countably equidecomposable via nonnegative
b.l.s.c.\ \fns.

\item[(ii)]
The proof here follows the same lines as for item (i), but instead of
addressing $f$ and $g$ alternatively we go in one direction, writing
$g$ as a series of \fns\ $\sim$ to a series dominated by $f$. 
\end{itemize}

The {\EM Proof of Lemma \ref{Lem:LSC}} has similarities with the proof of
Thm.\ \ref{Thm:USC}. As we did there, introduce here the notation:
If $\cA\subset\cC_{00}(\OM)^+$, then
$\tilde{\cA}$ is the set of all $h\in\cC_{00}(\OM)^+$
finitely equidecomposable via nonnegative \cts\ \fns\ with compact support
with some $f\in\cA$. For $A=\{f\}$ write $\tilde{f}:=\tilde{A}$.

Let $\mu$ be a positive Radon measure on $\OM$. In analogy with
Lemma \ref{Lem:USC},
consider the functional $\tilde{\mu}$ on $\cC_{00}(\OM)^+$, corresponding to
every $f\in\cC_{00}^+$ the {\em supremum} of $\mu$ on $\tilde{f}$,
In analogy with the proof of that lemma, $\mu$ (which may take the value
$+\I$) is additive. So, if it is finite and bounded on bounded sets in
$\cC_{00}^+$ (these are sets of \fns\ bounded uniformly on every compact), 
it extends to a positive \cts\ functional on $\cC_{00}(\OM)$,
i.e.\ an (invariant) positive Radon measure $\tilde{\mu}$.

\NOT{\EM continuation of the Proof of Thm.\ \ref{Thm:USC}:}%
For our l.s.c.\ \fn\ $f$, denote by $\cM_f$ the set of non-negative \cts\
\fns\ with compact support, strictly dominated on their support by $f$.
Define $\cM_g$ analogously.

A rephrasing of Lemma \ref{Lem:LSC} is: if 
$\LA \nu,f\RA\ge\LA\nu,g\RA$ for any $G$-invariant Radon measure $\nu$ on
$\OM$, then
\NOT{$\tilde{\cH_f}\subset \tilde{\cH_g}$.}
$\cM_g\subset\tilde{\cM_f}-\cC_{00}^+$.
\NOT{
Recall that $\tilde{\cH_g}$ is the set of \cts\ \fns\ $>$ some \fn\
finitely equidecomposable with $g$ via u.s.c.\ \fns. 
}

It is easily proved that $\tilde{\cM_f}-\cC_{00}^+$ is convex.
Let us prove that it is
open in $\cC_{00}(\OM)$ (note that $\cM_f-\cC_{00}^+$ need not be open --
$0$ is not an interior point if $f$ is not always $>0$).
Indeed if $h\in \cM_f$, then $f-h$ is l.s.c., $>0$ on $\supp h$ hence on
$\{f>0\}$. We have $h+\tilde{\cM_{f-h}}\subset\tilde{\cM_f}$.
Thus to prove $\tilde{\cM_f}$ open it suffices to prove

\begin{Lem} \label{Lem:LSC0Int} 
For $f$ nonnegative l.s.c.\ on $\OM$, s.t.\ 
any orbit of $G$ intersects the set $\{f>0\}$,
the function $0$ is an interior point of $\tilde{\cM_f}-\cC_{00}^+$.
\end{Lem}

\begin{Prff}
{\EM Proof of Lemma \ref{Lem:LSC0Int}}
By the definition of the topology in $\cC_{00}(\OM)$, we need to prove that
$\forall$ compact $K\subset\OM$ $\exists\eps>0$ s.t.\ all \cts\ \fns\
$h$ with support $\subset K$ and $\max |h|<\eps$ belong to
$\tilde{\cM_f}-\cC_{00}^+$.

The assumption that any orbit intersects $\{f>0\}$ implies that the open sets
$x\cdot\{u>0\},\:x\in G, u\in\cM_f$ cover $K$. 
Choose a finite subcovering\ $x_j\cdot\{u_j>0\},\: j=1,\ldots N$ and let
$u(\om)=\frac{1}{N}\sum_j u_j(x_j^{-1}\om)$. Then
$u\in\tilde{\cM_f}$ and $u>0$ on $K$. Now take $\eps=\min_{\om\in K}u(\om)$.
\end{Prff}

This allows us to use separation theorems for open convex sets in the
locally convex space $\cC_{00}(\OM)$. Suppose $f_0\in\cM_g$ but
$f_0\notin\tilde{\cM_f}-\cC_{00}^+$. Then $\exists$ a \cts\ functional
$\mu\ne0$ on $\cC(\OM)$ with $\mu<\mu(f_0)$ on $\tilde{\cM_f}-\cC_{00}^+$.
\NOT{Since $\tilde{\cM_f}-\cC_{00}(\OM)^+\subset\tilde{\cM_f}$} This implies
$\mu\ge0$. We wish to consider $\tilde{\mu}$ (as defined above).
We can assert that it is finite and \cts, since by Lemma \ref{Lem:LSC0Int},
the fact that $\mu<\mu(f_0)$ on $\tilde{\cM_f}$ implies $\tilde{\mu}$
is bounded above on $\tilde{\cM_f}-\cC_{00}^+$, a \nbd\ of $0$. Thus
$\tilde{\mu}$ is an invariant (positive) Radon measure on $\OM$ and it is not
$0$ since $\mu\le\tilde{\mu}$.
By definition $\tilde{\mu}\le\mu(f_0)$ on $\cM_f$, hence
$\tilde{\mu}(f)\le\mu(f_0)$ ($f$ is the pointwise supremum of $\cM_f$ and
a Radon measure commutes with suprema of l.s.c.\ \fns\ --
see Remark \ref{Rmk:USC-LSC}).
We have $\tilde{\mu}(g-f_0)>0$. Indeed, $f_0\in \cM_g$, hence $g-f_0>0$ on
$\{g>0\}$. Thus by Lemma \ref{Lem:LSC0Int} $\tilde{\cM_{g-f_0}}-\cC_{00}^+$
is a \nbd\ of $0$, in which $\tilde{\mu}$
is dominated by $\tilde{\mu}(g-f_0)$, hence the latter is $>0$. Thus
$\tilde{\mu}(f)\le\mu(f_0)\le\tilde{\mu}(f_0)<\tilde{\mu}(g)$.
This contradicts the assumption of the lemma, that 
$\LA \nu,f\RA\ge\LA\nu,g\RA$ for any $G$-invariant Radon measure $\nu$ on
$\OM$. 
\end{Prff}

\subsection{Averages}
\label{SS:Av}
The simplest case of equidecomposability is that of a \fn\ and its
averages.

\begin{Def}
Given a $G$-vector space $V$ (always viewed as real), (that is, the group $G$
acts on $V$ linearly).

Let $v\in V$. An {\EM average} of $v$ is any element of $V$ of the form
$$\sum_{x\in G}\la_x xv\quad\la_x\ge0,\:\sum_{x\in G}\la_x=1,
\mbox{ all but finitely many }\la_x=0.$$

That is, a member of the convex hull of the orbit of $v$.
\end{Def}

One may ask: to what extent can we characterize, say, equality of integral
w.r.t.\ all invariant measures, allowing just averages instead of any
equidecomposable \fn? In what sense are equidecomposable \fns\
``similar'' already in averages?

We use the following theorem, which is an infinite-dimensional version of
the famous Von~Neumann Minimax Theorem in Game Theory (\cite{vN},
\cite{vNM}, a standard textbook is \cite{Owen}). For completeness, we give
a proof:

\begin{Thm} \label{Thm:VN}
Let $\OM$ be a compact space.

Let there be given a convex set $\cF$ of\NOT{real-valued} u.s.c.\ \fns\
on $\OM$.

Denote by $\cM\OM$ the set of all p.m.'s on $\OM$.

Then
$$\inf_{f\in\cF}\max_{\om\in\OM}f(\om)=\max_{\mu\in\cM\OM}\inf_{f\in\cF}
\LA\mu,f\RA.$$

\end{Thm}

\begin{Prf}

Clearly
$$\inf_{f\in\cF}\max_{\om\in\OM}f(\om)\ge\max_{\mu\in\cM\OM}\inf_{f\in\cF}
\LA\mu,f\RA.$$
Now let $a<\inf_{f\in\cF}\max_{\om\in\OM}f(\om)$ and we prove
$a\le\max_{\mu\in\cM\OM}\inf_{f\in\cF} \LA\mu,f\RA$.

W.l.o.g.\ we may assume $a=0$ (replace each $f$ by $f-a$). Let
$$\eps=\inf_{f\in\cF}\max_{\om\in\OM}f(\om)>a=0.$$
Let $\cH$ be the set of
all \cts\ \fns\ dominating some member of $\cF$. Then the \fn\
$0$ is distanced at least $\eps$ from $\cH$ in the $\sup$-norm.
Moreover, $\cH$ is convex
since $\cF$ is so. By convex separation in the Banach space $\cC(\OM)$,
$\exists$ a $\mu_0\in \cC(\OM)^*,\:\mu\ne0$ positive on $\cH$. Since for
nonnegative \cts\ $f$ \quad $\cH+f\subset\cH$, $\mu_0$ must be
nonnegative on nonnegative \cts\ \fns. Hence we may and do assume
$\mu_0\in\cM\OM$. By construction $\mu_0$ is nonnegative on any \cts\
\fn\ dominating some $f\in\cF$.
By Remark \ref{Rmk:USC-LSC}
we have $\mu_0$ nonnegative on any
$f\in\cF$, i.e.\
$\inf_{f\in\cF} \LA\mu_0,f\RA\ge0$, implying
$0\le\max_{\mu\in\cM\OM}\inf_{f\in\cF} \LA\mu,f\RA$.
\end{Prf}

In case $\OM$ is a compact convex space, every
p.m.\ $\mu$ on $\OM$ has a barycenter $\om$, and the integral over $\mu$ of a
\cts\ (or u.s.c.) {\em affine} \fn\ $f$ on $\OM$ is $f(\om)$
(see \cite{Phelps}). Thus \ref{Thm:VN} takes the form:

\begin{Cor} \label{Cor:VNCnv}
Let $\OM$ be a compact convex space.

Let there be given a convex set $\cF$ of\NOT{real-valued} {\em affine}
u.s.c.\ \fns\ on $\OM$.

Then
$$\inf_{f\in\cF}\max_{\om\in\OM}f(\om)=\max_{\om\in\OM}\inf_{f\in\cF}
f(\om).$$

\end{Cor}

Now let us return to averages in a $G$-vector space. Suppose
a sublinear functional $p:V\to\bR$ {\em invariant under $G$} is given.
By sublinear is meant, as usual, that $p$ satisfies:
\BER{l}
p(u+v)\le p(u)+p(v)\quad u,v\in V\\
p(\la v)=\la p(v)\quad v\in V,\:\la>0
\EER

An example is: a $G$-compact space $\OM$, $V$ an invariant subspace of
$\cC(\OM)$, $p(f)=\max_{\om\in\OM}f(\om)$.

In fact, the general case of invariant sublinear
functionals reduces to the above example.

Indeed, Suppose $p$ is invariant sublinear as above. Let $\OM$ be the set
of all linear functionals $\om:V\to\bR$ satisfying $\om\le p$, with the
weak $V$-topology. Since
$$\om\le p\Rightarrow \om(v)=-\om(-v)\ge-p(-v)$$
$\OM$ is a $G$-compact space, and obviously any $v\in V$ induces a \cts\
affine \fn\ $v$ on $\OM$. By Hahn-Banach $p(v)=\max_{\om\in\OM}v(\om)$.

\begin{Thm} \label{Thm:MaxAv}
Given a $G$-vector space $V$. Let
$p:V\to\bR$ be sublinear and $G$-invariant.

Let $v\in V$, and denote by $\cA\cV(v)$ the set of its averages
(i.e.\ the convex hull of its orbit).

Denote by $\OM$ the set of all linear functionals $\om:V\to\bR$ satisfying
$\om\le p$.

Then

\begin{itemize}
\item[a.]
Suppose $G$ {\em amenable}. Then the infimum of $p$
over $\cA\cV(v)$ is equal to the maximum of $\LA\om,v\RA$ for all
the $G$-invariant $\om\in\OM$.

\item[b.]
For general $G$, the infimum of $p$ over $\cA\cV(v)$ is equal to the
maximum over $\om\in\OM$ of $\inf_{x\in G}\LA x\om,v\RA=
\inf_{x\in G}\LA \om,xv\RA$, i.e.\ to the maximum of the infimum of
$v$ on nonempty (convex) $G$-invariant subsets of $\OM$.

\end{itemize}
\end{Thm}

\begin{Prf}
Note first, that the fact that we indeed get maxima (and not just suprema)
follows from upper-semi-continuity of the maximized \fn\ on the compact
$\OM$.

In case $G$ is amenable, any compact convex invariant subset of $\OM$
contains an invariant point. Hence a.\ follows from b.

To prove b., identify $V$ with a set of affine \cts\ \fns\ on
the compact $\OM$. Then $p(v)=\max v$, and use Cor.\ \ref{Cor:VNCnv}:
\par\medskip
The infimum of ($p=\max$) over the convex $\cA\cV(v)$ =
\par
The maximum over $\om\in\OM$ of $\inf_{w\in\cA\cV(v)} w(\om)$ =
\par
The maximum over $\om\in\OM$ of $\inf_{x\in G} (xv)(\om)$
\par\medskip
and we are done.
\end{Prf}

\begin{Cor} \label{Cor:MaxAv1}
In case $G$ is {\em amenable} (retaining the assumptions of the theorem):

Two vectors giving the same value to all invariant $\om\in\OM$, in
particular two vectors finitely equidecomposable (via elements of $V$),
have the same infimum of $p$ over the set of their averages.
\end{Cor}

Cor.\ \ref{Cor:MaxAv1} says that, in the amenable case, the property of
two vectors to be finitely equidecomposable can be judged merely by their
averages. (In case $V$ is a $G$-normed space, this extends to one vector
approximable by vectors finitely equidecomposable with another).

Cf.\ Example \ref{Exm:NonAme} in \S\ref{SS:ERG}.

\begin{Rmk}
Compare Thm.\ \ref{Thm:MaxAv} a., for the case of $\cC(\OM)$ and $p=\max$,
with Cor.\ \ref{Cor:USC1}. In the amenable case, \fns\ finitely
equidecomposable with a given \fn\ are replaced by its averages,
on which there is much greater control. 
\end{Rmk}

Recall that Cor.\ \ref{Cor:USC1} led to Cor.\ \ref{Cor:USC-Boo}
for $G$-Boolean algebras. For $G$ amenable, Thm.\ \ref{Thm:MaxAv} gives:

\begin{Cor} \label{Cor:MaxAv-Boo}
Let $G$ be {\em amenable} and let $\cB$ be a $G$-Boolean algebra.

Let $E\in\cB$. If a constant $a$ is greater
than the supremum of the masses of $E$ w.r.t.\ all $G$-invariant
finitely additive p.m.\ on $\cB$ then $\exists$ a $p\in\bN$ and
$p$ $G$-translates of $E$, covering every part of $1$ no more than
$\lfloor ap\rfloor$ times.
\end{Cor}

Thus $E$ is in {\EM ``generalized Rokhlin position''}, where an $E\in\cB$
is in {\EM ``Rokhlin position''} if
$\exists$ a finite set $\{x_1,x_2,\ldots,x_N\}\subset G$ s.t.\
$x_iE,\:i=1,\ldots,N$ are disjoint. Then we are sure that for every
invariant finitely additive p.m.\ on $\cB$, the mass of $E$ $\le$ $1/N$.
Similarly, If every part is covered by 
no more than $k$ of the $x_iE,\:i=1,\ldots,N$ then we
know that for every such invariant p.m.\ $\mu$\quad$\mu(E)\le k/N$.

Cor.\ \ref{Cor:MaxAv-Boo} is a kind of converse statement.

\begin{Exm} \label{Exm:MaxAv-Boo1}
In Cor.\ \ref{Cor:MaxAv-Boo} (hence in Cor.\ \ref{Cor:MaxAv1})
the assumption that $G$ is amenable is essential. Just let $G$ be
non-amenable (discrete) and take $\cB=\cP(G)$, $G$ acting on itself by
left translations, $E=1(=G)\in\cB$. Note that the supremum of masses of
invariant finitely-additive p.m.'s is $0$, since no such p.m.\ exists.
\end{Exm}

\begin{Exm} \label{Exm:MaxAv-Boo2}
Cor.\ \ref{Cor:MaxAv-Boo} with $a$ {\em equal} to the supremum of the masses
(and Thm.\ \ref{Thm:MaxAv} a.\ with ``infimum''
over $\cA\cV(v)$ replaced by ``minimum'') may not hold.
As a counterexample take $G=\bZ$ acting on itself by translations,
$\cB=\cP(\bZ)$, $E=2\bZ\cup\{1\}$. 
\end{Exm}

\begin{Exm} \label{Exm:MaxAv-Boo3}
In Cor.\ \ref{Cor:MaxAv-Boo} $p$ may not be chosen ``in advance'': Let
$G=C_n^2$, $C_n$ being cyclic of order $n$. Let $G$ act on itself by
translations and let $\cB=\cP(G)$, $E=(C_n\times\{e\})\cup(\{e\}\times C_n)$.
Although the unique $G$-invariant p.m.\ gives to $E$ a measure $<2/n$,
$E$ is not in ``Rokhlin position'' even for $N=2$, since it intersects each
of its translates.
\end{Exm}

See \S\ref{SS:ERG} for a treatment of ``mean ergodic theorems'' for general
(discrete) groups, where ``convergence'' means: ``having averages arbitrarily
near the limit''. 

\newpage
\section{Equidecomposable Enhanced Functions}
\label{S:EDCT}

\subsection{Introduction and Inquiry}
\label{SS:EDCTIntr}
\NOT{Similarly to what was said in the beginning of}
In a similar manner to what was said in Section \ref{S:ED},
in the
continuous case too some notions are independent of the invariant measure
$\mu$ in an $\OM$ acted by $G$.
Such are $m$-chains and their vertices, in particular enhanced functions.
Two enhanced functions will be called {\EM equidecomposable}, with suitable
qualifications (such as, in the discrete case: finitely-, countably-, via
\cts\ \fns\ etc.) if a suitable $1$-chain exists, s.t.\ the two
enhanced \fns\ are, resp., its source and target. This notion too
does not depend on the invariant measure $\mu$ in $\OM$.

Note that $1$-chains obtained from weighted hypergraphs are not sufficient:
even to have one enhanced \fn\ equidecomposable with itself one needs
a $1$-chain supported on the diagonal of $G^2$. 

Thus, in applications of CHG, such as those in Section \ref{S:CT},
in particular \S\ref{SS:CtAsso}, the two enhanced \fns\ are, in fact,
{\em equidecomposable} (via measurable $1$-chains).

CVE tells us that equidecomposable enhanced \fns\ have the same expectation,
in other words if two functions, enhanced by $\La_1$ and $\La_2$ resp.,
are equidecomposable, then {\em for any invariant probability measure $\mu$
in $\OM$}, their integrals w.r.t.\ the corresponding infinitesimal
measures are the same. 

Thus the relation of two \fns\ to be equidecomposable after enhancement by
$\La_1$ and $\La_2$ appears as {\em dual} to the relation of two measures
to be infinitesimal measures obtained from the same invariant probability
$\mu$ by $\La_1$ and $\La_2$. One may try to take this duality as a
{\em defining property} of one of these two notions starting from the other.

To this end, one may try to prove analogs of the facts in Section \ref{S:ED}
for {\em enhanced functions}. It may be better to work in topological
$G$-spaces, and to restrict oneself to enhanced \fns\ which give u.s.c.\ or
l.s.c.\ $0$-chains, and also to equidecomposability understood as
being the source and target of a u.s.c.\ or l.s.c.\ $1$-chain. One may\NOT{
try to} define such attributes of chains by requiring that applying them to
every nonnegative \cts\ with compact support {\em test-\fn} be u.s.c.\
or l.s.c.\ In this setting, does one has an analog to Thm.\ \ref{Thm:LSC}?

Given a Borel-measurable $m$-chain (this means, as usual, that applying
it to a test-\fn\ on $G^{m+1}$ depends on $\om$ in a Borel manner).
Choose a countable collection $\cF$ of
test-\fns\ (i.e.\ non-negative \cts\ with compact support) on $G^{m+1}$
s.t.\ {\em every test-\fn\ is a non-decreasing limit of a sequence of members
of $\cF$} and s.t.\ {\em every member of $\cF$ is a positive combination of
convolutions of two test-\fns} -- see Rmk.\ \ref{Rmk:CntTstFns}.
\NOT{
(Take as $\cF$, e.g., all positive rational
finite combinations of the union of non-decreasing sequences of test-\fns,
each of them a convolution of two test-\fns, that converge to the
characteristic \fns\ of the members of a countable base to the topology in
$G^{m+1}$.)
}%
By \S\ref{SS:CP} we may assume $\OM$ is a dense Borel subset
of a $G$-metrizable compact space $K$ s.t.\ all applications of the $m$-chain
to members of $\cF$, hence to all test-\fns, extend to l.s.c.\ \fns\ on $K$.
(The application to a convolution of two test-\fn\ is a convolution of
a non-negative Borel \fn\ on $\OM$ and a test-\fn\ (i.e.\ $L^1$-\fn) on $G$,
and by \S\ref{SS:CP} any countable collection of such extend to l.s.c.\
for a suitable $K$.)
Thus the $m$-chain is l.s.c.\ on $\OM$ with the relative Lusin topology
(see \S\ref{SS:PRD}). Can the chain be extended to the compact $G$-space $K$ ?
Then it would be extendable to a b.l.s.c.\ $m$-chain on the 
canonical non-metrizable compactification $\cK$ in \S\ref{SS:CP}.
Can this extension be made unique?
Can this be used for the purposes mentioned above?

\NOT{
Referring to the compactification $\cK$ in \S\ref{SS:CP}, it seems that every
(nonnegative) measurable $m$-chain extends to a b.l.s.c.\ chain in $\cK$.
(Every convolution of a nonnegative measurable \fn\ on $\OM$ with a test-\fn\
on $G$ extends to a b.l.s.c.\ on $\cK$; hence for every nonnegative
measurable chain,
applying it to a test-\fn\ that is a convolution of two test-\fns\ is
b.l.s.c., and one may try to pass to a limit.) Can this extension be
made unique? Can this be used for the above purposes?
}

Before dealing with such questions, there is the question of the
{\em transitivity}
of the (qualified) equidecomposability relation, even for the discrete case.
This can be put as follows: suppose two $1$-chains have a common vertex
(that is, a common projection on $G$);
are the other two vertices the two vertices of some $1$-chain? 
or, if we are lucky, does there exist a $2$-chain having the given
$1$-chains as ``sides'', i.e.\ projections on $G^2$?
In the discrete case such a $2$-chain can be constructed, using some
canonical construction of a nonnegative matrix with given row- and column-
sums. But qualifications such as u.s.c.\ may be violated. Another approach
to the transitivity is via criteria for equidecomposability such as
Thm. \ref{Thm:LSC}. \S\ref{SS:Tame} below refers to the question of
transitivity.

Peculiar to the continuous case are question such as: given $\La_1$ and
$\La_2$ and $f_1$. Does there exist an $f_2$ s.t.\ $f_1$ enhanced by
$\La_1$ is equidecomposable with $f_2$ enhanced by $\La_2$?, and a related
question: can one reconstruct $\mu$ from the infinitesimal measure
$\frac{d\La}{d\la}\mu$? This is dealt with in \S\ref{SS:REC} below.%
\NOT{Here formulas such as (\ref{eq:NoDs2}) or
(\ref{eq:Ctsfs}), (\ref{eq:CtsfSs}) might be of help
(see also Rmk.\ \ref{Rmk:CtsFlow}).}

\subsection{Tame Invariant Chains and Transitivity of Equidecomposability}
\label{SS:Tame}

Our setting is a 2nd-countable locally compact group $G$ acting in a Borel
manner on a standard Borel space $\OM$.
In this \S, {\em measurable} will mean either ``Borel measurable''
or ``universally measurable''.
We use the notions about chains (in the ``\cts'' case)
mentioned in \S\ref{SS:CtsChn}.
For simplicity, we restrict ourselves to the case that
$G$ is {\em unimodular}, thus speak of invariant chains.
For non-unimodular $G$, replace ``invariant'' by
``$\Delta$-right-invariant'' and ``Haar measure'' by ``left Haar measure''.

\begin{Def} \label{Def:Tame}
A (non-negative) invariant ($\om$-dependent) $m$-chain is called
(Borel measurably, resp.\ universally measurably)
{\EM tame} if it can be written as a countable sum of non-negative
measurable invariant chains all of which are Radon measures (on $G^{m+1}$)
for each $\om\in\OM$ (i.e.\ are finite on compacta) -- the latter will be
referred to as (measurable, invariant) {\EM Radon chains}.

If two $0$-chains are source and target of the same {\em tame} $1$-chain, we
say that they are {\EM Borel-tamely}- resp.\ {\EM universally measurably
tamely} {\EM equidecomposable}.
\end{Def}

Check that in most of the examples in Section \ref{S:CT} the chains are tame.

A vertex, and more generally a face, of a tame $m$-chain is tame. This follows
from the fact that an invariant measurable Radon chain $\phi$ is a countable
sum of invariant measurable Radon chain with Radon projections -- just
write $d\phi=\sum_n h_n(x_0x_1^{-1},\ldots,x_0x_m^{-1})\,d\phi$ where the
$h_n$ are test-\fns\ on $G^m$ with sum $\equiv1$.

\begin{Rmk} \label{Rmk:RdnRN}
We shall use the fact that if $\nu_n$ is a sequence of Radon measures in a
2nd-countable locally compact $X$, then there exists a sequence $\al_n>0$
s.t.\ $\sum_n\al_n\nu_n$ is a Radon measure -- just take an increasing
sequence of compact sets $K_n\st X$ that eventually contains any fixed
compact $\st X$, and choose $\al_n>0$ s.t.\ $\sum_n\al_n\nu_n(K_n)<\I$.
Note that all the $\nu_n$ are absolutely \cts\ w.r.t.\ $\sum_n\al_n\nu_n$,
with bounded Radon-Nikodym derivatives.
\end{Rmk}

\begin{Lem} \label{Lem:RdnChnLtc}
Let $\phi_1$ and $\phi_2$ be measurable invariant ($\om$-dependent)
Radon chains. Then the (invariant) chain
$\phi_1\land\phi_2$ mapping each $\om\in\OM$ to the infimum of $\phi_1(\om)$
and $\phi_2(\om)$ in the lattice of non-negative Radon measures is also
measurable.
\end{Lem}

\begin{Prf}
The lemma follows from $\LA\phi_1(\om)\land\phi_2(\om),h\RA$
($h$ -- a {\em test-\fn} on $G^{m+1}$,
i.e. a non-negative \cts\ \fn\ with compact support)
being equal to the infimum of $\LA\phi_1,h_1\RA+\LA\phi_2,h_2\RA$ over the
pairs $\{(h_1,h_2):h_1,h_2\in\cF,\:h_1+h_2\ge h\}$, $\cF$ being a
countable collection of test-\fns, s.t.\ every test-\fn\ is the
limit of a non-decreasing sequence of members of $\cF$
(see Rmk.\ \ref{Rmk:CntTstFns}).
\NOT{
(Take as $\cF$, e.g., all positive rational finite combinations of the union
of non-decreasing sequences of test-\fns that converge to the characteristic
\fns\ of the members of a countable base to the topology in $G^{m+1}$.)
}
\end{Prf}

\begin{Lem} \label{Lem:RdnChnPrd}
Let $\phi_0$, $\phi_1$ and $\phi_2$ be measurable invariant ($\om$-dependent)
Radon chains, s.t.\ for each $\om\in\OM$ $\phi_1$ and $\phi_2$ are absolutely
\cts\ w.r.t.\ $\phi_0$.
Then the (invariant) chain $\phi_3$, s.t.\ for each $\om$
\begin{equation} \label{eq:RdnChnPrd}
d\phi_3:=\dfrac{d\phi_1}{d\phi_0}\dfrac{d\phi_2}{d\phi_0}\,d\phi_0
\end{equation}
is also measurable.
\end{Lem}

\begin{Prf}
If we had in (\ref{eq:RdnChnPrd}) $\land$, i.e.\ $\min$ instead of
multiplication of the Radon-Nikodym derivatives, the assertion would follow
from Lemma \ref{Lem:RdnChnLtc}. Therefore we would be done if we show 
how to express multiplication in $\BAR{\bR^+}$ using $\land$ and
``well-behaved'' limits (i.e.\ which preserve the measurability of the chain).
This is done by ($a,b\in\BAR{\bR^+}$):
$$ab=\lim_n2^{-n}\sum_{j,k\ge1}
\LQ(a\land(j\cdot2^{-n}))-(a\land((j-1)\cdot2^{-n}))\RQ\land
\LQ(b\land(k\cdot2^{-n}))-(b\land((k-1)\cdot2^{-n}))\RQ$$
the limit being uniform in $a,b$.

(To convince oneself of the validity of this formula, note that for all except
one of the summands, one of the expressions in $[\:]$ is $2^{-n}$,
the other being some $t\in[0,2^{-n}]$, and then 
$2^{-n}(2^{-n}\land t)=2^{-n}t$.)
\end{Prf}

\begin{Lem} \label{Lem:RdnChn2Seq}
Let $(\phi_i)_{i\ge1}$ and $(\phi'_j)_{j\ge1}$ be two sequences of
measurable invariant Radon $m$-chains s.t.\
$\sum_i\phi_i=\sum_j\phi'_j$. 
Then one can find measurable invariant Radon $m$-chains
$\phi_{ij}$ s.t.\ $\phi_i=\sum_j\phi_{ij}$ and $\phi'_j=\sum_i\phi_{ij}$.
\end{Lem}

\begin{Prf}
For $i\ge0$, $j\ge0$ write $\psi_i:=\sum_{1\le\ell\le i}\phi_\ell$,
$\psi'_j:=\sum_{1\le\ell\le j}\phi'_\ell$.
Let $\Psi:=\sum_i\phi_i=\sum_j\phi'_j$.
Then $\psi_i\uparrow\Psi$, $\psi'_j\uparrow\Psi$, where we have this
limit relation for the application of the measures to any $[0,\I]$-valued
Borel \fn\ on $G^{m+1}$. Construct a matrix of chains $(\phi_{ij})_{i,j\ge1}$
s.t.\ (cf.\ Lemma \ref{Lem:RdnChnLtc})
\begin{equation} \label{eq:RdnChn2Seq}
\sum_{1\le\ell\le i,1\le\ell'\le j}\phi_{\ell\ell'}=\psi_i\land\psi_j
\end{equation}
in other words,
$$\phi_{ij}=
\psi_i\land\psi'_j+\psi_{i-1}\land\psi'_{j-1}-
\psi_i\land\psi'_{j-1}-\psi_{i-1}\land\psi'_j$$
The measurability of $\psi_{ij}$ follows from Lemma \ref{Lem:RdnChnLtc}
provided we prove that the $\psi_{ij}$ are {\em non-negative}, i.e.\ that
$$\psi_i\land\psi'_{j-1}+\psi_{i-1}\land\psi'_j\le
\psi_i\land\psi'_j+\psi_{i-1}\land\psi'_{j-1}.$$
That follows from the following chain of assertions:
\BER{l}
\psi_i\land\psi'_{j-1}+\psi_{i-1}\land\psi'_j\le\psi_i+\psi_{i-1}\\
\psi_i\land\psi'_{j-1}+\psi_{i-1}\land\psi'_j\le\psi_i+\psi'_{j-1}\\
\psi_i\land\psi'_{j-1}+\psi_{i-1}\land\psi'_j\le\psi'_j+\psi_{i-1}\\
\psi_i\land\psi'_{j-1}+\psi_{i-1}\land\psi'_j\le\psi'_j+\psi'_{j-1}\\
\psi_i\land\psi'_{j-1}+\psi_{i-1}\land\psi'_j\le
\psi_i+\psi_{i-1}\land\psi'_{j-1}\\
\psi_i\land\psi'_{j-1}+\psi_{i-1}\land\psi'_j\le
\psi'_j+\psi_{i-1}\land\psi'_{j-1}\\
\psi_i\land\psi'_{j-1}+\psi_{i-1}\land\psi'_j\le
\psi_i\land\psi'_j+\psi_{i-1}\land\psi'_{j-1}
\EER
Being dominated by Radon chains, $\phi_{ij}$ are Radon chains.

It remains to prove the claims $\sum_j\phi_{ij}=\phi_i$ and
$\sum_i\phi_{ij}=\phi'_j$. By (\ref{eq:RdnChn2Seq}) these are equivalent to
$\lim_j\psi_i\land\psi'_j=\psi_i$, $\lim_i\psi_i\land\psi'_j=\psi'_j$.
Since $\lim_i\psi_i$ and $\lim_i\psi'_j$ are $\Psi$, and
$\psi_i\land\Psi=\psi_i$, $\Psi\land\psi'_j=\psi'_j$, all we have to show
is that in our case $\lim$ commutes with $\land$. This can be shown
for each $\om$ separately.
By Rmk.\ \ref{Rmk:RdnRN}, for any fixed $\om$ there is a Radon measure
s.t.\ all the $\phi_i$ and $\phi'_j$ are absolutely \cts\ w.r.t.\ it,
with bounded Radon-Nikodym derivatives.
Thus to prove $\lim$ commutes with $\land$ we can pass to the Radon-Nikodym
derivatives, which are \fns, and for them this is immediate.
\end{Prf}

\NOT{
\begin{Lem}
Let $(\phi_i)_{i\ge1}$ and $(\phi'_j)_{j\ge1}$ be two sequences of
measurable invariant Radon $m$- and $m'$-chains resp. s.t.\ a certain
$k$-face ($0\le k\le m,m'$) of the sum $\sum_i\phi_i$ is equal
to a certain $k$-face of the sum $\sum_j\phi'_j$. 
Then one can find measurable invariant Radon $m$-chains
$\phi_{ij}$ and $m'$-chains $\phi'_{ij}$, s.t.\ $\phi_i=\sum_j\phi_{ij}$,
$\phi'_j=\sum_i\phi'_{ij}$ and the corresponding $k$-faces in $\phi_{ij}$
and $\phi'_{ij}$ are equal.
\end{Lem}
}

\begin{Lem} \label{Lem:DecomFace}
Let $(\phi_i)_{i\ge1}$ be a sequences of measurable invariant Radon
$m$-chains s.t.\ a certain $k$-dimensional face ($0\le k\le m$) of the sum
$\sum_i\phi_i$ can be written as $\sum_j\phi'_j$ for some measurable
invariant Radon $k$-chains $(\phi'_j)_{j\ge1}$.
Then one can find measurable invariant Radon $k$-chains
$\phi_{ij}$ s.t.\ $\phi_i=\sum_j\phi_{ij}$ and $\phi'_j$ is the relevant
$k$-face of $\sum_i\phi_{ij}$.

Consequently, if a $k$-face of a tame $m$-chain is the sum of a sequence
of measurable invariant Radon $k$-chains, these can be written as the
$k$-faces of tame $m$-chains that sum to the given $m$-chain.
\end{Lem}

\begin{Prf}
Let $\tilde{\phi}_i$ be the relevant face of $\phi_i$, and by further
decomposing the $\phi_i$ one may assume $\tilde{\phi}_i$ are Radon.
Apply Lemma \ref{Lem:RdnChn2Seq} to the sequences $\tilde{\phi}_i$ and
$\phi'_j$, which have the same sum, to obtain measurable invariant
Radon $k$-chains $\phi'_{ij}$ s.t.\ $\phi'_j=\sum_i\phi'_{ij}$,
$\tilde{\phi}_i=\sum_j\phi'_{ij}$.

We have to write the $\phi'_{ij}$ as $k$-faces (i.e.\ projections)
of measurable invariant Radon $m$-chains $\phi_{ij}$ s.t.\ 
$\phi_i=\sum_j\phi_{ij}$. We do this as follows (fix $i,j$): $\phi'_{ij}$ is
dominated by $\tilde{\phi_i}$. Take the Radon-Nikodym derivative $r_{ij}$ of
$\phi'_{ij}$ w.r.t.\ $\tilde{\phi_i}$, which is a \fn\ on $G^{k+1}$,
expand it to a \fn\ on $G^{m+1}$ by composing it with the relevant projection,
and multiply it by $d\phi_i$ to obtain $\phi_{ij}$.

The only thing that needs further proof is that $\phi_{ij}$ is measurable.
Denote by $x'$ the variable in $G^{k+1}$ and by $(x,x')$ the variable in
$G^{m+1}$, $(x,x')\mapsto x'$ being the relevant projection (i.e.\ the
relevant $k$-face of the $m$-simplex $(x,x')$ is the $k$-simplex $x'$).
We have to prove that applying $\phi_{ij}$ to a fixed test-\fn\ $h(x,x')$
is measurable in $\om$. Since every test-\fn\ is a non-decreasing limit of
finite positive combinations of \fns\ of the form $f(x)g(x')$,
$f,g$ test-\fns\ (Rmk.\ \ref{Rmk:CntTstFns}), we may
assume $h(x,x')=f(x)g(x')$ is of that form. 
Now, for such test-\fn, applying $\phi_{ij}$ to it is the same as
applying $\phi_i$ and replacing $g$ by $g$ multiplied by the Radon-Nikodym
derivative $r_{ij}$. Thus if we keep $f$ fixed and concentrate on the
dependence on $g$, then the transition from $\phi_i$ applied to $f(x)g(x')$
to $\phi_{ij}$ applied to the same is given by the lattice operations in the
proof of Lemma \ref{Lem:RdnChn2Seq}, which preserve measurability
by Lemma \ref{Lem:RdnChnLtc}.
\end{Prf}

\begin{Thm} \label{Thm:TransitiveTame}
Let a 2nd-countable locally compact group act in a Borel manner on a
standard Borel space. Then if two (Borel- resp.\
universally-measurably) tame $1$-chains (Def.\ \ref{Def:Tame})
have a common vertex, then they are sides of the same (Borel- resp.\
universally-measurably) tame $2$-chain.

Consequently, the relation between (Borel- resp.\
universally-measurable) tame $0$-chains to be (Borel- resp.\
universally-measurably) tamely equidecomposable is transitive.
\end{Thm}

\begin{Prf}
Decomposing the (tame) common vertex into a sum of measurable invariant Radon
$0$-chains and applying Lemma \ref{Lem:DecomFace}, we may assume all the
chains are (measurable invariant) Radon.

Let $\phi_1$ (on $1$-simplices $(x_0,x_1)$) and
$\phi_2$ (on $1$-simplices $(x_0,x_2)$) be the two $1$-chains.
They have a common vertex $\phi_0$ (on $0$-simplices $(x_0)$).
For each fixed $\om\in\OM$, we construct the $2$-chain $\phi_3$ as follows: 
disintegrate the two $1$-chain w.r.t.\ the projections on $G$ defined by the
common vertex (about disintegration of measures see \cite{BourbakiINT}).
This gives families
$(\nu^{(1)}_{x_0})_{x_0\in G}$, $(\nu^{(2)}_{x_0})_{x_0\in G}$
of probability measures on $G$ s.t.\ for non-negative Borel $f(x_0,x_i)$
($i=1,2$) $x_0\mapsto\int_G f(x_0,x_i)\,d\nu^{(i)}_{x_0}(x_i)$ is Borel and
$$\phi_1=\int_G\delta_{x_0}\otimes\nu^{(1)}_{x_0}\,d\phi_0(x_0),\quad
\phi_2=\int_G\delta_{x_0}\otimes\nu^{(2)}_{x_0}\,d\phi_0(x_0)$$
(the integration of measures is defined, as usual, by applying test-\fns).
These $\nu^{(1)}_{x_0}$ and $\nu^{(2)}_{x_0}$ are determined by $\phi_1$
and $\phi_2$ (for fixed $\om$) up to change of the $\nu^{(i)}$'s in a
$\phi_0$-null set.

Now define, for each fixed $\om$:
$$\phi_3:=\int_G\delta_{x_0}\otimes\nu^{(1)}_{x_0}\otimes\nu^{(2)}_{x_0}
\,d\phi_0(x_0)$$
as a measure on $2$-simplices $(x_0,x_1,x_2)\in G^3$. It is standard to check
that this indeed defines a Radon measure, with projections $\phi_1$ and
$\phi_2$ and that the dependence on $\om$ defines an invariant $2$-chain.
It remains to prove that the chain $\phi_3$ is measurable.

As in the proof of Lemma \ref{Lem:DecomFace}, it suffices to check test-\fns\
of the form $h_0(x_0)h_1(x_1)h_2(x_2)$ where the $h_i$ are test-\fns.
Fix $h_1$ and $h_2$ and perform the integration on the $\nu^{(i)}$'s.
Then it is clear that, for the dependence on $h_0$, we have the situation of
Lemma \ref{Lem:RdnChnPrd} for
$\LA\phi_0,h_0\RA$, $\LA\phi_1,h_0\otimes h_1\RA$,
$\LA\phi_2,h_0\otimes h_2\RA$ and $\LA\phi_3,h_0\otimes h_1\otimes h_2\RA$.
Therefore applying that lemma gives the measurability of $\phi_3$.
\end{Prf}

\begin{Def} \label{Def:TameFn}
Let $\La$ be some right-invariant measure on $G$. A \fn\
$f:\OM\to\BAR{\bR^+}$ is called (Borel- resp.\ universally- (measurably))
$\frac{d\La}{d}$-{\EM tame} if $\forall\om\in\OM$
$\cO_\om f$ is $\La$-measurable, and the $0$-chain enhancement 
$f\frac{d\La}{d}$ is tame (Def.\ \ref{Def:Tame}).
$\frac{d\La}{d}$-{\EM Radon} \fns\ are defined analogously.
\end{Def}

Every non-negative Borel \fn\ is Borel-measurably $\frac{d\la}{d}$-tame for
$\la$ a Haar measure.

We also have:
if $f$ is a Borel \fn\ s.t.\ the set
$\{x\in G:\cO_\om f(x)\ne0\}$ has countable closure for each $\om\in\OM$,
then $f$ is Borel-measurably $\frac{d\cnt}{d}$-tame.
This follows from the following
``invariant'' way to decompose a countable close set $K$ in $G$ to
countably many discrete sets (on which, of course, the counting measure
is Radon): choose a right-invariant metric in $G$. Decompose the set of
isolated points $x\in K$ to the countably many discrete sets
$\{x\in K:2^{-(k+1)}<\dist(x,K\sm\{x\})\le2^{-k}\}$.
Then do the same to all derivatives of $K$.
(Note that there is a countable ordinal $\gamma$ s.t.\ the
$\gamma$-derivative of the (closed) support of $\cO_\om f$ is $\es$ for every
$\om\in\OM$. This follows from Prop.\ \ref{Prop:EGDelta} and from
the considerations in the footnote in \S\ref{SS:MEAS}).

Also, let an invariant probability measure $\mu$ be given in $\OM$
(and $\la$ a Haar measure on $G$). Then any 
$\frac{d\La}{d\la}\mu$-integrable (Borel) \fn\ $f$ can be changed on a Borel
$\frac{d\La}{d\la}\mu$-null set to become 
a Borel-measurably $\frac{d\La}{d}$-Radon (Borel) \fn $f'$.

To see this, note that, in the context of \S\ref{SS:EnhInf}, $f\in\Pre$, thus
the $0$-chain $f\frac{d\La}{d}$ is measurable and has a Radon expectation.
This implies that for any test-\fn\ $h$ on $G$, the set
$$\OM_h:=\{\om\in\OM: \int_Gf(x\om)h(x)\,d\La(x)<\I\}$$
is ($\mu$-)conull.
 
Choose a countable set $\cF$ of test-\fns\ on $G$ s.t.\ every test-\fn\ $h$
is majorized by some $h'\in\cF$.
Then the intersection $\OM'$ of all the $\OM_h$ for all test-\fns\ $h$ is
equal to $\cap_{h\in\cF}\OM_h$, hence is ($\mu$-)conull. Moreover, $\OM'$ is
$G$-invariant. At an $\om\in\OM'$, $\frac{d\La}{d}f$,
i.e.\ $f(x\om)\,d\La(x)$, is a Radon measure.
Pick a $G$-invariant {\em Borel} ($\mu$-)conull $\OM''\st\OM'$.
(Prop.\ \ref{Prop:BorelAlmInv}).
Replacing $f$ by $f'=f\cdot1_{\OM''}$ will do.

\subsection{Recovery of the Original Measure from the Infinitesimal
Measure}
\label{SS:REC}
In this \S\ we deal with two related problems mentioned
in\NOT{in item \ref{it:TwoQu} in the inquiry}
\S\ref{SS:EDCTIntr}: to recover the measure from the infinitesimal
measure, and to find a \fn\ which, enhanced by one invariant measure,
will be equidecomposable to a given \fn\ enhanced by another invariant
measure.

Our setting is a unimodular 2nd countable locally compact group $G$ 
(with Haar measure $\la$) acting in a Borel manner on a standard Borel
space $\OM$.
  
Occasionally, $\OM$ will be endowed with a (probability) measure $\mu$
such that the action of $G$ is measure-preserving, in other words,
$\mu$ is invariant.

Also, we shall sometimes have a right-invariant measure $\La$
on $G$. 

\begin{Rmk}
Suppose $G_0$ is an open subgroup of the acting group $G$. Consider the
restriction $\la_0:=\la|_{G_0}$ which is a Haar measure in $G_0$.
Then it is clear that for any invariant measure $\La$ in $G$,
with restriction $\La_0:=\La|_{G_0}$, enhancing a \fn\ or measure on $\OM$
by $\frac{d\La}{d\la}$ is the same as enhancing it by
$\frac{d\La_0}{d\la_0}$, referring to $G_0$ as the acting group.
Moreover, any $1$-chain for $G_0$ can be viewed as a $1$-chain for $G$,
hence two enhanced \fns\ equidecomposable w.r.t.\ $G_0$ are 
{\em ipso facto} equidecomposable w.r.t.\ $G$.

This implies that solving our above problems for $G_0$ will solve them for
$G$. That applies, in particular, to $G$ -- a Lie group and $G_0$ -- its
identity component. Thus we may assume our Lie groups are {\em connected}.

\end{Rmk}

\subsubsection{Reduction to the Counting Measure Case}
\label{SS:RdcCnt}

Note the two formulas in \S\ref{SS:NoDs}, which are, in fact, instances
of equidecomposability of enhanced \fns\ as in CHG, namely:

\begin{itemize}
\item
Let $f:\OM\to\BAR{\bR^+}$ be\NOT{Borel and}
Borel-measurably $\frac{d\La}{d}$-tame.
(Similar consideration will hold for universally-measurably tame \fns.)

Let $h:G\to\BAR{\bR^+}$ be non-negative Borel.

\NOT{
Suppose $\La$ has the universal measurability Fubini property
(Def.\ \ref{Def:MeasFub}). Then the $0$-chain $f\frac{d\La}{d}$ is
universally measurable.
}

Consider the weighted graph $F(x,y;\om)=h(xy^{-1})f(x\om)$ and the
$1$-simplex of measures $(\La,\la)$. That defines a Borel-measurably tame
$1$-chain\NOT{(note the $\sigma$-finiteness requirement of CHG holds)}
(consider $(x,y)\mapsto h(xy^{-1})$ as a test-\fn\ on $G^2$, and approach it
by sums of test-\fns\ with separated variables $x$ and $y$),
thus its source and target are equidecomposable. This means that
$$\LQ\int_G h(x^{-1})\,d\la(x)\RQ f(\om)\frac{d\La}{d\la}$$
is Borel-tamely equidecomposable to the \fn\
$$\om\mapsto\int_G h(x)f(x\om)\,d\La(x).$$
The latter \fn\ is Borel (consider $h$ as a test-\fn).

Consequently, every Borel-measurably $\frac{d\La}{d}$-tame \fn\
$f:\OM\to\BAR{\bR^+}$ is, enhanced by $\frac{d\La}{d\la}$, 
Borel-tamely equidecomposable with a Borel \fn.

If one takes $\La=\la$, one finds that any non-negative Borel \fn\ is
Borel-tamely equidecomposable to a \fn\ which is a ``convolution'' of a Borel
\fn\ on $\OM$ and an $\L^1$-\fn\ on $G$, hence is l.s.c.\ (lower semi-\cts)
for some Lusin topology on $\OM$ (see \S\ref{SS:CP}).

\item
Let $f:\OM\to\BAR{\bR^+}$ be%
\NOT{Borel with $\cO_\om f$ always zero outside a countable set (depending on
$\om$), in particular,} 
Borel-measurably $\frac{d\cnt}{d}$-tame.

Let $h:G\to\BAR{\bR^+}$ be Borel s.t.\ $h(x^{-1})$ is $\La$-integrable.

Consider the weighted graph $f(x\om)h(xy^{-1})$ and the $1$-simplex of
measures $(\cnt_G,\La)$.
Note the discussion of this graph in \S\ref{SS:NoDs}. It is noted there
that the $\sigma$-finiteness requirements of CHG hold for all $\om$, and
that discussion shows that the corresponding $1$-chain is Borel-tame.

\NOT{
with $\La$ having the universal measurability Fubini property
(Def.\ \ref{Def:MeasFub}).
Then\NOT{By Rmk.\ \ref{Rmk:CntMeasFub},} the corresponding $1$-chain is
universally measurably tame, and its source and target are equidecomposable.
}

Thus we have
$$\LQ\int_G h(x^{-1})\,d\La(x)\RQ f(\om)\frac{d\cnt}{d\la}$$
Borel-tamely equidecomposable with 
$$g(\om)\frac{d\La}{d\la}\quad\mbox{for }
g(\om):=\int_G f(x\om)h(x)\,d\cnt_G(x).$$
$g$ is Borel\NOT{, by Rmk.\ \ref{Rmk:CntMeasFub}.}
(consider $h$ as a test-\fn).

We conclude that for any invariant measure 
$\La$ that is not $\{0,\I\}$-valued (thus $\exists$ a $\La$-integrable
$h\ge 0$ with $\int h\,d\La\ne0$),
and for any non-negative Borel-measurably $\frac{d\cnt}{d}$-tame \fn\ $f$ on
$\OM$,%
\NOT{with $\cO_\om f$ always zero outside a countable set,}
$f$, enhanced by $\frac{d\cnt}{d\la}$,
is Borel-tamely equidecomposable to some Borel \fn\ enhanced by
$\frac{d\La}{d\la}$.

Therefore for any invariant $\mu$, $\frac{d\cnt}{d\la}\mu$ can be
reconstructed from $\frac{d\La}{d\la}\mu$ (Recall that any
$\frac{d\cnt}{d\la}\mu$-integrable Borel \fn\ is equal, outside a Borel
$\frac{d\cnt}{d\la}\mu$-null set, to a Borel-measurably
$\frac{d\cnt}{d}$-tame Borel \fn. -- see the end of \S\ref{SS:Tame})
\end{itemize}
 
To conclude, if one can recover $\mu$ from
$\frac{d\cnt}{d\la}\mu$ one can recover $\mu$ from any $\frac{d\La}{d\la}\mu$,
if $\La$ is not $\{0,\I\}$-valued; Borel-tame equidecomposability
allows us to pass from \fns\ enhanced by such $\La$ to usual \fns\
(even to \fns\ l.s.c.\ for a suitable Lusin topology)
and from \fns\ enhanced by $\cnt_G$ to \fns\ enhanced by such $\La$.

So, it remains to try to recover $\mu$ from $\frac{d\cnt}{d\la}\mu$, and to
try to find \fns\ which, enhanced by $\frac{d\cnt}{d\la}$ will be 
equidecomposable with given ordinary \fns, which may be assumed l.s.c.\ for
some Lusin topology.

\subsubsection{Finding $E$ with Discrete $\cO_\om E$ and $\spn E$ the Whole
Space -- Analogy to Ambrose-Kakutani}
 
From now on we restrict ourself to $G$ a unimodular {\em connected Lie group}.

As noted in Rmk.\ \ref{Rmk:CtsFlow}, in the case of
$\bR$-action (i.e.\ a flow) one may use the method of Ambrose and Kakutani
\cite{AmbroseKakutani} (see also \cite{Jacobs}, \cite{Nadkarni}),
to recover the original
measure from the infinitesimal measure (in fact, from the restriction of
the infinitesimal measure to suitable $E\st\OM$ s.t.\ the original
system has the structure of a ``flow under a \fn'').
This can be pursued for more general groups $G$ rather than $\bR$,
and one may try not to
refer to a particular invariant measure $\mu$, thus speaking about
G-Borel spaces. This is done in \cite{Kechris} and \cite{FeldmanHahnMoore}
(cf.\ also \cite{Wagh}).
Our approach will be differential-geometric, and
seems different from theirs (cf.\ \cite{Ramsay}).

Ambrose-kakutani teach us to look for Borel sets $E\st\OM$ with
$\cO_\om E$ {\em discrete} in $G$ and $\spn E$ conull, preferably
$\spn E=\OM$

One case when such $E$ does
not exists is when {\em the stabilizer of some $\om\in\OM$} (i.e.\ the
subgroup $\{x\in G: x\om=\om\}\st G$) {\em is a non-discrete closed subgroup
(i.e.\ of positive dimension)} -- 
as in Exm.\ \ref{Exm:InfMea} items \ref{it:InfMeaR2T} and \ref{it:InfMeaSO3}.
Note that we may (and do) assume that $\OM$ is a $G$-compact metric space
-- see \S\ref{SS:CP}, thus the stabilizer is always a {\em closed} subgroup.
In such cases one cannot expect to recover $\mu$ from
$\frac{d\cnt}{d\la}\mu$.

Take for example $G=SO(3)$ (with normalized Haar measure $\la$)
acting on the union of 
two concentric spheres in $R^3$ by rotations (compare Exm.\ \ref{Exm:InfMea}
item \ref{it:InfMeaSO3}). Taking as $\mu$ any non-trivial convex combination
of the normalized invariant areas on the spheres, one gets the same
$\frac{d\cnt}{d\la}\mu=\I\cdot\cnt$.

Note that since the stabilizer of $y\om$ is a conjugate of the stabilizer of
$\om$, the set of all $\om$ with non-discrete stabilizer is invariant.
Moreover, in the case that $\OM$ is a $G$-compact space this set is
{\em closed} (see \cite{Ramsay}).
Indeed, it is the projection on $\OM$ of the compact set
$$\{(\om,v)\in\OM\times T_eG: \|v\|=1,\:\forall 0\le t\le1\:
\exp(tv)\om=\om\}$$
where some norm on the tangent space $T_eG$ is understood. 

\NOT{Borel (Assume $\OM$ a $G$-compact metric space.
Denote appropriate metrics in $G$ and in $\OM$ by $d$.
Fix a sequence of countable open coverings of $G$ with mesh
($=$ maximum diameter of a member) $\to0$.
The stabilizer of $\om$ is discrete iff for some member $\cU$ of the
sequence, for any $U\in\cU$, $x\mapsto x\om$ is 1-1 on $U$. This will
happen iff for some fixed sequence $K_n$ of compacts with $\cup K_n=U$, 
$\forall\eps>0\forall n\exists\delta$ s.t.\
$x,y\in K_n,\:d(x,y)>\delta\Rightarrow d(x\om,y\om)>\eps$.)
}

But what if the stabilizers are discrete? we have the following analog of
Ambrose-Kakutani:

\begin{Thm} \label{Thm:LieAK}
(see \cite{Kechris} and \cite{FeldmanHahnMoore})
Let $G$ be a Lie group acting in a Borel manner on a standard Borel space
$\OM$.

Partition $\OM$ into two invariant Borel sets: the set $\OM_1$ of
$\om\in\OM$ with non-discrete stabilizer and its complement -- the set
$\OM_0$ of $\om\in\OM$ with discrete stabilizer.
Then the latter set contains a Borel $E\st\OM$ with $\cO_\om E$ discrete
for all $\om\in\OM$, and with $\spn E=\OM_0$.
\end{Thm}

\begin{Prf}

\NOT{We may and do assume $\OM=\OM_0$.}

\NOT{Construct a transfinite sequence $E_\al$,
indexed by the countable ordinals $\al$, s.t.\ $E_\al$ are Borel $\st\OM$
with $\cO_\om E_\al$ discrete for all $\om\in\OM$, and $\spn E_\al$ have
positive measure and are disjoint, as follows: 

Suppose the sequence had been defined for $\al<\al_0$. Delete
$\cup_{\al<\al_0}\spn E_\al$ from $\OM$. If the complement contains some
Borel $E$ with $\cO_\om E$ discrete for all $\om\in\OM$ and $\mu(\spn E)>0$,
define $E_{\al_0}:=E$. Otherwise the process terminates.

Now, the process must terminate, otherwise $\OM$ will contain an uncountable
disjoint family $\spn E_\al$, $\al$ runs over all countable ordinals, with
$\mu(\spn E_\al)>0$. When the process terminate, it is evident that 
$E=\cup E_\al$ will satisfy the requirements of the theorem, provided
$\cup\spn E_\al$ is conull. Thus (replacing $\OM$ by the invariant
$\OM\sm\cup\spn E_\al$ with relative measure, provided the latter
is not null) it remains to show the following: 

{\em If each $\om\in\OM$ has a discrete stabilizer, then $\exists$ $E\st\OM$
with $\cO_\om E$ discrete for all $\om\in\OM$ s.t.\ $\spn E$ has
positive measure}.
}

W.l.o.g.\ assume $\OM$ is a $G$-compact metric space (see \S\ref{SS:CP}).

Let $n$ be the dimension of the Lie group $G$, and choose a right Haar
measure $\la$ in $G$.

In the set of all closed subsets of $G$ or of all open subsets of $G$
we take the Effros Borel structure (see \S\ref{SS:MEAS}), obtained by
identifying each closed $F\st G$ (resp.\ each open $U\st G$) with the set of 
members $u$ of a fixed countable open base to the topology that satisfy 
$u\cap F=\es$ (resp.\ $u\st U$), this set being a member of $2^\bN$.

Our strategy will be to correspond to each $\om\in\OM_0$ a non-empty discrete
set $F(\om)\st G$ in a Borel and {\em equivariant} manner, the latter meaning
that $F(\om)=F(a\om)a,\:a\in G$. Then the set $E=\{\om:e\in F(\om)\}$
has the property
$$a\om\in E\Leftrightarrow e\in F(a\om)\Leftrightarrow a\in F(\om),
\quad a\in G$$
Thus $\cO_\om E$ is $F(\om)$, which is discrete non-empty. In particular,
$\spn E=\OM_0$.

Call a \fn\ $f:\OM\to\bR$ $\cC^\I$ if $f$ is \cts, 
$\forall\om\in\OM$ $\cO_\om f$ is
$\cC^\I$ on $G$ and $\om\mapsto\cO_\om f$ is \cts\ from $\OM$ to the
Fr\'echet space $\cC^\I(G)$ with the usual topology of uniform convergence
on compacta of all partial derivatives. For any \cts\ $f:\om\to\bR$ and
any $\cC^\I$-\fn\ with compact support $h:G\to\bR$, the convolution
$$\om\mapsto\int_G f(x\om)h(x)$$
is $\cC^\I$. Thus any \cts\ \fn\ on $\OM$ can be uniformly approximated by
$\cC^\I$ \fns.

We wish to be able to bring together differential geometric notions
pertaining to different points in $G$. We do this by equating the
tangent spaces via {\em right translations}. Namely, choose a
basis to the tangent space $T_eG$, and ``expand'' it to $n$ right-invariant
vector fields $X_i,\:i=1,\ldots,n$ forming a basis to the tangent space
at each point. We use the notation $Xf$, where $f$ is a $\cC^\I$ \fn\
on $G$ and $X$ is a tangent vector or vector field on $G$, understood as
acting on $f$ as usual (see \cite{Hicks}).

For any $\al\in\bR^n$ there is a unique differential form (which is $\cC^\I$)
$\tilde\al$ on $G$ satisfying 
$\LA\tilde\al,X_i\RA=\al_i,\:i=1,\ldots,n$
These $\tilde\al$ are the {\EM right Maurer-Cartan forms on $G$}
(see \cite{Cohn} \S4.4).

Let $f:\OM\to\bR$ be $\cC^\I$, let $\al\in\bR^n$ and consider the differential
form on $G$, depending on $\om\in\OM$, $d\cO_\om f-\tilde\al$. (Note that
this, as a member of the Fr\'echet space of the differential forms on $G$,
depends \cts{ly} on $(\om,\al)$,
thus there is no question about Borelness of the subsets of $\OM$ to be
considered below.)

Suppose we find an {\em invariant} (i.e.\ constant on orbits) Borel \fn\
$\om\in E'\mapsto\al(\om)\in\bR^n$
where $E'$ is a Borel invariant set $E'\st\OM$, s.t.\
$\forall\om\in E'$ the (closed) set $F_0(\om)\st G$ consisting of the points
where $d\cO_\om f-\tilde{\al(\om)}$ vanishes
(i.e.\ gives the zero element of the cotangent space)
contains some isolated points. Then the mapping $\om\in E'\mapsto F_0(\om)$
is Borel, and the right-invariant way in which $\tilde\al$ was
constructed (and the invariance of $\om\mapsto\al(\om)$) implies
$F_0(\om)=F_0(a\om)a\:a\in G$. To get a discrete non-empty
set $F(\om)$ out of $F_0(\om)$ in an equivariant way, just take:
$$F(\om):=\LB x\in G:\dist\LP x,F_0(\om)\setminus\{x\}\RP>
\frac12 \min\LP 1,\sup_{y\in F_0(\om)}\dist
\LP y,F_0(\om)\setminus\{y\}\RP\RP\RB$$ 

\NOT{
Let us speak of points $x\in F(\om)$ as {\EM far} points of $F_0(\om)$.

Thus, if $E\st\OM$ is the set
$$E:=\{\om\in E': e \mbox{ is a far point in the set where }
d\cO_\om f-\tilde{\al(\om)}
\mbox{ vanishes }\}$$
Then $E$ is Borel, and the right-invariant way in which $\tilde\al$ was
constructed (and the invariance of $\om\mapsto\al(\om)$) implies 
$$a\om\in E\Leftrightarrow d\cO_\om f-\tilde{\al(\om)}
\mbox{ vanishes at }a,\quad a\in G$$
Thus $\cO_\om E$ is just the vanishing set of $d\cO_\om f-\tilde{\al(\om)}$,
which is discrete, and clearly $\spn E=E'$.
}

Hence we shall be done if we can cover $\OM$ with a countable union
of invariant Borel sets $E'$ with
invariant (i.e.\ constant on orbits) Borel \fn\
$\om\in E'\mapsto\al(\om)\in\bR^n$ s.t.\ $\forall\om\in E'$ the set where
$d\cO_\om f-\tilde{\al(\om)}$ vanishes contains some isolated point
(reduce this covering to a disjoint countable covering and take the
union of the corresponding $E$).

Using the basis $X_i$, the differential form $d\cO_\om f-\tilde\al$
is described as a mapping $G\to\bR^n$, namely
\begin{equation} \label{eq:DiffMap}
x\mapsto\LP\LA d\cO_\om f-\tilde\al,X_i\RA\RP_i=
\LP X_i\cO_\om f-\al_i\RP_i
\end{equation}

This mapping has an invertible differential at some point of $G$ iff
its Jacobian matrix, which is, in fact, the Hessian of $f$ 
\begin{equation} \label{eq:Hessian}
H_{ij}(x;\om;f)=\LP X_j(X_i\cO_\om f)\RP_{ij}\quad x\in G,\om\in\OM 
\end{equation}
is non-singular at this point. If, for some $\om$ and $\al$, this Jacobian
is non-singular at some $a\in G$ where $d\cO_\om f-\tilde\al$ vanishes,
then at every such point $a$
(\ref{eq:DiffMap}) is locally 1-1, by the inverse function theorem,
hence every such $a$ is an isolated zero. 

Take as our $E'=E'(f)$ above the set of all $\om$ with
$H_{ij}(x;\om;f)$ (\ref{eq:Hessian}) non-singular somewhere on $G$.
$E'$ is open, its complement being the closed set of the $\om$ with
$H_{ij}$ singular everywhere in $G$.

In order to carry out our plan, we have to correspond to every $\om\in E'$,
in a Borel and invariant (i.e.\ constant on orbits) manner,
an $\al\in\bR^n$ s.t.\ $d\cO_\om f-\tilde\al$ vanishes somewhere
where $H_{ij}$ is non-singular.
This means that $\al$ belongs to the image $I(\om;f)\st\bR^n$ by
$$x\mapsto D(x):=\LP X_i\cO_\om f\RP_i$$
of the set of $x\in G$ where $H_{ij}(x)$ is non-singular.

But by the way $H_{ij}$ and $D(x)$ were defined, replacing $\om$ by $a\om$
($a\in G$) would change $H_{ij}(x)$ into $H_{ij}(xa)$ and $D(x)$ into
$D(xa)$. Therefore $\om\mapsto I(\om;f)$ is invariant
(i.e.\ constant on orbits). $I(\om;f)$ is open in $G$, and depends on
$\om$ in a Borel manner. Hence we can take as $\al(\om)$ the first element
belonging to $I(\om;f)$ in a fixed dense sequence in $\bR^n$. 

\NOT{
Note that $I_1(\om;f)$ is open in $G$, moreover
$$\tilde I_1(f):=\{(\om,\al):\al\in I_1(\om;f)\}\st\OM\times\bR^n$$
is open. Thus $\om\mapsto (I_1(\om))^c$ is Borel to the the Borel space of
closed subsets of $\bR^n$ (Hausdorff topology).
For $\om\in E'$, $I_1(\om;f)$ is nonempty.

By Sard's Lemma (see \cite{Schwartz} Ch.\ I) the image $I_0(\om;f)$
of the set of singular points has Lebesgue measure zero in $\bR^n$,
hence has empty interior. Being $F_\sigma$, it is meager. Also
$$\tilde I_0(f):=\{(\om,\al):\al\in I_0(\om;f)\}\st\OM\times\bR^n$$
is $F_\sigma$: one may write $I_1(\om)=\cup_{n\ge1}I'_k(\om)$ where
$I'_n(\om)\st G$ is an increasing sequence of closed sets, where
$\{(\om,\al):\al\in I'_n(\om)\}$ are closed
(take as $I'_n$ the image by $x\mapsto D(x)$ of the set of singular
points in a compact $K_n\st G$, with $K_n\uparrow G$).
So $\om\mapsto I'_n(\om)$ is Borel.

Baire category insures that $\forall\om\in E'$ $I_1\sm I_0$ is nonempty.
The way to choose in a Borel manner an $\al$ in it for each $\om\in E'$
(moreover, in an invariant, i.e.\ constant on orbits, way)
is given in the following lemma, whose proof is carried out by following a
constructive proof of the Baire Category Theorem (with some countability
assumptions):

\begin{Lem}
Let $X$ be a 2nd-countable locally compact space. Let $2^X$ be the space of
all closed subsets $\st X$ with the usual Borel structure.
Consider the countable product $X^{\{0,1,\ldots\}}$.
Then there is an $X$-valued Borel \fn\ on the Borel subset $\cS$ of the
product consisting of the sequences $\tilde{F}=(F_0,F_1,\ldots)$ with
$F_0\ne X$ and $F_1,F_2,\ldots$ nowhere dense, corresponding to each such
sequence a point in $F_0^c\sm\cup_{n\ge1}F_n$.
\end{Lem}

\begin{Prf}
Let $X\cup\{\I\}$ be the one-point (Alexandrov) compactification.
Fix a metric $d'$ giving the topology in $X\cup\{\I\}$. In $X$ take
the metric $d(x,y):=d'(x,y)+|1/d'(x,\I)-1/d'(y,\I)|$, which is
a complete metric inducing the topology in $X$.
Note that $\dist(x,F):X\times(2^X\sm\{\es\})\to\bR^+$ is Borel.

Fix a dense sequence in $X$.

Construct two sequences of \fns\ on $\cS$ $x_n(\tilde{F})\in X$,
$\eps_n(\tilde{F})>0$ as follows:
$x_n(\tilde{F})$ is the first member of the chosen dense sequence distanced
less than $\frac12\eps_j(\tilde{F})$ from $x_j(\tilde{F})$, $1\le j\le n-1$,
belonging to $F_0^c$ and not belonging to $F_j$, $1\le j\le n-1$.
$\eps_n(\tilde{F})$ is
$$\frac12\min\LP\eps_{n-1}(\tilde{F}),
\dist\LP x_n(\tilde{F}),\cup_{1\le j\le n}F_j\cup F_0\RP\RP.$$

This insures that the closed balls with centers $x_n$ and radii $\eps_n$ are
a decreasing sequence with the properties required by the usual proof of
Baire category, so that these closed balls intersect in a single point
$x(\tilde{F})$ and the $x_n(\tilde{F})$ converge to the Borel \fn\
$x(\tilde{F})$ with $x(\tilde{F})\in F_0^c\sm\cup_{n\ge1}F_n$.
\end{Prf}
}

Thus we can fulfill the requirements of the theorem for $E'=E'(f)$, that is,
find a Borel $E$ with discrete $\cO_\om E$ and $\spn E=E'$. Since $\OM$ is
compact metric, the union of the open $E'(f)$ for all $\cC^\I$-\fns\ $f$
is covered by a countable number of them, so the requirements of the theorem
hold also for this union. To conclude, we prove that this union is all of
$\OM_0$ -- the set of $\om$ with discrete stabilizer. (Certainly the union
is contained in $\OM_0$, since an $E$ as above cannot exist for
$\om\notin\OM_0$.)

So suppose $\om$ belongs to the complement of $E'(f)$ for all
$\cC^\I$-\fns\ $f$.
By the above, this means that for any such $f$, the Hessian
$H_{ij}(x;\om;f)$ (\ref{eq:Hessian}) is singular $\forall x\in G$.

Let $\cM$ be the set of matrix values attained by
$\LP H_{ij}(e;\om;f)\RP_{ij}$ at the unity $e$,
{\em for all $\cC^\I$-\fns\ $f$ on $\OM$}. $\cM$ is a vector subspace of
the vector space of $n\times n$ matrices. Since $f$ may be replaced by
$f(a\om)$ and $H_{ij}$ is defined in a right-equivariant manner, $\cM$ is
also the set of matrix values obtained by $H_{ij}(x;\om;f)$ for any
fixed $x$.

Now the Hessian $H_{ij}$ (\ref{eq:Hessian})
is symmetric, since (see \cite{Hicks}):
$$H_{ij}-H_{ji}= X_j(X_i\cO_\om f)-X_i(X_j\cO_\om f)=
X_j\LA d\cO_\om f,X_i\RA-X_i\LA d\cO_\om f,X_j\RA= 
\LA dd\cO_\om f,[X_j,X_i]\RA=0$$
In this situation, one has:

\begin{Lem} \label{Lem:Flanders}
(\cite{Flanders} Lemma 1.)
Let $\cM$ be a linear subspace of the space of all $n\times n$ symmetric
matrices over an infinite field $\bF$ with characteristic $\ne2$, with all
members of $\cM$ singular (i.e.\ having zero determinant).
Then there exists a non-zero vector $v\in\bF^n$ with
$v^TMv=0\:\:\forall M\in\cM$.
\end{Lem}

\begin{Prf}
Let $r$ be the maximum rank of members of $\cM$. By the assumption $r\le n-1$.
Let $A\in\cM$ be of rank $r$. Replacing, if necessary, every member $M$ of
$\cM$ by $Z^TMZ$, where $Z$ is a fixed non-singular matrix, we may assume
that $A$ is diagonal, with first $r$ diagonal entries $\ne0$ and the others
$0$. We prove that
$$\forall M\in\cM\:\forall i>r,j>r\:M_{ij}=0$$
Indeed, consider the linear family $M+tA$ of members of $\cM$, hence having
rank $\le r$. Consider in them the $(r+1)$-minor built from the first $r$
rows and columns, the $i$-th row and the $j$-th column. Its determinant,
which is an $r$ -th degree polynomial in $t$, is $0$ for every $t$. So its
coefficient of $t^r$, which is $A_{11}\cdots A_{rr}M_{ij}=0$, implying
$M_{ij}=0$.

In particular, $M_{nn}=0$ $\forall M\in\cM$, and the conclusion of the
lemma holds with $v=(0,\ldots,0,1)$.
\end{Prf}

So, by Lemma \ref{Lem:Flanders}, $\exists$ a $v\in\bR^n$, $v$ not the zero
vector, s.t.\ for our $\om$, $\forall$ $\cC^\I$-\fn\ $f:\OM\to\bR$,
$$\sum_{ij}H_{ij}(x;\om;f)v_iv_j=0$$
for all $x\in G$.
By (\ref{eq:Hessian}) this means:
\begin{equation} 
(\sum v_jX_j)\LQ(\sum v_iX_i)\cO_\om f\RQ=0\quad x\in G,\om\in\OM 
\end{equation}
This implies that if $t\mapsto x(t)=\exp(t\sum v_iX_i)$
is the integral curve of the
right-invariant non-zero vector field $\sum v_iX_i$ tracing\NOT{a right
coset of} the one-parameter subgroup corresponding to this vector, then
$\cO_\om f(x(t))$ has zero second derivative, hence is linear. Since $f$ is
bounded (\cts\ on a compact) $\cO_\om f(x)$, i.e.\ $f(x\om)$, is constant
on that one-parameter subgroup\NOT{any such right coset}. But we have noted
that $\cC^\I$-\fns\ are dense in $\cC(\OM)$. One concludes that $x\om=\om$
for every $x$ in the above one-parameter subgroup, i.e.\ the stabilizer of
$\om$ has positive dimension, thus $\om\notin\OM_0$.

This concludes the proof.
\end{Prf}

\subsubsection{Using $E$ with $\cO_\om E$ Non-Empty Countable Closed}
\label{SS:CntCls}

Having Thm.\ \ref{Thm:LieAK} at hand, one can address the problem
of recovery of $\mu$ from $\frac{d\cnt}{d\la}\mu$ and that of finding a \fn\
which, enhanced by $\frac{d\cnt}{d\la}$ will be equidecomposable to a given
usual \fn\ (which would solve the former problem). 

Thus assume $G$ is a unimodular connected Lie group, and that $E\st\OM$ is
Borel s.t.\ for all $\om\in\OM$ $\cO_\om E$ is a {\em non-empty countable
closed} subset of $G$.

We will use the ``nearest point'' construction of \S\ref{SS:Nearest},
in particular Prop.\ \ref{Prop:Nearest} (iv). One may note that,
at least for the recovery of $\mu$, we can avoid this construction by
applying CHG, instead of to the graph (\ref{eq:Nearest})
in \S\ref{SS:Nearest}
as below, to the Graph $\{(x,y)\in G^2: y\in\cO_\om E\}$. Here one uses
the fact that $\OM$ is a countable-to one Borel image of $E\times G$ by
$(\om\in E,x\in G)\mapsto x\om$.

We use the same notations as in \S\ref{SS:Nearest}. In particular,
for $\om\in\OM$, $\pi(\om)=\pi_E(\om)\in G$ is the unique nearest point
to $0$ in $\cO_\om E$, if it exists
(if there is no unique nearest point, $\pi_E(\om)$ is undefined);
for $\om\in E$, $P(\om)=P_E(\om)$ is the set of $x\in G$ s.t.\ $0$ is
the unique nearest point to $x$ in $\cO_\om E$.

By \S\ref{SS:Nearest}, the mapping $(\om,x)\mapsto x\om$ maps the set
$$\{(\om,x)\in E\times G:x\in P_E(\om)\}$$
in a (Borel) 1-1 manner onto the set
$$\OM'=\{\om\in\OM:\pi_E(\om)\mbox{ is defined}\}$$
i.e.
$$\OM'=\{\om\in\OM:\mbox{ there is a unique nearest point to }0\mbox{ in }
\cO_\om E\}.$$

Being a Borel 1-1 image of a Borel set, $\OM'$ is Borel. Also, $\OM'$ has
the property that $\cO_\om\OM'$ is (Haar-)conull for all $\om\in\OM$
(Prop.\ \ref{Prop:Geodesic} in \S\ref{SS:Nearest}). Thus, as far as a
Borel \fn\ on $\OM$ is considered as a chain (``enhanced'' by
$\frac{d\la}{d\la}$) it does not matter if we restrict it to $\OM'$.

Also, this 1-1 mapping $(\om,x)\mapsto x\om$ is a Borel isomorphism,
the inverse mapping $\OM'\to\{(\om,x)\in E\times G:x\in P_E(\om)\}$ being:
$$\om\mapsto\LP\pi(\om)\om,\pi(\om)^{-1}\RP.$$

Let $f:\OM\to\BAR{\bR^+}$ and $s:G\to\BAR{\bR^+}$ be Borel. Consider
the graph used to prove Prop.\ \ref{Prop:Nearest} (iii) with $s\equiv1$,
namely, the graph
$$\{(x,y)\in G^2: y \mbox{ is the unique nearest point to }x
\mbox{ in }\cO_\om E\}$$
weighted by $f(x\om)$, with the simplex of measures
$(\la,\cnt_G)$. As we have seen in \S\ref{SS:Nearest}, the obtained
$1$-chain is Borel-measurable. Moreover, it is Borel-tame -- see the
discussion following Def.\ \ref{Def:TameFn}.
The vertices, being hence Borel-tamely equidecomposable, are:

Source: the ordinary Borel \fn\ (restricted to $\OM'$, which, by the above,
does not matter):
$$\om\mapsto f(\om)$$

Target: the Borel function supported on $E$, enhanced by $\frac{d\cnt}{d\la}$:
\begin{equation} \label{eq:fPE}
\LQ\om\mapsto\int_{P_E(\om)}f(x\om)\,d\la(x)\RQ\,\frac{d\cnt}{d\la}
\end{equation}

\NOT{
So we conclude that for any \fn\ $g$ on $\OM$ s.t.\ $g(x\om)$ on
$\{(\om,x)\in E\times\OM:x\in P_E(\om)\}$ is of the form of a tensor:
$$(\om,x)\mapsto f(\om)s(x)$$ with $f$, $s$ Borel,
one can find a \fn\ (on $E$) Borel-tamely equidecomposable to it
when enhanced by $\frac{d\cnt}{d\la}$. This holds, of course, for any sum
of a series of such non-negative \fns, hence (\S\ref{SS:PRD} Prop.\
\ref{Prop:ProdFn}) for any \fn\ on $E\times\OM$ which is l.s.c.\ in a
product of two topologies in the two factors. Such will be $g(x\om)$, if
$g$ is a non-negative Borel \fn\ on $\OM$ which is l.s.c.\
for some Lusin topology on $\OM$ (with the given Borel structure) that
makes the action of $G$ \cts\ in the two variables. By what was said in
\S\ref{SS:REC}, every non-negative Borel \fn\ on $\OM$ is Borel-tamely
equidecomposable to a \fn\ l.s.c.\ for some such Lusin topology. Since
Borel-tame equidecomposability is transitive (Thm. \ref{Thm:TransitiveTame}),
we have:
}

Thus we have:

\begin{Thm} \label{Thm:ExsEqud}
Let a unimodular connected Lie group $G$ with Haar measure $\la$
act in a Borel manner on a standard Borel space $\OM$.
Let $E\st\OM$ be Borel with non-empty countable closed $\cO_\om E$ for every
$\om$.
Then any non-negative Borel $\fn$ $f$ on $\OM$
is Borel-tamely equidecomposable (Def. \ref{Def:Tame})
with a non-negative Borel $\fn$ on $E$ enhanced by $\frac{d\cnt}{d\la}$,
namely, with the enhanced \fn\ given by (\ref{eq:fPE}).

Consequently, if such an $E$ exists, 
then any invariant probability measure on $\OM$ can be recovered
from the corresponding $\frac{d\cnt}{d\la}$-enhanced infinitesimal measure,
even from its restriction to $E$.
\end{Thm}
\qed

\begin{Rmk}
\label{Rmk:2meas}
The ``flow under a \fn'' scene, as well as its generalization in this
\S, can be viewed as follows:

We have a standard space $\OM$ acted in a Borel
manner by a unimodular 2nd-countable locally compact group $G$ with Haar
measure $\la$. We have a Borel subset $E$ and we pick a right-invariant
measure $\La$ on $G$. On $E$ we are given a measure $\nu$, which
should hopefully be the restriction to $E$ of $\frac{d\La}{d\la}\mu$ for some
$G$-invariant (say, probability) measure $\mu$ on $\OM$.

Thus in the ``flow under a \fn'' construction we start with a standard
measure space $(E,\nu)$ with $\bZ$-action and a given positive measurable
\fn, we embed it into
a larger space $\OM$ with $\bR$-action s.t.\ the given $\bZ$-action is
given by $T_E$ and the given \fn\ is $\rho_E$ (as defined for the
$\bR$-case in \S\ref{SS:CtAsso}), and one finds a $\mu$,
invariant under the $\bR$-action, s.t.\ the initial measure $\nu$ in $E$ is
the restriction to $E$ of $\frac{d\cnt}{dt}\mu$.

For the general case
of an $E\st$ an $\OM$ acted by $G$ and a $\nu$ on $E$, we have at our
disposal the relation of two \fns\ on $E$ being tamely equidecomposable
when enhanced by $\frac{d\La}{d}$. A necessary condition for the existence
of a $\mu$ is that $\nu$ gives the same integral to such \fns, which we can
view as a substitute for the property of $G$-invariance for a $\nu$
``known'' only on $E$. If this condition is satisfied, and if moreover we have
the counterpart of Thm.\ \ref{Thm:ExsEqud} -- every ordinary Borel \fn\
on $\OM$ is tamely equidecomposable with some \fn\ supported on $E$,
enhanced by $\frac{d\La}{d\la}$, then we know the integral w.r.t.\
the sought-for 
$\mu$ of any Borel \fn\, and since this proposed integral is countably
additive on \fns, we have a unique $\mu$, where we have to check that
this $\mu$ is $\sigma$-finite (say, by showing that $1_\OM$ is a countable
sum of nonnegative \fns\ tamely equidecomposable with $\nu$-integrable
nonnegatve \fns\ on $E$). This view is of interest even
when $\La=\la$, so we are given an $E\st\OM$ and a $\mu$ on $E$ and seek
a $G$-invariant $\mu$ s.t.\ $\mu|_E=\nu$.

In the case of the ``flow under
a \fn'', $\nu$ is $T_E$-invariant, and it is easy to show that \fns\ on
$E$ are equidecomposable when the $\bR$-action on $\OM$ is considered and one
enhances the \fns\ by $\frac{d\cnt}{d\la}$ iff they are equidecomposable
w.r.t.\ the discrete $\bZ$-action given by $T_E$, so $\nu$ gives the same
integral to such \fns. Since, moreover, any Borel \fn\ $f$ on $\OM$ is
equidecomposable to a \fn\ on $E$ enhanced by $\frac{d\cnt}{dt}$, given,
say, by (\ref{eq:Ctsfs}) for $s\equiv1$:
$$
\LQ\int_0^{\rho(\om)}f(T^t\om)\,dt\RQ\,\frac{d\cnt}{dt},
$$
we are sure the required $\mu$
on $\OM$ exists and get a formula for $\int f\,d\mu$ in terms of $\nu$:
$$
\int_\OM f(\om)\,d\mu(\om)=
\int_E\LQ\int_0^{\rho(\om)}f(T^t\om)\,dt\RQ\,d\nu(\om).
$$

In the case dealt with in this \S\ -- $G$ is a Lie group and
$E$ is such that $\cO_\om E$ is always countable closed, one can speak of
a groupoid structure in $E$ instead of the nonexistent $T_E$, and \fns\
on $E$ are equidecomposable when enhanced by $\frac{d\cnt}{d\la}$ iff
they are equidecomposable w.r.t.\ to the (discrete) groupoid.
\end{Rmk}

\newpage
\section{Some Applications}
\label{S:APL}
\subsection{Comparison with a Result of G.~Helmberg;
          Persistence and Interruption of Patterns}
\label{SS:HE}

The continuous Kac Thm.\ \ref{Thm:CtsKac} bears some resemblance to a limit
theorem of G.~Helmberg \cite{Helmberg}
\NOT{
(\"Uber mittlere R\"uckkehrzeit unter einer masstreuen Str\"omung, {\it Z.\
Wahrsch.\ verw.\ Geb.\ 13, 165-179 (1969)}).
}
This theorem is formulated for
an $\bR^+$-action. Let us state it, partly using our notation:

Helmberg defines:
\begin{eqnarray}
\mbox{for }t>0 \quad E_t &:=& \{\om\in\OM : \exists r,s \mbox{ s.t. }
                              0\le r<s\le t \mbox{ \& } r\in \cO_\om E \mbox{ \& }
                              s\notin \cO_\om E\}\label{3.1}\\
r_E(\om) &:=& \inf\{s>0 : s\in \cO_\om E \mbox{ \& } \exists r \mbox{ s.t. }
                   0\le r<s \mbox{ \& } r\notin \cO_\om E\}\label{3.2}\\
\mbox{for }t>0 \quad E(t) &:=& E_t \cap \{\om : r_E(\om)<t\}\label{3.3}
\end{eqnarray}
He proves, using discrete approximation (i.e.\ the $\bZ^+$-actions $T^{nt}$
for $t>0$) the following

\begin{Thm} (G.~Helmberg) Suppose
$\LP\OM,\cB,\mu,(T^t)_{t\in\bR^+}\RP$
is a measure-preserving flow, and $E\in\cB$ is s.t.\ $E_t$ , $t>0$ and
$r_E$ are measurable, and s.t.\ $\spn E$ is conull. Suppose
   \begin{eqnarray}
     \lim_{s\to 0} \mu(E_s)=0 & \label{3.4}\\
     \lim_{s\to 0} {1\over s} \mu(E(s))=0 & \label{3.5}
   \end{eqnarray} 
Then
\begin{equation}
   \lim_{s\to 0} {1\over s} \int_{E_s} r_E(\om)\,d\mu(\om) = 1-\mu(E)
   \label{3.6}
\end{equation}
(we have added the inessential requirement that $\spn E$ 
is conull -- $\spn E$ being always measurable if the $E_t$'s are).%
\NOT{see the above-mentioned article of Helmberg}%
\footnote{In the second part of \cite{Helmberg}
Helmberg defines, for $E$ closed in a compact metric $\OM$ and continuous
action, the notions of $\Ex E$ and $\In E$ (different from ours), these
being subsets of $\OM$. He formulates requirements on them in this
topological setting that insure (\ref{3.4}) and (\ref{3.5})
hence (\ref{3.6}). $\Ex E$ is the set of $\om$ in $E$ which will visit $E^c$
in every $[0,t]$-future. Similarly for $\In E$.}
\end{Thm}

This theorem, for Borel $\bR$-action and $E$ Borel, can be proved using the
ideas of Section \ref{S:CT}. This will be presented here.

\NOT{
The $\bR^+$ case can be handled too, applying an
extension of our methods to semigroups -- see......
}

Let us remark first, that Helmberg's definitions may be described in the
context of ``persistence of a certain pattern''. Suppose we wish to speak
about ``moments of exit from $E$''. To do this we may consider the following
``pattern'': an interval of Time being the union of an interval of stay
outside $E$ and a subsequent interval of stay in $E$. More precisely: for any
$A\subset\bR$, consider the union of all open intervals $I\subset\bR$ such
that in the decomposition $I=(I\cap A)\cup(I\cap A^c)$ both intersections are
intervals (possibly empty) and $I\cap A^c$ precedes $I\cap A$, this being our
``pattern''. Let $\Gamma(A)$ be the {\em complement} of the union of all such
intervals, thus $\Gamma(A)$ is composed of Time moments when the pattern
is interrupted. For $A=\cO_\om E$
These time moments may be viewed naturally as ``moments of exit'' from $E$
(for a fixed $\om$). Similarly, the pattern leading to ``moments of entry''
will be obtained by requiring $I\cap A$ to precede $I\cap A^c$.

Assume we are in our usual setting of Borel $\bR$-action on a standard
space and $E$ is Borel.

Note that it is not hard to convince oneself, using
Prop.\ \ref{Prop:EGDelta} (i),
that for $\Gamma$ referring to ``moments of exit'' or to ``moments of entry'',
the closed set $\Gamma(\cO_\om E)$ is a measurable \fn\ of $\om$, being
a Borel \fn\ of $\BAR{\cO_\om E}$.

The graph used in the proof of the \cts\ Kac Thm.\ \ref{Thm:CtsKac} can be
constructed from the the closed set $C=\Gamma(\cO_\om E)$ instead of
$\BAR{\cO_\om E}$, the latter being the case where the pattern is:
``staying out of $A$''. ``Return time'' and ``arrival time'' can be defined
for general $\Gamma$, for example, the arrival time for the pattern of
``exit from $E$'' will have the meaning of ``the waiting time to exit $E$''.
An analog of Thm.\ \ref{Thm:CtsKac} can be proved.

Consequently, Helmberg's $r_E(\om)$ is a.e.\ ``the waiting time to enter $E$''
and the measurability of $\Gamma(\cO_\om E)$ as a \fn\ of $\om$ (for the
``entry'' pattern) implies Helmberg's condition that $r_E(\om)$ be measurable.

Also, for a.a.\ $\om$, $\om$ belongs to $E_t$ iff in the $[0,t]$-future
there is a ``moment of exit''; and it belongs to $E(t)$ iff this future
contains both a moment of exit and a moment of entry. 

Note that (\ref{3.6}) involves both the ``entry'' pattern ($r_E$) and the
``exit'' pattern ($E_t$).
 
\NOT{
(Indeed, an open $u$ is disjoint from $\Gamma(\cO_\om E)$ iff, say,
$(u\cap(\cO_\om E))\times(u\cap(\cO_\om E^c))$
is disjoint from $\{(t,s):t>s\}$. This is the complement of a Suslin set
$\st\OM$ which is a projection of a Borel set in $\OM\times\bR\times\bR$.) 
}

\par\medskip

The following lemma is in the spirit of Satz 3 in \cite{Helmberg}:

\begin{Main} {\EM Lemma}
For Borel action and $E$ borel, and assuming Helmberg's
condition (\ref{3.5}), let $\rho$ be a (nonnegative)
measurable ($\om$-dependent) invariant $0$-chain with Radon (i.e.\ Haar)
expectation. Then the integral $I(t)$ over $\om\in E(t)$ of the
integral of $1_{]0,t[}$ over $\rho$ is $o(t)$ as $t\to 0$.	
\end{Main}

\begin{Prf}
(see the proof of Satz 3 in \cite{Helmberg}):
Since $E(t)$ depends increasingly on $t$, $I(t)$ is an increasing
\fn\ of $t$. Hence it suffices to prove $I(t/n)=o(1/n)$ as $n\to\infty$
for fixed $t$, i.e.\ $nI(t/n)\to 0$ as $n\to\infty$.
To compute $I(t/n)$ we integrate $1_{]0,t/n[}$ over $\rho$, and then
integrate over the $\om$ with entry and exit (w.r.t.\ $E$) in the
$]0,t/n[$-future.
Substituting $T^{kt/n}\om$ for $\om$, we get the same by integrating
$1_{]kt/n,(k+1)t/n[}$ over $\rho$ and then integrating over the set of $\om$
with entry and exit in the $]kt/n,(k+1)t/n[$-future, the latter set
having measure $\mu(E(t/n))$. Summing this for $k=0,\ldots,n-1$,
we have $nI(t/n)\le$ the integral of $\int 1_{]0,t[}\,d\rho$ over a set
$\OM_1$ of $\om$'s with measure $\le n\mu(E(t/n))$. By (\ref{3.5}),
$\mu(\OM_1)\to 0$. Since $\rho$ is measurable with Haar expectation,
$\int 1_{]0,t[}\,d\rho$ belongs to $L^1(\OM)$. Therefore its
integral over $\OM_1$ tends to $0$ as $n\to\infty$, and we are done.
\end{Prf}

Now consider the following $\om$ dependent graph:
Let $\Gamma$ refer to the ``exit'' pattern. an arrow $(t_0,t_1)$ belongs
to the graph iff $t_1<t_0$, $t_0\notin\cO_\om E$, $t_1\in\Gamma(\cO_\om E)$
and $]t_1,t_0]\cap\cO_\om E=\es$. Take the $1$-simplex of invariant measures
$(dt,\cnt_\bR)$, and apply CHG (Check that its requirements are satisfied).
One finds:

Let
\BER{l}
 E'=\{\om\in\OM:0\in\Gamma(\cO_\om E)\mbox{ \& }
       \exists\eps>0\mbox{ s.t. }]0,\eps]\cap\cO_\om E=\es\}\\
 E''=\{\om\in\OM:\Big]\sup\LP]-\I,0]\cap\Gamma(\cO_\om E)\RP,0\Big]
       \mbox{ is nonempty and disjoint from }\cO_\om E\}
\EER
(Note that $\cO_\om E'\st\Gamma(\cO_\om E)$,
$E''\st E^c$, and $\cO_\om(E^c\sm E'')\st\Gamma(\cO_\om E)$.)

The source: a.e.\ the ordinary \fn\ $1_{E''}$.

The target: $r_E1_{E'}\frac{d\cnt}{dt}$.

Thus, by CHG, $r_E1_{E'}$, enhanced by $\frac{d\cnt}{dt}$, has expectation
$\mu(E'')$.

Consider the $0$-chain enhancement
$\rho=1_{E^c\sm E''}+r_E1_{E'}\frac{d\cnt}{d}$ i.e.\
$1_{\cO_\om (E^c\sm E'')}\,dt+\cO_\om\LP r_E1_{E'}\RP\,d\cnt_\bR$.
By the above, $\rho$ is supported on moments of exit for $\om$.

To obtain (\ref{3.6}), note that, by the above, $\rho$ has expectation
$1-\mu(E)$. Thus $1/t$ times the integral of $1_{]0,t[}$ over $\rho$ has
expectation $1-\mu(E)$. Since $\rho$ is supported on moments of exit, this
integral is $0$ for $\om\notin E_t$. For $\om\in E_t\setminus E(t)$
(i.e.\ exit but no entry in the $]0,t[$-future, hence only one moment of
exit then) the integral differs from $r_E(\om)$ by no more than $t$,
and one notes (\ref{3.4}). The $E(t)$ part is disposed of by the lemma
and $r_E(\om)$ being $\le t$ there, and we are done.\qed

\subsection{Aaronson and Weiss's Kac Functions}
\label{SS:AW}
In \cite{AaronsonWeiss}
J.~Aaronson and B.~Weiss make crucial use of a {\EM Kac \fn}
of a subset $E\subset\OM,\;\mu(E)>0$ of a measure space $(\OM,\mu)$
on which the group $G=Z^d$ acts measure-preservingly. In $Z^d$ take
the $l^\I$-norm.

They define a {\em Kac \fn} as a measurable
\fn\ $\phi:E\to\bZ^+$ satisfying
\begin{eqnarray} \label{eq:AWKacFn11}
\bigcup_{n=1}^{\I}\bigcup_{\|x\|\le n} T^x\left\{\om:\phi(\om)=n\right\}=
\OM\mbox{ mod }\mu\\ 
\label{eq:AWKacFn12}
\int_E\phi^d\,d\mu<\I
\end{eqnarray}

That is, to every $\om\in E$ one corresponds a 
``radius'' $\phi(\om)$ in
$G=\bZ^d$ s.t.\ the union, for $x\in E_{\om}$, of the cubes with radius
$\phi(T^x\om)$ and center $x$ is a.e.\ the whole of $G$, while the ``volume''
of the cube has finite expectation.

In their above article they seem to
promise to prove in a future article that in the ergodic case such a
Kac \fn\ alway exists, while in the present article they prove its
existence (in the ergodic case) in a random sense,
\NOT{ (we prefer to say:
non-deterministic sense, in order not to confuse with ``random sets''),}%
i.e.\ when one passes to an appropriate extension of the dynamical system
(in their case, a product) and lets $\phi$ depend on the points of the
extension.
\NOT{Note that the
existence of a Kac \fn\ for every ergodic $\bZ^d$-action is equivalent
to its existence in the particular case of an (invariantly) ergodic random
subset of $\bZ^d$.}

In the sequel, a way
to get such a ``random'' Kac \fn\
\NOT{(for the ergodic case)}%
(for the case $\spn E$ conull)
is presented, which is essentially Aaronson and Weiss's method.
Instead of the pointwise ergodic theorem for $Z^d$-action, which they apply,
the HG theorem is used. 

Consider the compact group $\bZ_2$ of the {\em dyadic integers}.
Every element of $\bZ_2$ has an infinite dyadic expansion with
ascending powers of $2$.

Take the normalized (i.e.\ probabilistic) Haar measure in $\bZ_2^d$,
denoted by $\la$.

One has the {\em $\bZ^d$-dyadic odometer}, which is the set $\bZ_2^d$ acted
upon measure-preservingly by $\bZ^d$ via addition:

$$T^x\al=x+\al\quad x\in\bZ^d,\al\in\bZ_2^d$$

An element $\al\in\bZ_2^d$ can be identified with a
hierarchy of dyadic partitions of
$\bZ^d$, where in the $n$-th step $\bZ^d$ is partitioned into cubes with
side $2^n$ (to be called $n$-cubes of the hierarchy), given by inequalities
of the form
$$k_i\le x_i<k_i+2^n,\quad i=1,\ldots,d$$
In the hierarchy corresponding to an $\al\in\bZ_2^d$, each $n$-cube will
consist of those $x\in\bZ^d$ with common $n+1,n+2,\ldots$'th digits in the
dyadic expansion of the $d$ coordinated of $x+\al$. When the odometer is
identified with the set of hierarchies of partitions of $\bZ^d$, $\bZ^d$
acts on these hierarchies by shift.

Now consider the product dynamical system $\OM\times\bZ_2^d$.
Our aim is to prove the existence of a Kac \fn\ $\phi(\om,\al)$
on $\OM\times\bZ_2^d$, i.e.\ $\phi$ has to satisfy

\begin{eqnarray} \label{eq:AWKacFn21}
\bigcup_{n=1}^{\I}\bigcup_{\|x\|\le n}
T^x\left\{(\om,\al):\phi(\om,\al)=n\right\}=
\OM\mbox{ mod }\mu\times\la\\
\label{eq:AWKacFn22}
\int_E\phi^d\,d(\mu\times\la)<\I
\end{eqnarray}

\NOT{
{\EM Remark before proof} The product $2^G\times\bZ_2^d$ need not be
ergodic. If it is not, the sub-$\sigma$-algebra of invariant sets in it
defines a factor (all spaces being standard) on which $G$ acts trivially,
and the measure disintegrates into measures supported on fibers, these measures
giving {\it ergodic} dynamical systems.
Also, a.a.\ of them will satisfy the analogues of (5.6') and (5.7'). So we
obtain satisfactory {\it ergodic} extensions.
}

We correspond to every $(\om,\al)$ a graph $F(\om,\al)\subset G^2$,
in a measurable and invariant way, so that a.e.\ from every $x\in G$
emanates a unique edge, which terminates in $\cO_\om E$,
while the set $S=S(\om,\al)$ of $y$'s from which
an edge goes to $0$ is a.e.\ (in $E$) ``thick'' in the sense that
\begin{equation} \label{eq:AWThick}
\mbox{card}(S)\ge M\cdot\LP\diam(S\cup\{0\})\RP^d
\end{equation}
$M$ being a constant (which may depend on $d$).

By HG, the expectation of $\mbox{card}(S(\om,\al))$ is $1$,
hence if $\phi(\om,\al)$ is the radius of the smallest cube centered at
$0$ and containing $S(\om,\al)$, (\ref{eq:AWKacFn22}) is satisfied.
(\ref{eq:AWKacFn21}) is satisfied too: 
Indeed, we need prove that a.a.\ $(\om',\al')\in \OM\times\bZ_2^d$ is a
member of the left-hand side of (\ref{eq:AWKacFn21}). 
There is a edge emanating from $0$ $(0,-x)\in F(\om',\al')$. Invariance
(\ref{eq:Invr}) implies $(x,0)\in F(T^{-x}\om',\al'-x)$, hence
$x\in S(T^{-x}\om',\al'-x)$, thus $\|x\|\le n$ where
$$n=\phi(T^{-x}\om',\al'-x)=\max_{z\in S(T^{-x}\om',\al'-x)}\|z\|$$
and $(\om',\al')=T^x(T^{-x}\om',\al'-x)\in T^x\{\phi=n\}$,
verifying (\ref{eq:AWKacFn21}).

\NOT{
That is, to every $\om\in E$ one corresponds a cube of radius $\phi(\om)$
in $G=\bZ^d$ s.t.\ the union, for $x\in E_{\om}$, of the cubes with
radius $\phi(T^x\om)$ and center $x$ is a.e.\ the whole of $G$, while the
``volume'' of the cube has finite expectation.
}

We proceed to construct a graph $F(\om,\al)$ with the ``thickness''
property (\ref{eq:AWThick}).
(The graph should depend on $(\om,\al)$ invariantly and
measurably.)
We have to specify the target of the edge emanating from
some $x\in \bZ_2^d$. As mentioned above, $\al$ may be viewed as a hierarchy
of partitions of $\bZ_2^d$ into cubes, the cubes of the $n$-th partition
(which have side $2^n$) will be called $n$-cubes.

The set $\cO_\om E = \{x \in G : x\om\in E\}$ is a.e.\ non-void, $\spn E$
being conull.

We may restrict ourselves to $\al$ belonging to the conull set 
$(\bZ_2\sm\bZ)^d$. For such $\al$, for any $x\in\bZ^d$ the union of cubes
containing $x$ for all the partitions in the hierarchy is $G=\bZ^d$.

The target of the edges emanating from the points $x\in G$ will be
found in steps: In the $n$-th step, we review the cubes of the $n$-th
partition. For each such cube that intersects $\cO_\om E$, take, say, the
first $z$ in
the intersection w.r.t.\ the lexicographic ordering, and let all points
in the cube at which the target had not been defined in previous steps
be given the target $z$.

For $\om\in E$, this makes $S(\om,\al)$, i.e.\ the set of $x\in G$ with edge
going to $0$, a subset of some $n$-cube which contains an $n-1$ -subcube
($n$ being the last step in which $0$ was designated as a target --
for $\al\in(\bZ_2\sm\bZ)^d$ there is always such a last step if $0$ is not
the first element of $\cO_\om E$ lexicographically, which happens a.a.\ by
the argument of Poincar\'e's recurrence (see \S\ref{SS:EqiChn})).

Thus (\ref{eq:AWThick}) is satisfied with $M$ depending only on $d$,
and we are done.

{
\subsection{Applications to The Renewal Theorem}
\label{SS:Rnwl}
\newtheorem{Asm}{Assumption}
\newcommand{\laa}{{\tilde{\la}}}

We shall sketch how formulas obtained from CHG can be applied to prove
some classical renewal limit theorems (see \cite{Feller}, \cite{Durrett},
\cite{Lindvall}, \cite{NeveuP}, \cite{DaleyVereJones}).%
\footnote{I am indebted to Prof.\ Jon Aaronson for suggesting to me
the possibility of applying methods of this work to renewal problems.}
The proofs will be presented for the classical case of stationary independent
renewal times with finite expectation, both for nonperiodic
continuous Time and for discrete Time (where these two cases will be treated
completely analogously). It seems that one can extend the proofs to apply to
the stationary case with more general ``Time'' (i.e.\ acting groups other
than $\bZ$ and $\bR$) and assume much less than independence, but in these
settings one has to make some technical assumptions to make the proof work.
Anyhow, we shall make restrictive assumptions only when we need them.

Let $G$ be the acting group, assumed {\em abelian},
which eventually will be $\bZ$ or $\bR$.
Denote by $\cnt$ the counting measure on $G$ and by $dx$ a Haar measure
in $G$: Lebesgue in the case $G=\bR$ and $\cnt_\bZ$ in the case $G=\bZ$.
The mass of a measurable set $S\st G$ w.r.t.\ $dx$ will be denoted by
$|S|$.

We are given a probability space $(E,\nu)$ where $E=\LP\bR^+\RP^\bZ$ with
the usual Borel structure, and $\nu$ is stationary, i.e.\ shift-invariant.
Denote by $X_n$ the stochastic variable equal to the $n$-th coordinate,
and define $S_n=\sum_{1\le j\le n}X_j$, 
$S_{-n}=\sum_{1\le j\le n}X_{-j}$ for $n\ge0$.
The $X_n$ are the {\EM renewal times}
and {\EM $S_n$ is the time of the $n$'th renewal}.
Assume also that each $X_n$ is a.s.\ positive, and that a.s.\
$\lim_{n\to\I}S_n=+\I$, $\lim_{n\to-\I}S_n=-\I$.
\NOT{the sum of positively
indexed (resp.\ negatively indexed) coordinates is a.s.\ $+\I$.}
Identifying $\om\in E$ with $\LB S_n:n\in\bZ\RB$,
$E$ can be identified (after taking away a null subset, if necessary)
with the set of all discrete subsets of $G=\bR\mbox{ or }\bZ$
unbounded below and above and containing $0$ and as such it is a subset of
$\OM=$ the set of all discrete subsets of $G$ unbounded below and above.
On the latter $G$ acts by translations: $x\om:=\om-x$ and there is
a unique $\bR$-invariant\NOT{$\sigma$-finite} measure $\mu$ on $\OM$
(not necessarily probability or finite, but it will be shown below that
$\mu$ is $\sigma$-finite), s.t.\
$\nu=\frac{d\cnt}{dx}\mu|_E$. In other words, $\nu$ is the Palm
measure for $\mu$ (see \S\ref{SS:Palm}). $\mu$ can be obtained as in
Rmk.\ \ref{Rmk:2meas} or as in \cite{Neveu} (see the proof of f.\ in
Thm.\ \ref{Thm:JointDist}).

Note that since $E=\{\om\in\OM:0\in\om\}$, we have that $\forall\,\om\in\OM$
$\om$, as a subset of $\bR$, is identical with $\cO_\om E$.

The renewal limit theorems that we have in mind state that,
under some assumptions,
the $\nu$-integral of the sum over $\cO_\om E(=\om)$ of a translation
$x\mapsto h(a^{-1}x)$ of a fixed Borel \fn\ $h:E\to\BAR{\bR^+}$ tends
as $a\to+\I$ to $(\mu(\OM))^{-1}\int h(x)\,dx$. (Where if $\mu(\OM)=\I$
it is understood that $(\mu(\OM))^{-1}=0$.)
The fact that for $G=\bR$ this does not hold for every Borel $h$ is 
clear if one considers the case when the ranges of all the $X_n$ are
countable.

\subsubsection{The Case of Mixing $G$-action}

Let us play with CHG to get some formulas (recall that we assumed $G$ 
abelian)\NOT{(keep in mind that $\cO_\om E=\om$)}:

\begin{itemize}
\item
Let $h(x)$ be a Borel $\BAR{\bR^+}$-valued \fn\ on $G$.

Take the $\om$-dependent graph, consisting of the $1$-simplices
$(x,y)\in G^2$ with $y\om\in E$, weighted by $h(yx^{-1})$.
The $1$-simplex of measures on $G$ will be taken as $(dx,\cnt)$.

The source is the \fn\
$$H(\om):=\sum_{y\in\cO_\om E}h(y).$$

The target is the enhanced \fn\
$$\LP\int h(x)\,dx\RP1_E(\om)\cdot\frac{d\cnt}{dx}.$$

Thus by CHG these two have the same $\mu$-integral, i.e.\ the $\mu$-integral
of $H$ on $\OM$ equals $\int h\,dx$ times $\nu(E)=1$.

\item
Let $h_1(x)$ and $h_2(x)$ be Borel $\BAR{\bR^+}$-valued \fns\ on $G$.

Take the $\om$-dependent $2$-hypergraph, consisting of the $2$-simplices
$(x,y,z)\in G^3$ where $y\om,z\om\in E$, weighted by
$h_1(yx^{-1})h_2(zx^{-1})$. The $2$-simplex of measures on $G$ will be
$(dx,\cnt,\cnt)$.

The $0$-vertex is the \fn\
$$\LP\sum_{y\in\cO_\om E}h_1(y)\RP\cdot\LP\sum_{z\in\cO_\om E}h_2(z)\RP.$$

The $1$-vertex is the enhanced \fn\
$$\LP\sum_{z\in\cO_\om E}(h_1*h_2)(z)\RP1_E(\om)\cdot\frac{d\cnt}{dx},$$
where $h_1*h_2$ is the convolution:
$(h_1*h_2)(z):=\int_G h_1(x)h_2(zx)\,dx.$

Thus by CHG these two have the same $\mu$-integral.

For $\om\in\OM$, let $H_i(\om):=\sum_{y\in\cO_\om E}h_i(y)$, $i=1,2$.

Assume $\mu(\OM)<\I$.
Take $h_1$ fixed but replace $h_2$ by the translated
$z\mapsto h_2(a^{-1}z)$, $a\in G$.
One finds that the $\nu$-integral on $\om\in E$ of the sum over $\cO_\om E$
($=\om$) of the translated $z\mapsto (h_1*h_2)(a^{-1}z)$
is equal to the inner product w.r.t.\ $\mu$ on $\OM$ of $H_1$ and the
translated $\om\mapsto H_2(a\om)$. If the $G$-action
on $\OM$ is {\em mixing} this inner product tends, as $a\to\I$, to
$(\mu(\OM))^{-1}\int H_1\,d\mu\cdot\int H_2\,d\mu$
(assuming, say, that $h_1$ and $h_2$ are
bounded with compact support). But by the previous item, the latter equals
$(\mu(\OM))^{-1}\int h_1,d\nu\cdot\int h_2,d\nu=
(\mu(\OM))^{-1}\int h_1*h_2\, d\nu$. Thus we get the
well-known fact (see \cite{Delasnerie}):
\par\medskip
{\em If $\mu(\OM)<\I$ then the $G$-action on $\OM$ being mixing implies
holding of the renewal limit theorem for \fns\ $h$ on $G$ which are
convolutions of two bounded Borel \fns\ with compact support.}
\end{itemize}

\subsubsection{The Nonperiodic Case (for $\bE X_1<\I$)}

We use the same notation as in \S\ref{SS:Nearest} and \S\ref{SS:CntCls},
but with ``nearest'' understood not as in \S\ref{SS:Nearest} but as an
arbitrary but fixed Borel measurable translation-invariant way to correspond
to any closed subset $F\st G$ and point $x\in G$ a point in $F$.
In fact, our way to do so (for $G=$ $\bR$ or $\bZ$) will always be to take
the point $\max(F\cap\,]-\I,x])$. Thus, according to the notation there,
for $\om\in\OM$ $\pi(\om)=\pi_E(\om)$ is the last point in
$\cO_\om E$ ($=\om$) not bigger than $0$
and for $\om\in E$ $P(\om)=P_E(\om)$ is the set of $x\in G$
with $0$ as the last point in $\cO_\om E$ not bigger than $x$, that is,
$P(\om)=\left[0,\min\LP\,]0,\I]\cap\cO_\om E\RP\right[\,=[0,X_1(\om)[$.
But for the time being, assume ``nearest'' to be interpreted in an arbitrary
way as above.

Since we assume $G$ is abelian, we use additive notation for $G$, compelling
us to write $T^x\om$ instead of $x\om$.

Let us ``play'' with CHG:

Take the $\om$-dependent graph (cf.\ (\ref{eq:Nearest})) 
\BER{l}
\{(x,y)\in G^2: y \mbox{ is the ``nearest'' point to }x
\mbox{ in }\cO_\om E(=\om)\} =\\=
\{(x,y)\in G^2: y-x=\pi_E(T^x\om)\}=
\{(x,y)\in G^2: T^y\om\in E,\,x-y\in P_E(T^y\om)\},
\EER
and the $1$-simplex of measures on $G$\quad$(dx,\cnt)$.

The source is the \fn\ $1$. The target is the enhanced \fn\ 
$|P_E(\om)|1_E(\om)\frac{d\cnt}{dx}$. Hence by CHG:
\begin{equation} \label{eq:RnwlK}
\mu(\OM)=\int_E|P_E(\om)|\frac{d\cnt}{dx}\,d\mu(\om)=
\int_E|P_E(\om)|\,d\nu(\om)
\end{equation}
Note that both sides may be $+\I$, yet if we take in the above graph only
the edges $(x,y)$ with $y-x\in K_n$ where $K_n$ are compact with union $G$,
we deduce, using CHG, that $\mu$ on $\OM$ is {\em $\sigma$-finite}:
$\OM$ is the union of the sets with finite $\mu$-measure
$\{\om:\pi(\om)\in K_n\}$.

\NOT{
Choosing for ``nearest'' our usual meaning $P(\om)=[0,X_1(\om)[$
the last expression is $<\I$ we have the
(Kac) formula
$$\mu(\OM)=\int_EX_1(\om)\,d\nu(\om)$$
We shall assume from now on that the last expression is $<\I$
}

We shall assume from now on that $\mu(\OM)=\int_E|P(\om)|\,d\nu(\om)$
is {\em finite}. For our meaning of ``nearest''
$P(\om)=[0,x_1(\om)[$ this and (\ref{eq:RnwlK}) mean that
$\int_EX_1(\om)\,d\nu(\om)=\mu(\OM)<\I$

Now for $a\in G$, Take the $\om$-dependent graph
\BER{l}
\{(x,y)\in G^2: x\om\in E, y \mbox{ is the ``nearest'' point to }x-a
\mbox{ in }\cO_\om E(=\om)\} =\\=
\LB(x,y)\in G^2: T^x\om\in E, y-x=-a+\pi_E\LP T^{x-a}\om\RP\RB =\\=
\LB(x,y)\in G^2: T^x\om,T^y\om\in E,\,x-y\in a+P_E(T^y\om)\RB,
\EER
and the $1$-simplex of measures $(\cnt,\cnt)$.

The source is the enhanced \fn\ $1_E\frac{d\cnt}{dx}$.
The target is the enhanced \fn\ 
$$\#\om\cap\LP a+P_E(\om)\RP\frac{d\cnt}{dx}.$$

Thus CHG tells us that
\begin{equation}\label{eq:RnwlT}
\int_E\#\om\cap\LP a+P_E(\om)\RP\,d\nu=1
\end{equation}
i.e.\ that substituting in the renewal limit theorem for $h(x)$ the 
$\om$-dependent $1_{p(\om)}(x)$ gives an identity, rather than a limit:
summing the translated \fn\ on $\om=\cO_\om E$ and taking the $\nu$-expectation
gives the correct value (at least when $\mu(\OM)<\I$):
$$1=(\mu(\OM))^{-1}\int_E\LP\int1_{P(\om)}(x)\,dx\RP\,d\nu(\om)$$
(see (\ref{eq:RnwlK})).

From this we shall deduce the limit theorem in a ``Tauberian'' manner.

\NOT{
Let $\gamma(\om)$ be\NOT{the $\nu$-expectation of} the $\om$-dependent measure
on $G$ equal to the counting measure on $\om=\cO_\om E\st G$, i.e.\ 
giving, for each subset of Time $G$ the \# of renewals in this subset.
}

Denote by $\cC_{00}=\cC_{00}(G)$ the space of \cts\ \fns\ on $G$ with compact
support.

Our goal is to show that:
\begin{equation}\label{RnwlE}
\lim_{a\to+\I}\int_E\LP\delta_{-a}*\sum_{x\in\om}\delta_x\RP
\,d\nu(\om)=(\mu(\OM))^{-1}\,dx,
\end{equation}
the limit taken in the weak topology w.r.t.\ $\cC_{00}$.
This will imply holding of the renewal limit assertion
for $h(x)$ bounded Riemann-integrable with compact support.

Note, that the measure on $G$\quad$\sum_{x\in\om}\delta_x$ gives, 
for each subset of Time $G$, the \# of renewals in this subset.

We shall make assumptions needed to carry out our proof and show that they
hold for the case considered by us.

\begin{Asm} \label{Asm:B}
The Borel \fns\
$\om\mapsto\int_G h(x)\,d\LP\delta_{-a}*\sum_{x\in\om}\delta_x\RP(x)=
\sum_{x\in\om}h(x-a)$,
$a\in G$
form a weakly compact family in $L^1(E,\nu)$ for each fixed $h\in\cC_{00}$.
\end{Asm}

This assumption holds when $G=\bR$ or $\bZ$ and the $X_n$ are independent.
Indeed, one may consider $h=1_{[0,c]}$ instead of $h\in\cC_{00}$.
Now, if $S_m$ is the first one not less than $a$, then the \# of renewals
in $[a,a+c]$ is $\le$ than the \# of renewals in $\LQ S_m,S_m+c\RQ$,
which by independence has the same distribution \fn\ as the \# of renewals
in $[0,c]$.
The latter does not depend on $a$, and is an $L^1$-\fn\ of $\om\in E$, since
the probability that the \# of renewals in $[0,c]$ $\ge n$ is
\BER{l}
\Pr(S_n\le c)\ \le \Pr\LP\sum_{1\le j\le n}\min(X_j,1)\le c\RP=\\=
\Pr\LP\sum_{1\le j\le n}\min(X_j,1)-n\bE\LP\min(X_1,1)\RP\le
c-n\bE\LP\min(X_1,1)\RP\RP=\\=
O(n^{-2}).
\EER

By Assumption \ref{Asm:B} the set of measures
$\int_E\LP\delta_{-a}*\sum_{x\in\om}\delta_x\RP\,d\nu(\om)$, $a\in G$
is bounded on every compact $\st G$, hence is contained
in a compact metrizable set in the weak topology w.r.t.\ $\cC_{00}$.
Thus assume that $a_j\to+\I$ such that
$\int_E\LP\delta_{-a_j}*\sum_{x\in\om}\delta_x\RP\,d\nu(\om)$
converges in the weak topology w.r.t.\ $\cC_{00}$ to some measure $\laa$
on $G$ (necessarily a Radon measure that is uniformly bounded on the
translations of any fixed compact set).
It suffices to prove that for each such $(a_j)$\quad
$\laa=(\mu(\OM))^{-1}\,dx$.

\begin{Asm} \label{Asm:01}
The Hewitt-Savage 0-1 law holds, namely, every permutable event in $E$, i.e.\
every event not changed by any permutation of the $X_n$'s that moves only
a finite number of indices, has probability $0$ or $1$.
\end{Asm}

This holds when the $X_n$'s are independent -- see \cite{Durrett} \S3.1.

Now, Assumption \ref{Asm:01} implies that for every $h\in\cC_{00}$,
every limit of a subsequence of
$$
\om\mapsto\int_G h(y)\,d\LP\delta_{-a_j}*\sum_{x\in\om}\delta_x\RP(y)
=\sum_{x\in\om}h(x-a_j)
$$
in the weak topology of $L^1(E,\nu)$, being measurable w.r.t.\ the
$\sigma$-algebra of permutable events, must be a.e.\ constant, necessarily
equal to the limit of the $\nu$-integrals $\LP\int_G h(x)\,d\laa(x)\RP\cdot1$.
By Assumption \ref{Asm:B}, the latter is the weak $L^1$-limit of
$\om\mapsto\sum_{x\in\om}h(x-a_j)$.
This means that for every \fn\ $F$ on $E\times G$ which is of the form
$(\om,x)\mapsto f(\om)h(x)$, $f\in L^\I(E)$, $h\in\cC_{00}$,
one has

\begin{equation} \label{eq:RnwlF}
\lim_j\int_E\sum_{x\in\om}F(\om,x-a_j)\,d\nu(\om)=
\int_G\int_E F(\om,x)\,d\nu(\om)\,d\laa(x).
\end{equation} 

We would like to have (\ref{eq:RnwlF}) for $F(\om,x)=1_{a+P(\om)}(x)$, so that
we may compare (\ref{eq:RnwlT}) and (\ref{eq:RnwlF}).

In the set of open (resp.\ closed) subsets of $\OM$ we take the Effros Borel
structure (see \S\ref{SS:MEAS}), given by identifying each open set $S\st G$
(resp.\ each closed $S\st G$) with the set $\{u\in\cU:u\st S\}\st2^\cU$
(resp.\ $\{u\in\cU:u\cap S=\es\}\st2^\cU$),
where $\cU$ is some countable base to the topology.

\begin{Asm} \label{Asm:R}
For $\nu$-a.a.\ $\om\in E$ $|\partial(P(\om))|=0$ and the mappings sending
$\om\in E$ to the interior $P(\om)^\circ$, resp.\ to the closure
$\BAR{P(\om)}$, are measurable.
\end{Asm}

This assumption clearly holds for our case $P(\om)=[0,X_1(\om)[$.

From Assumption \ref{Asm:R} one deduces that the \fn\ on $E\times G$\quad
$(\om,x)\mapsto1_{P(\om)^\circ}(x)$ is a supremum of countably many
finite linear combinations of \fns\ of the form
$(\om,x)\mapsto 1_S(\om)h(x)$, $S\st E$ measurable and $h\in\cC_{00}$,
hence is a limit of an increasing sequence $(F_k)$ of such combinations.
(Note that the collection of such linear combinations is stable w.r.t.\ the
lattice operations $\min$ and $\max$.)
By (\ref{eq:RnwlF}) we have, for any $a\in G$:
$$
\lim_j\int_E\sum_{x\in\om}F_k(\om,x-a-a_j)\,d\nu(\om)=
\int_G\int_E F_k(\om,x-a)\,d\nu(\om)\,d\laa(x).
$$
Going to the limit in $k$ one obtains:
$$
\lim_j\int_E\#\LP\om\cap(a+a_j+(P(\om))^\circ\RP\,d\nu(\om)\ge
\int_G\int_E1_{P(\om)^\circ}(x-a)\,d\nu(\om)\,d\laa(x)
$$
and taking into account (\ref{eq:RnwlT}), one has
$$
1\ge\int_G\nu\{\om:x-a\in P(\om)^\circ\}\,d\laa(x).
$$
The last expression, as a \fn\ of $a$, is the convolution of $\laa$ with
the \fn\ $x\mapsto\nu\{\om:x\in P(\om)^\circ\}$. By Assumption \ref{Asm:R},
this \fn\ is equal $dx$-a.e.\ to $x\mapsto\nu\{\om:x\in P(\om)\}$,
therefore its convolution with $\laa$ is equal $dx$-a.e.\ to
$\laa*(x\mapsto\nu\{\om:x\in P(\om)\})$. Thus we conclude:
\begin{equation} \label{eq:RnwlI}
1\ge\laa*(x\mapsto\nu\{\om:x\in P(\om)\})\qquad dx-\mbox{a.e.} 
\end{equation}

If we knew that $\forall\om P(\om)\st$ a fixed compact $K$, we could 
similarly prove the opposite inequality by considering $U\sm P(\om)$
instead of $P(\om)$ where $U$ is fixed open, relatively compact and contains
$K$. But that need not be the case. Clearly, we shall still have the opposite
inequality if the following assumption holds:

\begin{Asm} \label{Asm:C}
For every $\eps>0$ $\exists$ is a relatively compact set $K$ s.t.\ for all
$a\in G$ close enough to $+\I$
\begin{equation} \label{eq:RnwlC}
\int_E\#\om\cap\LP a+\LP P(\om)\sm K\RP\RP\,d\nu\le\eps
\end{equation}
\end{Asm}
Using CHG, one can transform the left-hand side of (\ref{eq:RnwlK}):
take the $\om$-dependent graph
\BER{l}
\LB(x,y)\in G^2: T^x\om,T^y\om\in E,\,x-y-a\in P(T^y\om)\sm K\RB =\\=
\LB(x,y)\in G^2: T^x\om\in E, y-x+a=\pi_E\LP T^{x-a}\om\RP\notin-K\RB
\EER
and the $1$-simplex of measures $(\cnt,\cnt)$.

The source is the enhanced \fn\
$$1_{\{\om\in E:\pi(T^{-a}\om)\notin-K\}}\frac{d\cnt}{dx};$$
the target is the enhanced \fn\ 
$$\#\om\cap\LP a+\LP P(\om)\sm K\RP\RP\frac{d\cnt}{dx};$$
and by CHG we can write (\ref{eq:RnwlC}) in the equivalent form
\begin{equation}\label{eq:RnwlC1}
\nu\{\om\in E:\pi(T^{-a}\om)\notin-K\}<\eps
\end{equation}

Let us show that if $G=\bR$ or $\bZ$ and the $X_n$ are independent, and
$\mu(\OM)=\int_EX_1(\om)\,d\nu(\om)<\I$, and one
takes $P(\om)=[0,X_1(\om)[$ then Assumption \ref{Asm:C} holds, where we use
the form (\ref{eq:RnwlC1}). Indeed, what we have to prove amounts to showing
that the $\nu$-probability that the last renewal before time $-a$ was
even before $-a-b$ (i.e.\ that $\om\cap[-a-b,-a[=\es$),
tends to $0$ when $b\to+\I$ uniformly in $a$ for $a$ near $+\I$ (we take
$a>0$). But by independence of the $X_n$,
\BER{l}
\Pr(\om\cap[-a-b,-a[=\es)\le\sum_{k\in\bZ^+}\sum_{n\in\bZ^+}\Pr
\LP S_{-n}(\om)\in[-a+k,-a+k+1[\:\land\:X_{-(n+1)}>b+k\RP\le\\ \le
\sum_{k\in\bZ^+}\LP\bE\#\om\cap[-a+k,-a+k+1[\RP\Pr(X_1>b+k)\le
\mbox{const}\cdot\sum_{k\in\bZ^+}\Pr(X_1>b+k)\to_{b\to+\I}0
\EER
where in the last inequality we used Assumption \ref{Asm:B} to get the bound
$\mbox{const}$ and also the fact that $X_1$ is integrable.

Thus we finally have:

\begin{equation} \label{eq:RnwlI1}
\laa*(x\mapsto\nu\{\om:x\in P(\om)\})=1\qquad dx-\mbox{a.e.} 
\end{equation}

Note that the \fn\ on $G$\quad$x\mapsto\nu\{\om:x\in P(\om)\}$ has
$dx$-integral equal to $\int_E|P(\om)|\,d\nu(\om)=\mu(\OM)<\I$.
Suppose the following assumption holds:

\begin{Asm} \label{Asm:I1}
The fourier transform of the nonnegative $L^1$ \fn\ on $G$\quad
$x\mapsto\nu\{\om:x\in P(\om)$ never vanishes.
\end{Asm}

Then this assumption together with (\ref{eq:RnwlI1}) would imply that
$\laa$ is $dx$ multiplied by the reciprocal of the value at $0$ of the
Fourier transform of $f(x)=\nu\{\om:x\in P(\om)\}$, namely by the reciprocal
of $\int_G\nu\{\om:x\in P(\om)\}\,dx=\mu(\OM)$, and we would be done.
This implication obtains as follows:

We wish to infer from the fact that $\laa$ and $(\mu(\OM))^{-1}\,dx$ have
the same convolution with our $f\in L^1(G,dx)$,
whose Fourier transform $\hat{f}(t)$ never vanishes,
that $\laa=(\mu(\OM))^{-1}\,dx$. Now, $\laa$ can be weakly
approximated by its convolutions with \fns\ in $\cC_{00}$, which are bounded
\cts\ \fns\ on $G$ (Recall that $\laa$ is uniformly bounded on the translates
of any fixed compact $\st G$). Thus, it suffices to prove that if $\ell(x)$ is
\cts\ bounded on $G$ and $\ell*f=1$ then $\ell\equiv(\mu(\OM))^{-1}\cdot1$.
This follows from the fact that the closed ideal in $L^1(G)$ (where
the multiplication is convolution and we take the norm topology) generated
by $f$ is the whole $L^1(G)$ (this is Wiener's Tauberian Theorem).
This follows from Fourier transform considerations:
Indeed, the said ideal contains the \fns\ $k(x)\in L^1(G)$ whose fourier
transforms are of the form $\hat{k_1}(t)\hat{f}(t)$ where $k_1\in L^1(G)$
with $\hat{k_1}\in\cC_{00}(\hat{G})$.
Since $\hat{f}$ never vanishes, every \fn\ in $L^1(G)$ with Fourier transform
in $\cC_{00}(\hat{G})$ can be approximated by such $k$, hence belongs to the
ideal and the latter \fns\ are dense in $L^1(G)$.

Thus it remains to ensure that Assumption \ref{Asm:I1} holds when we take
$G=\bR$ or $\bZ$ and $P(\om)=[0,X_1(\om)[$. But then 
$\nu\{x\in P(\om)\}=\Pr\{X_1(\om)>x\ge0\}$ which is $0$ for $x<0$ and
is nonincreasing nonnegative for $x\ge0$. If $t$ is a nonzero element
of the dual group $\hat{G}$, we have, if $G=\bR$, $\hat{G}=\bR$ and
$dx$ is the Lebesgue measure:
$$
\int_0^\I\exp(2\pi ixt)\Pr\{X_1>x\}\,dx=
(2\pi it)^{-1}\int_0^\I\LP 1-\exp(2\pi ixt)\RP\,d\Pr\{X_1>x\},
$$
and if $G=\bZ$, $\hat{G}=\bR/\bZ$ and $dx$ is the counting measure:
$$
\sum_{0\le x<\I}\exp(2\pi ixt)\Pr\{X_1>x\}=
\LP 1-exp(2\pi t)\RP^{-1}\sum_{0\le x<\I}\LP 1-exp(2\pi ixt)\RP
\LP\Pr\{X_1>x-1\}-\Pr\{X_1>x\}\RP.
$$
Therefore the Fourier transform of $x\mapsto\Pr\{X_1>x\}$ can vanish
at $t$ only if the distribution of $X_1$ is concentrated in $t^{-1}\bZ$.
Thus, in the nonperiodic case, i.e.\ when there is no proper closed
subgroup $H\st G$ s.t.\ a.s.\ $X_1\in H$, Assumption \ref{Asm:I1} holds
and we have the renewal limit assertion.
}

\newpage
\appendix
\section{Appendices} \label{S:AP}
\subsection{Generation of Measures via Given ``Preintegrable'' Functions}
\label{SS:PRE}
We describe a way to obtain a measure on a set $\OM$,
\NOT{which will be used in the sequel.}%
which we use in \S\ref{SS:EnhInf}. This method seems convenient when a
measure has to be constructed by some ``integration'' of a family of given
measures.

The starting point is a set $\Pre$ of $[0,\I]$-valued \fns\ on $\OM$, called
{\EM preintegrable}, with a functional ({\EM integral}) $\cI:\Pre\to[0,\I[$,
satisfying the following assumptions:
 
\begin{enumerate}
\item 
$\Pre$ is a cone, i.e.\ $\Pre$ contains $0$ and is stable w.r.t.\ addition
and multiplication by finite nonnegative real constants, and $\cI$ is
additive and non-negatively linear.

\item 
If $f,g:\OM\to\BAR{\bR^+}$ s.t.\ $f,f+g\in\Pre$ then $g\in\Pre$.
(Consequently, $f,g\in\Pre,\:f\ge g\Rightarrow\cI(f)\ge\cI(g)$.)

\item 
If $f_n\in\Pre,\:n\in\bN$, $f_n\uparrow$ and $\cI(f_n)$ is bounded, then
$\lim f_n\in\Pre$ and $\cI(\lim f_n)=\lim(\cI(f_n))$.

\item \label{it:comp} 
If $f_0\in\Pre$ and $\cI(f_0)=0$ then any \fn\ $f\le f_0$ belongs to $\Pre$.

\end{enumerate}

Now say that a set $E\subset\OM$ is {\EM measurable} if
$$\forall f\in\Pre\; f\cdot 1_E\in\Pre.$$

One proves easily, using the above assumptions, that the measurable sets
form a $\sigma$-algebra. {\EM Measurable \fns} will be \fns\ measurable
w.r.t.\ this $\sigma$-algebra. Note that if $g$ is measurable and bounded,
then
$$\forall f\in\Pre\; f\cdot g\in\Pre$$.

To define the measure $\mu$ on this $\sigma$-algebra, a measurable set $E$
will have finite measure iff $1_E\in\Pre$ and then $\mu(E):=\cI(1_E)$.
Otherwise $\mu(E)=\I$. The assumptions on $\Pre$ and $\cI$ imply readily
that $\mu$ is $\sigma$-additive.
(Moreover, by \ref{it:comp}.\ $\mu$ is complete, i.e.\ every subset of a set
of measure $0$ is measurable.)
Thus $\int\:d\mu$ is defined. As usual, a \fn\ $f\ge0$
is {\EM integrable} if it is measurable and has finite integral.

An important fact is that any integrable \fn\ $f\ge0$ is preintegrable, and
any measurable preintegrable \fn\ $f\ge0$ is integrable, and then $\cI(f)$
and $\int f\,d\mu$ coincide. (Thus, to find the integral of a measurable
\fn\ $f\ge0$ one just checks if $f$ is in $\Pre$. If it is, its integral is
$\cI(f)$, otherwise $\int f\,d\mu=\I$). 

Indeed, Note first that any $\{0,\I\}$-valued $f\in\Pre$ must have
$\cI(f)=0$, since $\cI(f)$ is finite and $2\cI(f)=\cI(2f)=\cI(f)$.
Therefore if $f\ge0$ is measurable preintegrable, then $\{f=+\I\}$ is
null.  
Now an integrable $f\ge0$ can be obtained from characteristic \fns\ of
sets of finite measure by addition and increasing limits with bounded
integral, hence it is in $\Pre$. On the other hand, if $f\ge0$ is measurable 
preintegrable, and $a>0$, then $E=\{\I>f>a\}$ is measurable. $1_E$ is of the
form $f\cdot g$ where $g$ is measurable bounded, hence $1_E$ is 
preintegrable, implying $\mu(E)=\cI(1_E)<\I$. From that one easily deduces
$f$ integrable and $\int f\,d\mu=\cI(f)$.

\subsection{Standard Borel Spaces and Products of Two Standard Spaces}
\label{SS:PRD}
Recall that a {\EM standard Borel space} is a Borel space which is isomorphic,
as a Borel space, to a Lusin topological space (recall that in any
topological space the Borel structure understood is the $\sigma$-algebra of
``ordinary'' Borel subsets). We shall use facts about Lusin and
Polish spaces -- see  \cite{BourbakiTOP} Ch. IX \S6, \cite{KechrisD}, 
\cite{Kuratowski} (where the terminology is slightly different).
Recall, in particular, that topological spaces where the topology can be
given by a complete saparable metric are called {\EM Polish spaces};
Topological spaces which are \cts\ 1-1 images of Polish spaces are called
{\EM Lusin spaces};
a subset of a Lusin space is Borel iff it is Lusin in its relative topology;
a subset of a Polish space is $G_\delta$ iff it is Polish in its
relative topology;
for any 1-1 Borel mapping between standard Borel spaces the image is Borel
and the mapping is an isomorphism of the Borel structures (with the image);
any Lusin space is a 1-1 continuous image of a $0$-dimensional polish
space, i.e.\ a Polish space with a 
base to the topology consisting of clopens;
any Polish space can be \cts{ly} embedded in a metric
compact space, the latter can be chosen $0$-dimensional if the former is.

In fact, two standard spaces of the same cardinality are isomorphic as
Borel spaces. Thus the only isomorphism types of standard Borel spaces are:
finite sets, the type of a countable set with the $\sigma$-algebra of all
subsets, and the unique type of a standard space of the cardinality of
the continuum.
Thus for many purposes one may assume the latter is the unit interval
with Borel subsets. However, for the purposes below it is preferable
to consider the totality of all Lusin topologies in the standard space,
having in mind, of course, topologies s.t.\ their Borel structure is
the given one.

In this vein, one notes that for every countable Boolean algebra of subsets
(of a standard Borel space) which separates points,
the obvious mapping to $2^\bN$ is 1-1 Borel, hence a Borel isomorphism
with the image, which is Borel in $2^\bN$ hence Lusin in the relative
topology. This implies that for every countable set of Borel sets there is
a Polish topology where all are clopen, consequently for every countable
collection of Borel bounded real-valued \fns\ one may find a Polish topology
where all are continuous, and for every countable collection of
$[0,\I]$-valued Borel \fns\ there is a Polish topology where all are
l.s.c.\ (lower semi-continuous).

Note that by using the diagonal in a countable product, one proves that
for every countable family of Lusin (resp.\ Polish, resp.\ $0$-dimensional
Polish) topologies there is a topology of the same kind finer than all of
them.

Thus when one is confronted with, say, a non-negative Borel \fn\ on a
standard space, one may assume that it is l.s.c.\ for some Polish topology
there.

Matters are not so simple if one deals with {\em a product of two standard
spaces $X$ and $Y$} and one may choose topologies in $X$ and $Y$, but
in the product one always take the {\em product topology}.

\begin{Prop} \label{Prop:ProdSet}
Let $X$ and $Y$ be standard Borel spaces. Let $E\st X\times Y$ be Borel.
T.f.a.e:
\begin{itemize}
\item[(i)]
$E$ is a disjoint union of countably many ``Borel rectangles'':
products of Borel sets in $X$ and $Y$
\item[(ii)]
$E$ is open in some product of Lusin topologies in $X$ and $Y$.
\end{itemize}
\end{Prop}

\begin{Prf}
{\EM (i) $\Rightarrow$ (ii)}: take in $X$ and $Y$ topologies making all
sides of the rectangles open.
\par\medskip

{\EM (ii) $\Rightarrow$ (i)}: Since every Lusin space is a \cts\ 1-1 image of
a Polish space, one may assume the topologies are Polish. Choose countable
bases to the topologies and consider the countable Boolean algebras generated
by the bases. $E$ is a countable union of ``rectangles'' with sides belonging
to the Boolean algebras, hence a countable disjoint union of such.
\end{Prf}

As an example of a set which does not satisfy (i) and (ii) in the previous
proposition, take the diagonal in $[0,1]\times[0,1]$.

\begin{Prop} \label{Prop:ProdFn}
Let $X$ and $Y$ be standard Borel spaces. Let $f:X\times Y\to[0,\I]$ be
Borel. T.f.a.e:
\begin{itemize}
\item[(i)]
$f$ can be represented as a series:
$$f(x,y)=\sum_{i\ge1}g_i(x)h_j(y)\quad x\in X,\:y\in Y$$
where $g_i:X\to[0,\I]$, $h_i:Y\to[0,\I]$ are Borel.
\item[(ii)]
$f$ is l.s.c.\ (lower semi-continuous) for the product of some Lusin
topologies in $X$ and $Y$.
\item[(iii)]
Every set $\{f>a\}$, $a\in[0,\I[$ satisfies the requirements of the previous
proposition.
\end{itemize}
\end{Prop}

\begin{Prf}
{\EM (ii) $\Rightarrow$ (iii)} is obvious.
\par\medskip

{\EM (i) $\Rightarrow$ (ii)}: there are Lusin topologies in $X$ and $Y$
making all $g_i$ and $h_i$ l.s.c., thus making $f$ l.s.c.
\par\medskip

{\EM (iii) $\Rightarrow$ (i)}: take all sets $\{f>a\}$ for {\em dyadic}
$a$, describe them as disjoint unions of ``Borel rectangles'' and 
consider the countable Boolean algebras in $X$ and $Y$ generated by all
their sides. $f$ is a supremum of countable positive combinations of
characteristic
\fns\ of countable unions of rectangles with sides in the Boolean algebras,
hence a supremum of countable positive combinations of characteristic
\fns\ of single rectangles, and since the latter combinations are stable
w.r.t.\ lattice operations and subtraction, $f$ is a sum of a series of
positive multiples of characteristic \fns\ of rectangles. 
\end{Prf}

Again, the characteristic \fn\ of the diagonal in $[0,1]\times[0,1]$ does
not satisfy the requirements of Prop.\ \ref{Prop:ProdFn}.

\subsection{Converting Measurable Action to Continuous Action} 
\label{SS:CP}
Our setting is a 2nd-countable locally compact group $G$ acting in a Borel
manner on a standard Borel space $(\OM,\cB)$,
thus making it into a {\em standard $G$-space}.

A special case of the above is a {\EM metrizable compact $G$-space}, where
one takes usual Borel sets and $(x,\om)\mapsto x\om$ is assumed
{\em \cts\ in the two variables}.

There is a well-known method (\cite{Varadarajan}, see also, e.g.\
\cite{AmbroseKakutani}, \cite{Doob}, \cite{Mackey} where the idea of
mapping $\om$ to the function on $G$ $\cO_\om f$ is employed)
to embed any standard $G$-space as a Borel subset of a metrizable compact
$G$-space (this is done in a definitely {\em non-unique} way).

Choose any countable set $\cF$ of Borel \fns\ $f$ with $|f|\le1$,
separating points in $\OM$ (this exists by standardness). For any
$\om\in\OM$ and any $f\in\cF$ we have
$$\cO_\om f:= x\mapsto f(x\om), \mbox{ thus } \cO_\om f:G\to\bR$$
This may be considered as an element of the unit ball of $\cB(L^\I(G))$,
the latter taken w.r.t.\ (right or left) Haar measure, and is endowed with
the $w^*$-topology from $L^1$. This unit ball is metrizable compact and the
mapping
$$\om\mapsto\LP\cO_\om f\RP_{f\in\cF}$$
maps $\OM$ into the compact metrizable $K=\cB(L^\I(G))^\cF$, mapping the
$G$-action into right translation in any coordinate $\cB(L^\I(G))$.

Now we can verify some facts:

The action of $G$ by right translation in any $\cB(L^\I(G))$ (hence in a
power $\cB(L^\I(G))^\cF$) is \cts\ in the two variables.

Our map from $\OM$ into the power is Borel. Indeed, the Borel structure
in $\cB(L^\I(G))$ is defined by some countable set of ``coordinates''
$h\in L^1(G)$, and for such $h$, the function ($\la$ is some Haar measure)
$$\om\mapsto\LA\cO_\om f,h\RA=\int f(x\om)h(x)\,d\la$$
is Borel.

This map is 1-1. Indeed, suppose $\om_1$ and $\om_2$ map to the same
$$\LP\cO_{\om_1} f\RP_{f\in\cF}=\LP\cO_{\om_2} f\RP_{f\in\cF}.$$
(where equality of the $\cO_\om f$ means equality in $L^\I$, that is
equality a.e.\ w.r.t.\ Haar). This means that for a.a.\ $x\in G$
$\forall f\in\cF\:(\cO_{\om_1}f)(x)=(\cO_{\om_2}f)(x)$, i.e.
$\forall f\in\cF\:f(x\om_1)=f(x\om_2)$. Hence this holds for {\em some} $x$,
which means, since $\cF$ separates points, that $x\om_1=x\om_2$, implying
$\om_1=x^{-1}x\om_1=x^{-1}x\om_2=\om_2$.

Since the map is 1-1 Borel between standard spaces, its image is a Borel
subset of the metrizable compact $K$ (see \S\ref{SS:PRD}).

Instead of the above $K$, we may and do take as our $K$ the closure of
the image of $\OM$.

Note that the relative topology from $K$ is a Lusin topology in $\OM$
itself (see \S\ref{SS:PRD}) making the $G$-action \cts\ in the two
variables.

Note that if $f\in\cF$ and one takes a convolution of $f$ with an 
$L^1$-\fn on $G$ (i.e.\ one considers $\om\mapsto\int_G f(x\om)\,d\nu(x)$,
when $\nu$ is a finite measure absolutely \cts\ w.r.t.\ Haar), then the
latter (extends to) a \cts\ \fn\ on $K$. If $G$ is discrete one may take
$\nu=\delta_{e}$, and any $f\in\cF$ extends to a \cts\ \fn\ on $K$. Thus
{\em for discrete $G$} $K$ may be tuned so that any countably many given
bounded Borel \fns\ extend to \cts\ \fns\ on $K$, hence so that any
countably many given $[0,\I]$-valued Borel \fns\ extend to l.s.c.\
(lower semi-\cts) on $K$ (here the extension need not be unique).

If $\OM$ was from the start a dense invariant Borel subset of a
metrizable compact $G$-space $K_0$, we can take all the $f\in\cF$ \cts\ on
$K_0$, and then we have a mapping from $K_0$ into the power
$K=\cB(L^\I(G))^\cF$
which is 1-1 continuous, hence an isomorphism with a compact sub-$G$-space.
This shows that by this construction one gets all embeddings of $\OM$ as
a dense subset of a metrizable compact $G$-space up to isomorphism.

These metrizable compact $G$-spaces in which we embedded a standard
Borel $G$-space are not canonical, since a countable set of $f$'s has
to be chosen. One may get a {\em canonical, but not metrizable} 
$G$-compact by taking the set $\cF$ of {\em all} Borel \fns\ $f$ with
$|f|\le1$, and considering the closure $\cK(\OM)$ of the image of $\OM$ in
the power $\cB(L^\I(G))^\cF$. The Borel sets in $\OM$ will be
the intersection of $\OM$ with the {\em Baire} sets in $\cK$.

For $G$ countable discrete, $\cK$ does not depend on $G$ and
is just the Stone space of the $\sigma$-algebra of Borel sets in
$\OM$. In particular, in this case every bounded Borel \fn\ on $\OM$
extends to a \cts\ \fn\ on $\cK$, and every unbounded non-negative Borel
\fn\ extends to a b.l.s.c.\ \fn\ on $\cK$ (see \S\ref{SS:ULSC}).

For general $G$ this is not the case: a necessary condition for a Borel
$f$ on $\OM$ to extend to a \cts\ \fn\ on $\cK$ is that $\cO_\om f$ is
\cts\ on $G$ for all $\om\in\OM$. On the other hand, any convolution of
a Borel \fn\ on $\OM$ with some $L^1$ \fn\ on $G$ does extend to a \cts\
\fn\ on $\cK$.

\subsection{Mean Ergodic Theorems for General (Discrete) Groups}
\label{SS:ERG}
The aim of this \S\ is to give a treatment, in the spirit of \S\ref{SS:Av},
of the well-known derivation of mean ergodic theorems for general groups,
using weak compactness and Ryll-Nardzewski's fixed point Theorem.
(see \cite{BirkhoffAlaoglu}, \cite{GlicksbergDeLeeuwA},
\cite{GlicksbergDeLeeuwD}, \cite{Greenleaf}, \cite{Jacobs} \S2,
\cite{TroallicE}, \cite{TroallicS}). 
{
\newcommand{\toAA}{{\buildrel\mbox{AA}\over\to}}
\newcommand{\limAA}{{\lim\atop\mbox{AA}}}

\NOT{
\begin{Thm}
Let $V$ be a {\EM $G$-normed space}, i.e.\ $V$ is a $G$-vector space
(that is: $G$ acts on $V$ linearly),
normed by an {\em invariant} norm $\|\|$.

\end{Thm}
}%

$G$ is a (discrete) group, in general {\em non-amenable}.

Consider a {\EM $G$-normed space} $V$, i.e.\
$V$ is a $G$-vector space (that is: $G$ acts on $V$ linearly),
normed by an invariant norm $\|\|$. Recall the definition of averages
in \S\ref{SS:Av}

\begin{Def} \label{Def:AA} (cf.\ \cite{BirkhoffAlaoglu})
We say that a $v\in V$ converges to a $v_0\in V$ in the {\EM Accumulating
Averages (AA)} sense or that $v_0$ is the AA-limit of $v$
(denote: $v \toAA v_0$ or $v_0= \limAA v$),
if {\em $v_0$ is $G$-invariant} and

\begin{center}
for every average $v'$ of $v$, every neighbourhood of $v_0$ contains an
average of $v'$.
\end{center}
\end{Def}

It is straightforward that the AA-limit of $v$ is {\em unique}.

Also, $v\toAA v_0\Leftrightarrow v-v_0\toAA 0$.

We have: if $\limAA v$ exists, then for any $G$-invariant $v^*\in V^*$,
$\LA v^*,v\RA=\LA v^*,\limAA v\RA$.

\NOT{
Note that
$$v\toAA v_0\Rightarrow vx^{-1}\toAA v_0\Rightarrow v\toAA xv_0
\Rightarrow xv_0=v_0$$
Thus, {an AA-limit \em $v_0$ must be invariant}.
}

We use the following notation:
$\cA\cV(v)$ is the set of averages of a vector $v$,
\BER{c}
\|v\|_A=\inf_{w\in\cA\cV(v)}\|w\| \\
\|v\|_{AA}=\sup_{w\in\cA\cV(v)}\|w\|_A
\EER
Thus, $v\toAA v_0\Leftrightarrow \|v-v_0\|_{AA}=0$. 

Clearly, $\|v\|_A\le\|v\|_{AA}$.

If $\limAA v$ exists then $\|\limAA v\|\le\|v\|_A$

\begin{Prop} \label{Prop:AAIneq}
\parsk

\begin{itemize}
\item[a.]
$\|v+v'\|_A\le\|v\|_{AA}+\|v'\|_A$
\item[b.]
$\|v+v'\|_{AA}\le\|v\|_{AA}+\|v'\|_{AA}$
\end{itemize}
\end{Prop}

\begin{Prf}
Call an operator of the form $v\mapsto\sum_x\la_xxv$ (finite sum) where
$\forall x\:\la_x\ge0\:\sum_x\la_x=1$ an {\em averaging operator (a.o.)}.
These operators form a semigroup.
In the rest of the proof, $L$, $L'$ and $L''$ will refer to a.o.'s.
We have
$\forall L, \|Lv\|\le\|v\|$
and
$\forall L,\|L(v+v')\|\le\|L(v)\|+\|L(v')\|$

Proof of a.:  

Choose $L$ with $\|Lv\|\le\|v\|_A+\eps$. Choose $L'$ with
$\|L'Lv'\|\le\|v'\|_{AA}+\eps$. Then 
$\|L'L(v+v')\|\le\|v\|_A+\|v'\|_{AA}+2\eps$.
\par\medskip

Proof of b.:

$\forall$ $L''$, Choose $L$ with $\|LL''v\|\le\|v\|_{AA}+\eps$.
Choose $L'$ with
$\|L'LL''v'\|\le\|v'\|_{AA}+\eps$. Then 
$\|L'LL''(v+v')\|\le\|v\|_{AA}+\|v'\|_{AA}+2\eps$.
\end{Prf}

Thus, $\|\|_{AA}$ is a semi-norm dominated by $\|\|$.

This implies that 
AA-convergence has the desired properties:
\NOT{there can be only one AA-limit, and}
if $v\toAA v_0$ and $w\toAA w_0$ then
$\alpha v+\beta w\toAA\alpha v_0+\beta w_0$.

\NOT{
(Given an average of $v+w$, choose an further average $\tilde{v}$, ``near''
$v_0$, of the corresponding average of $v$, then choose a ''near'' further
average of $w$ of the average corresponding to $\tilde{v}$.
Here the invariance of the norm is used:
if an average is ``near'', further averages will stay ``near''.).
}
That much cannot be said for the property: every neighbourhood of $v_0$
contains an average of $v$ (i.e.\ $\|v-v_0\|_A=0$).
Note that while for abelian $G$ the averaging operators commute, hence
any two averages have a common average, which implies immediately
$\|v\|_{AA}=\|v\|_A$,
for non-abelian groups two averages need
not have a common average.
\NOT{
$v$ has averages arbitrarily close to $v_0$ $\Rightarrow$ $v\toAA
v_0$.
}

\begin{Exm} 
$G=$ the infinite dihedral group $D_\I$ (which is {\em amenable}),
realized as the set of the transformations
of $\bR$ generated by $x_1:t\to t+1$ and $x_{-}:t\to -t$.
$V=$ the space of polynomials of degree $\le3$. 
$v=t(t-1)(t+1)$. $x_{-}v=-v$ so
$w=\frac{1}{2}(v(t+1)-v)=\frac{3}{2}t(t-1)$ and $-w=\frac{1}{2}(-v(t+1)+v)$
are averages of $v$ which have no common average (any average of $w$ has
coefficient $\frac{3}{2}$ at $t^2$).
\end{Exm}

\begin{Exm} \label{Exm:NonAme} 
Let $G$ be a non-amenable group. Let $K$ be the set of averages on
$\ell^\I(G)$, i.e.\ the set of positive functionals on $\ell^\I(G)$
giving the value $1$ to the constant sequence $1$. In $\OM$, take the
$w^*$-topology from $\ell^\I$. $K$ is a $G$-convex compact space.
Since $G$ is non-amenable, there is no $G$-fixed point in $K$, hence
there is a minimal non-empty convex subset $K_0\subset K$. $K_0$ 
is not a singleton. It satisfies: for every $\om\in K_0$, the closed
convex hull of its orbit is the whole $K_0$ (if $G$ is countable,
$K_0$ has a compact {\em metrizable} factor with the same property).

Let $f:K_0\to\bR$ be non-constant \cts\ affine with minimum $0$ and
maximum $a>0$. By Thm.\ \ref{Thm:MaxAv} b.\ (for $\cC(K_0)$ and $p=\max$),
$f$ has averages with maximum arbitrarily close to $0$, and averages with
minimum arbitrarily close to $a$ (a property inherited by all their
averages). Thus $\|f\|_{\cC(K_0),A}=0$ while $\|f\|_{\cC(K_0),AA}=a$. 

Note that by Thm.\ \ref{Thm:AAFixed}, for {\em amenable $G$} we have always
$\|v\|_A=\|v\|_{AA}$.
\end{Exm}

\begin{Rmk} \label{Rmk:AABanach}
We have seen that $v\to\limAA v$ is a linear operator,
defined on a linear subspace of
$V$ and is norm-continuous there. Also it is clearly a {\em closed} operator.
Thus in case $V$ is a {\em Banach} space, its domain of definition
(i.e.\ the set of vectors having AA-limit) is {\em closed}.
\end{Rmk}

\begin{Prop} \label{Prop:AACnvx}
Let $B^*$ be the unit ball of $V^*$ with the $w^*$-topology. Let $v\in V$.
Then:
\begin{itemize}
\item[a.]
$\|v\|_A$ equals the maximum over $G$-invariant $K\subset B^*$
(or, if one wishes,
 over convex $w^*$-compact $G$-invariant $K\subset B^*$) of
$\inf_{v^*\in K}\LA v,v^*\RA$.
\item[b.]
$\|v\|_{AA}$ equals the supremum over $w^*$-convex compact invariant
$K\subset B^*$ of the minimum over $G$-invariant $K'\subset K$
(or, if one wishes,
 over minimal convex $w^*$-compact $G$-invariant $K'\subset K$)
of $\sup_{v^*\in K'}\LA v,v^*\RA$.
\end{itemize}
\end{Prop}

\begin{Prf}
$K$, $K'$ will refer to $G$-invariant convex $w^*$-compact subsets of $B^*$.
\begin{itemize}
\item[a.]
Follows from Thm.\ \ref{Thm:MaxAv} b.\ for $p=\max$.
\item[b.]
By a.,
\BER{l}
\|v\|_{AA}=\sup_{w\in\cA\cV(v)}\|w\|_A=\\
=\sup_{w\in\cA\cV(v)}\sup_K\min_{v^*\in K}\LA w,v^*\RA=\\
=\sup_K\sup_{w\in\cA\cV(v)}\min_{v^*\in K}\LA w,v^*\RA=\\
=\sup_K\min_{K'\subset K}\max_{v^*\in K'}\LA v,v^*\RA
\EER
The last equality following from Thm.\ \ref{Thm:MaxAv} b.\ (for $p=\max$)
applied to $-v$ instead of $v$.
\end{itemize}
\end{Prf}

For the next theorem, we use the Ryll-Nardzewski Fixed-Point Theorem
(see \cite{BourbakiEVT} Ch.~IV App.\ for a proof):

\begin{Main}
{\EM Ryll-Nardzewski Fixed-Point Theorem}: Let $V$ be a normed space,
and $K$ a convex non-empty weakly-compact subset of $V$. Let $G$ be
a group on affine norm-isometries of $K$. Then $G$ has a fixed point in
$K$
\end{Main}

\begin{Thm} \label{Thm:AAFixed}
Let $V$ be a $G$-normed space, i.e.\ a normed space on which a group $G$
acts linearly isometrically.
Suppose either {\em $V$ is reflexive} or {\em $G$ is amenable}.

Then for every $v\in V$, $\|v\|_A=\|v\|_{AA}$ $=$ the maximum of
$|\LA v^*,v\RA|$ over {\em $G$-invariant} $v^*\in V^*$ with norm $\le1$.

Consequently,\par
$v\toAA 0$ $\Leftrightarrow$ \par
$\Leftrightarrow$
$v$ has averages with arbitrarily small norm
$\Leftrightarrow$ \par
$\Leftrightarrow$
$v$ is annulled by all $G$-invariant $v^*\in V^*$.
\end{Thm}

\begin{Prf}
The two assumptions: {\em $V$ is reflexive} or {\em $G$ is amenable},
have in common that they imply that in every bounded $w^*$-closed convex 
subset $K\subset V^*$ there is a $G$-fixed point. (for the case of
$V$ reflexive this follows from Ryll-Nardzewski's Thm.) Hence every minimal
convex $w^*$-compact $G$-invariant subset of $V^*$ is a singleton.
Having said this, the assertion of the theorem follows from
Prop.\ \ref{Prop:AACnvx} (note that the maximum of $|\LA v^*,v\RA|$
over $G$-invariant $v^*$'s in the unit ball of $V^*$
equals the maximum of $\LA v^*,v\RA$). 
\end{Prf}


Let $G$ act measure-preservingly on a probability space $(\OM,\cB,\mu)$.

The \fn\ spaces $L^p(\OM)$, $1\le p\le\I$ are $G$-normed spaces.
For $1<p<\I$ they are reflexive and Thm.\ \ref{Thm:AAFixed} applies.
Hence every $L^p$-\fn\ with conditional expectation $0$ w.r.t.\
the Boolean-$\sigma$-algebra of almost-invariant%
\footnote{see ``Notations''}%
subsets, in other words, any
$L^p$-\fn\ annulled by the $G$-invariant members of the dual space,
converges AA in $L^p$ norm to $0$.

$G$-Invariant members of $L^p$ converge AA to themselves.

Now, any continuous functional $\phi\in L^{p'}$
annulled both by the \fns\ with conditional expectation $0$ and by the
invariant \fns\ must be $0$, and the set of members of a $G$-Banach space
with AA-limit is closed (Remark \ref{Rmk:AABanach}). Hence in our case it is 
the whole space and by closeness we are allowed to say the same about
$L^1$. 

To conclude, we have the following ``mean ergodic theorem'', valid for
{\em any group $G$} and not referring to F{\o}lner sequences.

\begin{Thm} \label{Thm:AAErg} 
Let $G$ act measure-preservingly on a probability space $(\OM,\cB,\mu)$.
Let $1\le p<\I$.

Then every $f\in L^p$ converges in the AA (accumulating
averages) sense in $L^p$-norm to
an invariant $f_0\in L^p$. $f_0$ is the conditional expectation of
$f$ w.r.t.\ the $\sigma$-algebra of almost-invariant%
\footnote{see ``Notations''}%
subsets. (in the $L^2$ case, $f_0$ is the orthogonal projection of $f$ on
the space of $G$-invariant \fns). In other words, every (finite) average of
$f$ has (finite) averages arbitrarily close to $f_0$ in $L^p$ norm.
\end{Thm}
\qed

}


\newpage
\begin{thebibliography}{WW-W}
\addcontentsline{toc}{section}{References}

\bibitem[AK]{AmbroseKakutani}
Ambrose~W., Kakutani~S., {\em Structure and continuity of measurable
flows}, Duke Math.\ J.\ {\bf 9} 25-42 (1942).

\bibitem[AW]{AaronsonWeiss}
Aaronson J., Weiss B.,
{\em A $\bZ^d$ ergodic theorem with large normalising constants},
in: Convergence in Ergodic Theory and Probability, de Gruyter, 1996.%
\NOT{preprint (1994).}

\bibitem[Bi]{Birkhoff}
Birkhoff~G.~D., {\em Proof of a recurrence theorem for strongly transitive
systems}, Proc.\ Nat.\ Acad.\ Sci.\ {\bf 17}, No.~12, 650-660. 
Birkhoff: Collected Mathematical Papers, Vol.~2, 398-408. (1931).

\bibitem[BA]{BirkhoffAlaoglu}
Birkhoff~G., Alaoglu~L., {\em General ergodic theorems}, Annals of Math.\
{\bf 41} No.~2, 293-309 (1940).

\bibitem[Bl]{Blanchard}
Blanchard~F., {\em K-flots et th\'eor\`eme de renouvellment}.,
Z.\ Wahrsch.\ verw.\ Geb.\ {\bf 36}, 345-358 (1976).

\bibitem[Bo-A]{BourbakiALG}
Bourbaki~N., {\em Alg\`ebre}. Ch. 4 \`a 7, Masson, Paris 1981.

\bibitem[Bo-E]{BourbakiEVT}
Bourbaki~N., {\em Espaces Vectoriels Topologiques}. Masson, Paris, 1981. 

\bibitem[Bo-I]{BourbakiINT}
Bourbaki~N., {\em Int\'egration}. Hermann, Paris, 1965. 

\bibitem[Bo-T]{BourbakiTOP}
Bourbaki~N., {\em Topologie G\'en\'erale}. Ch. 5 \`a 10,
Diffusion C.C.L.S.\ Paris, 1974.

\bibitem[Br]{Breiman}
Breiman~L., {\em Probability}. SIAM, Philadelphia, 1992.

\bibitem[Co]{Cohn}
Cohn~P.~M., {\em Lie Groups}. Cambridge University Press, 1957.

\bibitem[DV]{DaleyVereJones}
Daley~D.~J., Vere-Jones~D., {\em An Introduction to the Theory of Point
Processes}. Springer-Verlag, 1988.

\bibitem[De]{Delasnerie}
Delasnerie~M., {\em Flot m\'elangeant et mesures de Palm},
Ann.\ Inst.\ Henri Poincar\'e {\bf XIII} No.\ 4, 357-369 (1977). 

\bibitem[Do]{Doob}
Doob~J.~L., {\em One-parameter families of transformations}, Duke Math.\
J.\ {\bf 4} 752-774 (1938).

\bibitem[Du]{Durrett}
Durrett~R., {\em Probability: Theory and Examples}. 2nd Ed., 
Duxbury Press at Wadsworth Publishing Co., 1996.

\bibitem[El]{Ellis}
Ellis~R., {\em Topological dynamics and ergodic theory}, Ergod.\ Th.\ \&
Dynam.\ Sys.\ {\bf 7}, 25-47 (1987).

\bibitem[Fed]{Federer}
Federer~H., {\em Geometric Measure Theory}. Springer-Verlag, 1969.

\bibitem[Fel]{Feller}
Feller~W., {\em An Introduction to Probability Theory and its Applications}.
3rd ed., Vol.\ I, John Wiley \& Sons, 1968.

\bibitem[Fl]{Flanders}
Flanders~H., {\em On spaces of linear transformations with bounded rank},
J.\ London Math.\ Soc.\ {\bf 37} 10-16 (1962).

\bibitem[FHM]{FeldmanHahnMoore}
Feldman~J., Hahn~P., Moore~C.~C., {\em Orbit structure and countable
sections for actions of continuous groups}, Advances in Math.\
{\bf 28}, 186-230 (1978).

\bibitem[FM1]{FeldmanMoore1}
Feldman~J., Moore~C.~C., {\em Ergodic equivalence relations, cohomology
and Von Neumann Algebras I}, Trans.\ Am.\ Math.\ Soc.\ {\bf 234} 
No.~2, 289-324 (1977).

\bibitem[FM2]{FeldmanMoore2}
Feldman~J., Moore~C.~C., {\em Ergodic equivalence relations, cohomology
and Von Neumann Algebras II}, Trans.\ Am.\ Math.\ Soc.\ {\bf 234} 
No.~2, 325-359 (1977).

\bibitem[GD-A]{GlicksbergDeLeeuwA}
Glicksberg~I., DeLeeuw~K.,
{\em Applications of Almost Periodic Compactifications},
Acta Mathematica {\bf 105} 63-97 (1961).

\NOT{
\bibitem[GD-S]{GlicksbergDeLeeuwS}
Glicksberg~I., DeLeeuw~K.,
{\em Almost Periodic Functions on Semigroups},
Acta Mathematica {\bf 105} 99-140 (1961).
}

\bibitem[GD-D]{GlicksbergDeLeeuwD}
Glicksberg~I., DeLeeuw~K.,
{\em The Decomposition of Certain Group Representations},
J.~d'Analyse Math.\ {\bf 15} 135-192 (1965).

\bibitem[Gr]{Greenleaf}
Greenleaf~F.~P., {\em Invariant Means on Topological Groups and Their
Applications}.
Van Nostrand, 1969. 

\bibitem[Ha-M]{HalmosMEA}
Halmos~P.~R., {\em Measure Theory}. Van Nostrand, 1950.

\bibitem[Ha-B]{HalmosBOO}
Halmos~P.~R., {\em Lectures on Boolean Algebras}. Van Nostrand, 1963.

\bibitem[He]{Helmberg}
Helmberg~G.,
{\em \"Uber mittlere R\"uckkehrzeit unter einer masstreuen
Str\"omung},
Z.\ Wahrsch.\ verw.\ Geb.\ {\bf 13}, 165-179 (1969).

\bibitem[Hi]{Hicks}
Hicks~N.~J., {\em Notes on Differential Geometry}. 
Van Nostrand, 1965.

\bibitem[J]{Jacobs}
Jacobs~K., {\em Neuere Methoden und Ergebnisse der Ergodentheorie}.
Springer-Verlag, 1960.

\bibitem[Kac]{Kac}
Kac~M., {\em On the notion of recurrence in discrete stochastic processes},
Bull.\ Am.\ Math.\ Soc.\ {\bf 53}, 1002-1010 (1947).

\bibitem[Kas]{Kastelyn}
Kastelyn~P.~W., {\em Variations on a theme by Marc Kac}, Journal of Statistical
Physics, {\bf 46} Nos.~5/6, 811-827 (1987).

\bibitem[Ke]{Kechris}
Kechris~A.~S., {\em Countable sections for locally compact group actions},
Ergod.\ Th.\ \& Dynam.\ Sys.\ {\bf 12}, 283-295 (1992).

\bibitem[Ke-D]{KechrisD}
Kechris~A.~S., {\em Classical Descriptive Set Theory}. Springer-Verlag, 1995.

\bibitem[Ku]{Kuratowski}
Kuratowski~K., {\em Topology}. Vol.~1. Academic Press, New York and London,
and Pa\'{n}stwowe Wydawnictwo Naukowe, Warszawa, 1966.

\bibitem[Li]{Lindvall}
Lindvall~T., {Lectures on the Coupling Method}. John Wiley \& Sons, 1992.

\bibitem[Loo]{Loomis}
Loomis~L.~H., {\em An Introduction to Abstract Harmonic Analysis}.
Van Nostrand, 1953.

\bibitem[Lu]{Lusin}
Lusin~N., {\em Le\c{c}ons sur les Ensembles Analytiques et Leurs
Applications}. Gauthier-Villars, Paris, 1930.

\bibitem[Ma]{Mackey}
Mackey~G.~W., {\em Point realizations of transformation groups}, Ill.\ J.\
Math.\ {\bf 6}, 327-335 (1962).

\bibitem[Me]{Mecke}
Mecke~J., {\em Station\"are Zuf\"alige Masse auf Lokalkompakten Abelschen
Gruppen}, Z.\ Wahrsch.\ verw.\ Geb.\ {\bf 9}, 36-58 (1967).

\bibitem[Mo]{Moschovakis}
Moschovakis~Y.~N., {\em Descriptive Set Theory}. North-Holland Publishing
Company, 1980.

\NOT{
\bibitem[MS]{MuhlySolel}
Muhly~P.~S., Solel~B., {\em Subalgebras of Groupoid $C^*$-algebras},
J.\ reine angew.\ Math.\ {\bf 402}, 41-75 (1989).
}

\bibitem[Na]{Nadkarni}
Nadkarni~M.~G., {\em Basic Ergodic Theory}. Hindustan Book Agency, 
Delhi, 1995.

\bibitem[Ne]{Neveu}
Neveu~J., {\em Sur la structure des processus ponctuels stationnaires},
C.\ R.\ Acad.\ Sc.\ Paris, {\bf 267} 561-564 (1968).

\bibitem[Ne-P]{NeveuP}
Neveu~J., {\em Processus Ponctuels}. In: Lecture Notes in Math., No.\ 598,
Springer-Verlag, 1977.

\bibitem[Ow]{Owen}
Owen~G., {\em Game Theory}. 2nd ed., Academic Press, 1982.

\bibitem[Pa]{Paterson}
Paterson~A.~L., {\em Amenability}. Mathematical Surveys and Monographs No.~29,
Am.\ Math.\ Soc., 1988.

\bibitem[Pe]{Petersen}
Petersen~K., {\em Ergodic Theory}. Cambridge University Press, 1989.

\bibitem[Ph]{Phelps}
Phelps~R.~R., {\em Lectures on Choquet's Theorem}.
Van Nostrand, 1966.

\NOT{
\bibitem[Phi]{Phillips}
Phillips~W.~J., {\em Flow under a function and discrete decomposition of
properly infinite $W^*$-algebras}, Pacific J.\ of Math.\ {\bf 114} No.~1,
221-234 (1984).
}

\bibitem[Ra]{Ramsay}
Ramsay~A., {\em Local product structure for group actions}, Ergodic Th.\
\& Dynam.\ Sys.\ {\bf 11}, 209-217 (1991).

\bibitem[Sch]{Schwartz}
Schwartz~J.~T., {\em Differential Geometry and Topology}. Gordon and
Breach, 1968.

\bibitem[Ta]{TarskiCA}
Tarski~A., {\em Cardinal Algebras}. Oxford University Press, 1969.

\bibitem[Tr-E]{TroallicE}
Troallic~J.~P., {\em Espaces fonctionelles et th\'eor\`emes de I.~Namioka},
Bull.\ Soc.\ Math.\ France, {\bf 107}, 127-137 (1979).

\bibitem[Tr-S]{TroallicS}
Troallic~J.~P., {\em Semigroupes semitopologiques et
presque-p\'eriodicit\'e}, in: Recent Developments in the Algebraic,
Analytical and Topological Theory of Semigroups, Lecture Notes No.~998,
Springer-Verlag, 1981.

\bibitem[Va]{Varadarajan}
Varadarajan~V.~S., {\em Groups of automorphisms of Borel spaces}, Trans.\
Am.\ Math.\ Soc.\ {\bf 109}, 191-220 (1963).

\bibitem[vN]{vN}
Von~Neumann~J., {\em Zur Theorie der Gesellschaftsspiele}, Mathematische
Annalen, {\bf 100}, 295-320 (1928).

\bibitem[vNM]{vNM}
Von~Neumann~J., Morgenstern~O., {\em Theory of Games and Economic Behavior}.
3rd ed., Princeton University Press, 1953.

\bibitem[Wg]{Wagh}
Wagh~V.~M., {\em A descriptive version of Ambose' representation theorem
for flows}, Proc.\ Indian Acad.\ Sci.\ (Math.\ Sci.) {\bf 98} No. 2-3,
101-108 (1988).

\bibitem[Wa]{WagonBT}
Wagon~S., {\em The Banach-Tarski Paradox}. Cambridge University Press,
1986.

\bibitem[We0]{WehrungHB}
Wehrung~F., {\em Th\'{e}or\`{e}me de Hahn-Banach et paradoxes continus et
discrets}, C.~R. Acad. Sci. Paris {\bf 310} I, 303-306 (1990).

\bibitem[We1]{WehrungInj1}
Wehrung~F., {\em Injective positively ordered monoids I}, J. Pure Appl. Alg.
{\bf 83}, 43-82 (1992).

\bibitem[We2]{WehrungInj2}
Wehrung~F., {\em Injective positively ordered monoids II}, J. Pure Appl. Alg.
{\bf 83}, 83-100 (1992).

\bibitem[Wl]{Weil}
Weil~A., {\em Basic Number Theory}. Springer-Verlag, 1967.

\bibitem[Wn]{Weinstein}
Weinstein~A., {\em Groupoids: unifying internal and external symmetry},
Notices of the AMS, {\bf 43} No. 7, 744-752, July 1996.
\end{thebibliography}
\end{document}